\documentclass[a4paper]{amsart}

\usepackage{mathrsfs,amssymb}
\usepackage{multicol}\multicolsep=0pt
\usepackage{pstricks,pst-node}

\usepackage[sort]{cite}

\def\crulefill{\leavevmode\leaders\hrule height 1pt\hfill\kern 0pt}
\long\def\QUERY#1{%
\leavevmode\newline%
\noindent$\star\star\star$\thinspace\textsf{Comment/Query}\crulefill\newline%
   \space #1\newline\hbox to 120mm{\crulefill}$\star\star\star$\newline
}

\numberwithin{equation}{section}
\theoremstyle{definition}
\newtheorem{Defn}[equation]{Definition}
\newtheorem{Example}[equation]{Example}
\newtheorem{Remark}[equation]{Remark}
\theoremstyle{plain}
\newtheorem{Prop}[equation]{Proposition}
\newtheorem{Theorem}[equation]{Theorem}
\newtheorem{Assumption}[equation]{Assumption}
\newtheorem{Lemma}[equation]{Lemma}
\newtheorem{Cor}[equation]{Corollary}
\newtheorem{Point}[equation]{}

\newtheorem{THEOREM}{Theorem}

\def\Case#1{\medskip\noindent\textbf{Case #1}:\leavevmode\newline}

\makeatletter
\def\enumerate{\begingroup\ifnum\@enumdepth>3\@toodeep\else
      \advance\@enumdepth\@ne
      \edef\@enumctr{enum\romannumeral\the\@enumdepth}%
      \topsep\z@\parskip\z@
      \list{\csname label\@enumctr\endcsname}
        {\@nmbrlisttrue\let\@listctr\@enumctr
         \parsep\z@\itemsep\z@\topsep\z@
         \setcounter{\@enumctr}{0}
         \def\makelabel##1{\hss\llap{\rm ##1}}
       }\fi}

\makeatother

\let\bar=\overline
\let\epsilon=\varepsilon
\def\({\big(}
\def\){\big)}

\def\C{\mathbb C}
\def\N{\mathbb N}
\def\R{\mathbb R}
\def\Q{\mathbb Q}
\def\Z{\mathbb Z}
\def\ZC{\mathcal Z}

\def\0{\underline{0}}
\def\bu{\mathbf u}
\def\B{\mathscr B}
\def\BDiag{\mathcal B}
\def\Bcal{\mathfrak B}
\def\Dcal{\mathcal D}
\def\Ef{{\mathcal E}_f}
\def\G{\mathcal G}
\def\H{\mathscr H}
\let\proj=\varepsilon
\def\Nrf{\mathbb N_r^{(f)}}
\def\Sym{\mathfrak S}
\def\W{\mathscr W}
\def\Waff{\W^{\text{aff}}}
\def\Wlam{\W_{r,n}^{\gdom\lambda}}
\def\Wlambda{\W_{r,n}^{\gedom\lambda}}

\DeclareMathOperator{\Char}{char}

\def\simk{\overset{k}\sim}
\def\simkk{\overset{k+1}\sim}
\def\siml{\overset{l}\sim}
\DeclareMathOperator{\End}{End}
\DeclareMathOperator{\GL}{GL}
\DeclareMathOperator{\Rad}{Rad}
\DeclareMathOperator*{\Res}{Res}

\def\half{\frac12}

\let\gdom\rhd
\let\gedom\unrhd
\def\a{\mathfraxk a}

\def\floor#1{\lfloor\tfrac#1\rfloor}
\def\UPD{\mathscr{T}^{ud}}
\def\Std{\mathscr{T}^{std}}

\def\s{\mathfrak s}
\def\ts{\tilde\s}
\def\t{\mathfrak t}
\def\u{\mathfrak u}
\def\v{\mathfrak v}

{\catcode`\|=\active
  \gdef\set#1{\mathinner{\lbrace\,{\mathcode`\|"8000%
                                   \let|\midvert #1}\,\rbrace}}
  \gdef\seT#1{\mathinner{\Big\lbrace\,{\mathcode`\|"8000%
                                   \let|\midverT #1}\,\Big\rbrace}}
}
\def\midvert{\egroup\mid\bgroup}
\def\midverT{\egroup\,\Big|\,\bgroup}

\def\Set[#1]#2|#3|{\Big\{\ #2\ \Big| \
           \vcenter{\hsize #1mm\centering #3}\Big\}}

\def\map#1#2{\,{:}\,#1\!\longrightarrow\!#2}

\def\BrauerDiagram#1#2{\vcenter{\hbox{%
        \begin{psmatrix}[colsep=2mm,rowsep=1mm,emnode=dot,dotscale=1.3,nodesep=0pt]%
              #1&\\#1&#2%
        \end{psmatrix}%
}}}

\title{Cyclotomic Nazarov--Wenzl algebras}

\author{Susumu Ariki}
\address{Research Institute for Mathematical Sciences, %
         Kyoto University, Kyoto 606-8502. Japan.}
\email{ariki@kurims.kyoto-u.ac.jp}

\author{Andrew Mathas}
\address{School of Mathematics and Statistics F07, %
         University of Sydney, NSW 2006, Australia}
\email{a.mathas@maths.usyd.edu.au}

\author{Hebing Rui}
\address{Department of Mathematics, East China Normal University, %
         200062 Shanghai, P.R. China.}
\email{hbrui@math.ecnu.edu.cn}

\thanks{\kern-\parindent The first and second named authors enjoyed
the hospitality of the Bernoulli Center at the EPFL Lausanne when they
visited the program ``Group representation theory'' in June 2005. The
paper was completed during their visit and the authors acknowledge
partial support from the program. The second and third authors were
supported, in part by Fellowships from the Japanese Society for the
Promotion of Science. All three authors were supported by ARC grant
``Quantized representation theory'' during a visit to Sydney in
January 2005. The third author was also supported by the Natural Science Foundation in
China.}

\begin{document}

\begin{abstract}
Nazarov~\cite{Nazarov:brauer} introduced an infinite dimensional
algebra, which he called the \textit{affine Wenzl algebra}, in his
study of the Brauer algebras. In this paper we study certain
``cyclotomic quotients'' of these algebras. We construct the
irreducible representations of these algebras in the generic case and
use this to show that these algebras are free of rank $r^n(2n-1)!!$
(when $\Omega$ is $\bu$--admissible). We next show that these algebras
are cellular and give a labelling for the simple modules of the
cyclotomic Nazarov--Wenzl algebras over an arbitrary field. In
particular, this gives a construction of all of the finite dimensional
irreducible modules of the affine Wenzl algebra.
\end{abstract}

\sloppy \maketitle

\centerline{\textit{On the occasion of Professor George Lusztig's
                      60$^{\text{th}}$ birthday}}

\section{Introduction}

The \textsf{Brauer algebras} were introduced by Richard
Brauer~\cite{BrauerAlg} in his study of representations of the
symplectic and orthogonal groups. In introducing these algebras Brauer
was motivated by Schur's theory (see \cite{Green}), which links the
representation theory of the symmetric group $\Sym_n$ and the general
linear group $\GL(V)$ via their commuting actions on ``tensor space''
$V^{\otimes n}$, where $\Sym_n$ acts by place permutations.
Analogously, the Brauer algebras are the centralizers of the image of
a symplectic or orthogonal group in $\End(V^{\otimes n})$, where $V$
is the defining representation of the group.

The Brauer algebras have now been studied by many authors and they
have applications ranging from Lie theory, to combinatorics and knot
theory; see, for example,
\cite{BirmanWenzl,Brown:brauer,DoranHanlonWales,Enyang,FiGr,HalRam:basic,HW2,HW1,J:Brauer,LeducRam,Martin:partition,Nazarov:brauer,Ram:chBrauer,Rui:ssbrauer,Terada:brauer,Wenzl:ssbrauer,Weyl,CCXi:BMW}.
In this paper we are interested not so much in the Brauer algebra
itself but in affine and cyclotomic analogues of it. Our starting
point is a (special case of) Nazarov's~\cite{Nazarov:brauer} affine
Wenzl algebra $\Waff_n(\Omega)$, an algebra which could legitimately
be called the degenerate affine BMW algebra. Nazarov introduced the
affine Wenzl algebra when studying the action of ``Jucys--Murphy
operators'' on the irreducible representations of the Brauer algebras.
Nazarov's idea was that the affine Wenzl algebra should play a similar
role in the representation theory of the Brauer algebras to that
played by the affine degenerate Hecke algebra of type~A in the
representation theory of the symmetric group.

Let $R$ be a commutative ring. The representation theory of the affine
Wenzl algebras $\Waff_n(\Omega)$, where $\Omega=\set{\omega_a\in
R|a\ge0}$, has not yet been studied.  Motivated by the theory of the
affine Hecke algebras and the cyclotomic Hecke algebras of type
$G(r,1,n)$ \cite{Ariki:can,DJM:cyc,Klesh:book} we introduce a
``cyclotomic'' quotient $\W_{r,n}(\mathbf
u)=\Waff_n(\Omega)/\langle\prod_{i=1}^r(X_1-u_i)\rangle$ of
$\Waff_n(\Omega)$, which depends on an $r$--tuple of parameters
$\bu=(u_1,\dots,u_r)\in R^r$. We call $\W_{r,n}(\bu)$ a
\textit{cyclotomic Nazarov--Wenzl algebra}. This paper develops
the representation theory of the algebras $\W_{r,n}(\bu)$.

The first question that we are faced with is whether the cyclotomic
Nazarov--Wenzl algebra $\W_{r,n}(\bu)$ is always free as an $R$--module. The
Brauer algebra $\B_n$ is free of rank
$(2n-1)!!=(2n-1)\cdot(2n-3)\dots\cdot3\cdot 1$. We expect that the
cyclotomic Nazarov--Wenzl algebra $\W_{r,n}(\bu)$ should be free of rank
$r^n(2n-1)!!$. In section~3, a detailed study of the representation
theory of $\W_{r,2}(\bu)$ shows that, in the semisimple case,
$\W_{r,2}(\bu)$ has rank $r^n(2n-1)!!|_{n=2}$ if and only if $\Omega$
is $\bu$--admissible (Definition~\ref{u-admissible}). This constraint
on $\Omega$ involves Schur's $q$--functions. Our first main result is
the following.

\begin{THEOREM}\label{RANK}
Let $R$ be a commutative ring in which $2$ is invertible.  Suppose
that $\bu\in R^r$ and that $\Omega$ is $\bu$--admissible.  Then the
cyclotomic Nazarov--Wenzl algebra $\W_{r,n}(\bu)$ is free as an
$R$--module of rank $r^n(2n-1)!!$.
\end{THEOREM}

The proof of this result occupies a large part of this paper. The idea
behind the proof comes from \cite{AK}: for ``generic'' $R$ we
explicitly construct a class of irreducible representations of
$\W_{r,n}(\bu)$ and use them to show that the dimension of
$\W_{r,n}(\bu)/\Rad\W_{r,n}(\bu)$ is at least $r^n(2n-1)!!$. It is
reasonably easy to produce a set of $r^n(2n-1)!!$ elements which span
$\W_{r,n}(\bu)$, so this is enough to prove Theorem~\ref{RANK}.  We
construct these irreducible representations by giving ``seminormal
forms'' for them (Theorem~\ref{seminormal}); that is, we give explicit
matrix representations for the actions of the generators of
$\W_{r,n}(\bu)$.  The main difficulty in this argument is in showing
that these matrices respect the relations of $\W_{r,n}(\bu)$, we do
this using generating functions introduced by
Nazarov~\cite{Nazarov:brauer}. There is an additional subtlety in that
we have to work over the real numbers in order to make a consistent
choice of square roots in the representing matrices.

The next main result of the paper shows that $\W_{r,n}(\bu)$ is a
cellular algebra in the sense of Graham and Lehrer~\cite{GL}. This
gives a lot of information about the representations of the cyclotomic
Nazarov--Wenzl algebras. For example, cellularity implies that the
decomposition matrix of $\W_{r,n}(\bu)$ is unitriangular.

\begin{THEOREM}\label{CELLULAR}
Suppose that $2$ is invertible in $R$ and that $\Omega$ is
$\bu$--admissible. Then the cyclotomic Nazarov--Wenzl algebra
$\W_{r,n}(\bu)$ is a cellular algebra.
\end{THEOREM}

We prove Theorem~\ref{CELLULAR} by constructing a cellular basis for
$\W_{r,n}(\bu)$. We recall the definition of a cellular basis in
section~\ref{cellular}; however, for the impatient experts we mention
that the cell modules of $\W_{r,n}(\bu)$ are indexed by  ordered pairs
$(f,\lambda)$, where $0\le f\le\floor{n2}$ and $\lambda$ is a
multipartition of $n-2f$, where $0\le f\le\floor{n2}$, and the bases
of the cell modules are indexed by certain ordered triples which are
in bijection with the $n$--updown $\lambda$--tableaux.

Finally we consider the irreducible $\W_{r,n}(\bu)$--modules over a
field~$R$. The cell modules of $\W_{r,n}(\bu)$ have certain quotients
$D^{(f,\lambda)}$, where $0\le f\le\floor{n2}$ and $\lambda$ is a
multipartition of $n-2f$, which the theory of cellular algebras says
are either zero or absolutely irreducible. Now, the cyclotomic
Nazarov--Wenzl algebra $\W_{r,n}(\bu)$ is filtered by two sided ideals
with the degenerate Hecke algebras $\H_{r,n-2f}(\bu)$ of type
$G(r,1,n-2f)$ appearing as the successive quotients for $0\le
f\le\floor{n2}$. In section~6 we show that the algebras
$\H_{r,m}(\bu)$ are also cellular (in fact, this is the key to proving
Theorem~\ref{CELLULAR}); as a consequence, the irreducible
$\H_{r,m}(\bu)$--modules are the non--zero modules $D^\lambda$,
where~$\lambda$ is a multipartition of~$m$.  

Using the results of the last paragraph we can construct all of the
irreducible representations of the cyclotomic Nazarov--Wenzl algebras
when~$\Omega$ is admissible. This enables us to construct all finite
dimensional representations of the affine Wenzl algebras over an
algebraically closed field when~$\Omega$ is admissible
(Theorem~\ref{finite irreds}). In the special case when $\omega_0\ne0$
we also have the following classification of the irreducible
$\W_{r,n}(\bu)$--modules.

\begin{THEOREM}\label{simple classification}
Suppose that $R$ is a field in which $2$ is invertible, that $\Omega$ is
$\bu$--admissible and that $\omega_0\ne0$. Then
$\set{D^{(f,\lambda)}|0\le f\le\floor{n2}, \lambda\vdash n-2f
       \text{ and } D^\lambda\ne0}$
is a complete set of pairwise non--isomorphic irreducible
$\W_{r,n}(\bu)$--modules.
\end{THEOREM}

As an application of Theorem~\ref{simple classification} we give
necessary and sufficient conditions for~$\W_{r,n}(\bu)$ to be
quasi--hereditary when $R$ is a field and $\omega_0\ne0$.

We note that Orellana and Ram \cite{OrelRam}, building on
\cite{LeducRam},  gave explicit formulae for the seminormal
representations of the affine and cyclotomic BMW algebras
\cite[Theorem~6.20]{OrelRam}.  As the (degenerate) affine Wenzl algebra
is the degeneration of the affine BMW algebra, it is natural to expect
that we should be able to derive the seminormal representations of
$\Waff_n(\Omega)$ from the results of \cite{OrelRam}. Unfortunately,
this is not possible because Orellana and Ram construct the seminormal
representations only for a very restrictive class of cyclotomic BMW
algebras corresponding to certain specializations of the parameters
(see \cite[Theorem~6.17(c-d)]{OrelRam}). These parameter choices are
sufficient for the purposes of \cite{OrelRam}, however, it is not
clear that ``enough'' of these parameter choices are
$\bu$--admissible so we are unable to exploit~\cite{OrelRam}. We also
remark we had to work quite hard to ensure that we had made a
consistent choice of square roots in our seminormal representations
(cf.~Assumption~\ref{be}), and that it is not clear to us that
\cite{OrelRam} have made a coherent choice of roots in their seminormal
representations.

Another possible approach to the construction of the seminormal
representations in this paper is via Jones' ``basic
construction''~\cite{Wenzl:ssbrauer}.  Wenzl constructed the
semisimple irreducible representations of the Brauer algebras this
way. The final step of his argument used the non--degeneracy of the
Markov trace to show that all partitions of $n$ are ``permissible''
(see \cite[Theorem~3.4]{Wenzl:ssbrauer}). We were unable to extend
Wenzl's arguments to the cyclotomic Nazarov--Wenzl algebras because we
were unable to prove that we have an analogous non--degenerate trace
form.  Note also that, a priori, it is not clear that
$\W_{r,n-1}(\bu)$ is a subalgebra of~$\W_{r,n}(\bu)$.

Finally, we remark that other variants of signed and cyclotomic Brauer
algebras, $G$-Brauer algebras and cyclotomic BMW
have been studied previously in the papers
\cite{CGW:BMW,GoodmanHauschild,HO:cycBMW,ParvathiKamaraj:signedBrauer,ParvathiKamaraj:matrixUnits,ParvathiSavithri,RuiYu}.

\def\NI[#1|#2|#3]{%
  \rlap{$#1$}\hspace*{12mm}%
  \hbox to 46mm{#2\dotfill}%
  \hspace*{4mm}\llap{\pageref{#3}}%
}
\medskip
{\small
\multicolsep=1mm
\columnsep=4mm
\begin{center}\parindent=0pt%
\begin{multicols}{2}[\textsc{\normalsize Notational index}]\noindent%
\NI[a_\t(k)|Seminormal matrix entry|ab]
\NI[b_\t(k)|Seminormal matrix entry|ab]
\NI[\B_n(\omega)|Brauer algebra|Brauer presentation]
\NI[\BDiag(n)|Set of Brauer diagrams|Brauer diagrams]
\NI[B_\gamma|Element of $\Waff_n(\lambda)$ or $\W_{r,n}$|B gamma]
\NI[c_\t(k)|Content of $k$ in $\t$|content]
\NI[C^{(f,\lambda)}_{(*)(*)}|Cellular basis element of $\W_{r,n}$|W cell basis]
\NI[\deg|Degree function on $\W_{r,n}$|degree]
\NI[E_i|Generator of $\Waff_n(\Omega)$|Waff relations]
\NI[E^f|$=E_{n-1}E_{n-3}\dots E_{n-2f+1}$|Ef]
\NI[\Ef|$=\W_{r,n-2f}E_1\W_{r,n-2f}$|epsilon f]
\NI[e_{\t\u}(k)|Seminormal matrix entry|e]
\NI[f^{(n,\lambda)}|$=\#\UPD_n(\lambda)$|updown]
\NI[\H_{r,n}(\bu)|Degenerate Hecke algebra|degenerate H]
\NI[m_{\s\t}|Cellular basis element of $\H_{r,n}$|H cellular]
\NI[M_{\s\t}|Element of $\W_{r,n}(\bu)$|M_st]
\NI[S_i|Generator of $\Waff_n(\Omega)$|Waff relations]
\NI[s_{\t\u}(k)|Seminormal matrix entry|x-equal]
\NI[T_i|Generator of $\H_{r,n}(\bu)$|degenerate H]
\NI[\Std(\lambda)|Standard $\lambda$-tableaux|tableaux]
\NI[\UPD_n(\lambda)|Updown $\lambda$-tableaux|tableaux]
\NI[\bu|Parameters for $\W_{r,n}(\bu)$|cyclo NZ]
\NI[\W_{r,n}(\bu)|Cyclotomic Nazarov--Wenzl alg.|cyclo NZ]
\NI[\W_{r,n}^f|$=\W_{r,n}E^f\W_{r,n}$|Ef]
\NI[\Wlambda|Two--sided ideal in $\W_{r,n}^f$|Wlambda]
\NI[\Waff_n(\Omega)|Affine Wenzl algebra|Waff relations]
\NI[\tilde W_1(y)|$=\sum_{a\ge0}\omega_ay^{-a}$|W_1 remark]
\NI[W_k(y)|$=\frac12-y+(y-\frac{(-1)^r}2)\prod_{\alpha}%
    \frac{y+c(\alpha)}{y-c(\alpha)}$|rational-W]
\NI[X_j|Generator of $\Waff_n(\Omega)$|Waff relations]
\NI[X^\alpha|A monomial $X_1^{\alpha_1}\dots X_n^{\alpha_n}$|X alpha]
\NI[Y_j|Generator of $\H_{r,n}(\bu)$|degenerate H]
\NI[\Delta(\lambda)|Seminormal representation|seminormal rep]
\NI[\sigma_f|A map $\H_{r,n-2f}\to\W^f_{r,n}/\W_{r,n}^{f+1}$|sigma f]
\NI[\omega_a,\Omega|Parameters for $\Waff_n(\Omega)$|Waff relations]
\NI[\omega_k^{(a)}|Central element|tilde W]
\NI[\simk|Equivalence relation on $\UPD_n(\lambda)$|tilde]
\NI[\t\ominus\u|$\t\setminus\u$ or $\u\setminus\t$|ab]
\end{multicols}
\end{center}
}

\section{Affine and cyclotomic Nazarov--Wenzl algebras}

In~\cite{Nazarov:brauer}, Nazarov introduced an affine analogue of the
Brauer algebra which he called the (degenerate) \textsf{affine Wenzl}
algebra. The main objects of interest in this paper are certain
``cyclotomic'' quotients of Nazarov's algebra. In this section we
define these algebras and prove some elementary results about them.

Fix a positive integer $n$ and a commutative ring $R$ with
multiplicative identity~$1_R$. \textit{Throughout this paper we will
assume that $2$ is invertible in~$R$}.


\begin{Defn}[\protect{Nazarov~\cite[\S4]{Nazarov:brauer}}]
\label{Waff relations}
Fix $\Omega=\set{\omega_a|a\ge0}\subseteq R$.  The 
(degenerate) \textsf{affine Wenzl algebra} $\Waff_n=\Waff_n(\Omega)$ is 
the unital associative $R$--algebra with generators
$\set{S_i,E_i, X_j|1\le i<n \text{ and }1\le j\le n } $
and relations
\begin{multicols}{2}
\begin{enumerate}
    \item (Involutions)\newline
$S_i^2=1$, for $1\le i<n$.
    \item (Affine braid relations)
\begin{enumerate}
\item $S_iS_j=S_jS_i$ if $|i-j|>1$,
\item $S_iS_{i+1}S_i=S_{i+1}S_iS_{i+1}$,\newline for $1\le i<n-1$,
\item $S_iX_j=X_jS_i$ if $j\ne i,i+1$.
\end{enumerate}
    \item (Idempotent relations)\newline
$E_i^2=\omega_0E_i$, for $1\le i<n$.
    \item (Commutation relations)
\begin{enumerate}
\item $S_iE_j=E_jS_i$, if $|i-j|>1$,
\item $E_iE_j=E_jE_i$, if $|i-j|>1$,
\item $E_iX_j=X_jE_i$,\newline if $j\ne i,i+1$,
\item $X_iX_j=X_jX_i$,\newline for $1\le i,j\le n$.
\end{enumerate}
    \item (Skein relations)\newline
        $S_iX_i-X_{i+1}S_i=E_i-1$ and\newline
        $X_iS_i-S_iX_{i+1}=E_i-1$,\newline for $1\le i<n$.
    \item (Unwrapping relations)\newline
        $E_1X_1^aE_1=\omega_aE_1$, for $a>0$.
    \item (Tangle relations)
\begin{enumerate}
\item $E_iS_i=E_i=S_iE_i$,\newline for $1\le i\le n-1$,
\item $S_iE_{i+1}E_i=S_{i+1}E_i$,\newline for $1\le i\le n-2$,
\item $E_{i+1}E_iS_{i+1}=E_{i+1}S_i$,\newline for $1\le i\le n-2$.
\end{enumerate}
    \item (Untwisting relations)\newline
        $E_{i+1}E_iE_{i+1}=E_{i+1}$ and\newline
        $E_iE_{i+1}E_i=E_i$, for $1\le i\le n-2$.
    \item (Anti--symmetry relations)\newline
        $E_i(X_i+X_{i+1})=0$ and\newline
        $(X_i+X_{i+1})E_i=0$, for $1\le i<n$.
\end{enumerate}
\end{multicols}
\end{Defn}

Our definition of $\Waff_n$ differs from Nazarov's in two respects.
First, Nazarov considers only the special case when $R=\C$; however,
as we will indicate, most of the arguments that we need
from~\cite{Nazarov:brauer} go through without change when $R$ is an
arbitrary ring. More significantly, Nazarov considers a more general
algebra which is generated by the elements $\set{S_i,E_i, X_j,
\hat\omega_a|1\le i<n, 1\le j\le n \text{ and } a\ge0}$ such that the
$\hat\omega_a$ are central and the remaining generators satisfy the
relations above. For our purposes it is more natural to define the
elements $\omega_a$ to be elements of~$R$ because without this
assumption the cyclotomic quotients of~$\Waff_n$ would not be finite
dimensional.

Note that $E_iE_{i+1}S_i=E_iE_{i+1}E_iS_{i+1}=E_iS_{i+1}$
and $S_{i+1}E_iE_{i+1}=S_iE_{i+1}E_iE_{i+1}=S_iE_{i+1}$. Thus
a quick inspection of the defining relations shows that $\Waff_n$ has
the following useful involution.

\begin{Point}\label{* involution}
There is a unique $R$--linear anti--isomorphism $*\map{\Waff_n}\Waff_n$ such that
$$ S_i^*=S_i,\qquad E_i^*=E_i\quad\text{and}\quad X_j^*=X_j,$$
for all $1\le i<n$ and all $1\le j\le n$. Moreover, $*$ is an
involution.
\end{Point}

Using the defining relations it is not hard to see that $\Waff_n$ is
generated by the elements $S_1,\dots,S_{n-1},E_1,X_1$. There
is no real advantage, however, to using this smaller set of generators
as the corresponding relations are more complicated.

\begin{Lemma}[\protect{cf. \cite[(2.6)]{Nazarov:brauer}}]\label{SX^a}
Suppose that $1\le i<n$ and that $a\ge1$. Then
$$S_iX_i^a=X_{i+1}^aS_i+\sum_{b=1}^aX_{i+1}^{b-1}(E_i-1)X_i^{a-b}.$$
\end{Lemma}

\begin{proof} We argue by induction on $a$. When $a=1$ this is relation
    \ref{Waff relations}(e) . If $a\ge1$ then, by induction, we have
\begin{align*}
S_iX_i^{a+1}&=S_iX_i^aX_i
     =\Big\{X_{i+1}^aS_i+\sum_{b=1}^aX_{i+1}^{b-1}(E_i-1)X_i^{a-b}\Big\}X_i\\
    &=X_{i+1}^aS_iX_i+\sum_{b=1}^aX_{i+1}^{b-1}(E_i-1)X_i^{a+1-b}.\\
    \intertext{Now, by the skein relation \ref{Waff relations}(e),
    $S_iX_i=X_{i+1}S_i+E_i-1$, so}
S_iX_i^{a+1} &=X_{i+1}^a\big(X_{i+1}S_i+E_i-1\big)
                 +\sum_{b=1}^aX_{i+1}^{b-1}(E_i-1)X_i^{a+1-b}\\
             &=X_{i+1}^{a+1}S_i
                 +\sum_{b=1}^{a+1}X_{i+1}^{b-1}(E_i-1)X_i^{a+1-b},\\
\end{align*}
as required.
\end{proof}

\begin{Cor}\label{constraint}
Suppose that $a\ge0$. Then
$$
\omega_{2a+1}E_1=\frac12\Big\{-\omega_{2a}
            +\sum_{b=1}^{2a+1}(-1)^{b-1}\omega_{b-1}\omega_{2a+1-b}\Big\}E_1.$$
\end{Cor}

\begin{proof}
Take $i=1$ and multiply the equation in Lemma~\ref{SX^a} on the
left and right by $E_1$. Since $S_1E_1=E_1=E_1S_1$, this gives
\begin{align*}
E_1X_1^aE_1&= E_1X_2^aE_1+\sum_{b=1}^aE_1X_2^{b-1}(E_1-1)X_1^{a-b}E_1.\\
\intertext{Since $E_1X_1^cE_1=\omega_cE_1$, $E_1(X_1+X_2)=0$ and
$X_1X_2=X_2X_1$ we can
rewrite this equation as}
\omega_aE_1&=(-1)^a\omega_a E_1
         +\sum_{b=1}^a (-1)^{b-1}E_1X_1^{b-1}(E_1-1)X_1^{a-b}E_1\\
   &=(-1)^a\omega_a E_1+\sum_{b=1}^a (-1)^{b-1}
        \big(E_1X_1^{b-1}E_1X_1^{a-b}E_1-E_1X_1^{a-1}E_1\big).\\
   &=(-1)^a\omega_a E_1
         +\sum_{b=1}^a (-1)^{b-1}(\omega_{b-1}\omega_{a-b}-\omega_{a-1})E_1.\\
   &=(-1)^a\omega_a E_1
         +\sum_{b=1}^a (-1)^{b-1}\omega_{b-1}\omega_{a-b}E_1
         +\sum_{b=1}^a(-1)^b\omega_{a-1}E_1.\\
\end{align*}
Setting $a=2a'+1$ proves the Corollary.
\end{proof}

If we assume that $E_1\ne0$ in $\Waff_n$ and that $\Waff_n$ is torsion
free then this result says that the $\omega_a$, for $a$ odd, are
determined by the $\omega_b$, for $b$ even.

\begin{Remark} If $a>0$ then the proof of the Corollary also gives the identity
$$0=\Big\{\sum_{b=1}^{2a}(-1)^{b-1}\omega_{b-1}\omega_{2a-b}\Big\}E_1.$$
However, this relation holds automatically because
\begin{align*}
    \sum_{b=1}^{2a}(-1)^{b-1}\omega_{b-1}\omega_{2a-b}
        &=\sum_{b=1}^a(-1)^{b-1}\omega_{b-1}\omega_{2a-b}
        +\sum_{b=a+1}^{2a}(-1)^{b-1}\omega_{b-1}\omega_{2a-b} \\
        &=\sum_{b=1}^a(-1)^{b-1}\omega_{b-1}\omega_{2a-b}
        +\sum_{b'=1}^{a}(-1)^{2a-b'}\omega_{2a-b'}\omega_{b'-1} \\
        &=0.
\end{align*}
\end{Remark}

Before we define the cyclotomic quotients of $\Waff_n$, which are the
main objects of study in this paper, we recall some standard
definitions and notation from the theory of Brauer algebras and some
of Nazarov's results.

A \textsf{Brauer diagram} on the $2n$ vertices $\{1,\dots,n,\bar
1,\dots,\bar n\}$ is a graph with $n$ edges such that each vertex lies
on a (unique) edge. Equivalently, a Brauer diagram is a partitioning
of $\{1,\dots,n,\bar1\dots,\bar n\}$ into $n$ two element subsets. Let
$\BDiag(n)$ be the set of all Brauer diagrams on
$\{1,\dots,n,\bar 1,\dots,\bar n\}$. Then
$\#\BDiag(n)=(2n-1)!!$.\label{Brauer diagrams}

Let $\gamma\in\BDiag(n)$ be a Brauer diagram. A \textsf{vertical} edge
in $\gamma$ is any edge of the form $\{m, \bar m\}$, where $1\le m
\le n$.  \textsf{Horizontal} edges are edges of the form $\{m,p\}$, or
$\{\bar m, \bar p\}$, where $1\le m<p\le n$.

For $i=1,\dots,n-1$ let $\gamma(i,i+1)$ be the Brauer diagram with
edges $\{i,\bar{i+1}\}$, $\{i+1,\bar i\}$ and all other edges being
vertical. Similarly, let $\gamma_i$ be the Brauer diagram with edges
$\{i,i+1\}$, $\{\bar i,\overline{i+1}\}$, and with all other edges
being vertical. We set $s_i=b_{\gamma(i,i+1)}$ and $e_i=b_{\gamma_i}$.
We also let $\gamma_e$ be the graph with edges
$\set{\{i,\bar i\}|1\le i\le n}$.

Brauer diagrams can be represented diagrammatically as in the following
examples. The vertices in the first rows are labelled from left to
right as $1$ to $4$, and the vertices in the second row are labelled
$\bar 1$ to $\bar 4$.

$$\gamma_e \quad =\quad
    \BrauerDiagram{&&&}%
       { \ncline{1,1}{2,1}\ncline{1,2}{2,2}\ncline{1,3}{2,3}\ncline{1,4}{2,4} },
   \quad \gamma(1,2) \quad =\quad
    \BrauerDiagram{&&&}%
       { \ncline{1,1}{2,2}\ncline{1,2}{2,1}\ncline{1,3}{2,3}\ncline{1,4}{2,4} },
   \quad\text{and}\quad \gamma_2 \quad =\quad
    \BrauerDiagram{&&&}%
       { \ncline{1,1}{2,1}\ncline{1,2}{1,3}\ncline{2,2}{2,3}\ncline{1,4}{2,4} }.
$$

Given two Brauer diagrams $\gamma, \gamma'\in\BDiag(n)$ we define their
product to be the diagram $\gamma\bullet\gamma'$ which is obtained by
identifying vertex $\bar i$ in $\gamma$ with vertex $i$ in $\gamma'$, for
$1\le i\le n$. Let $\ell(\gamma,\gamma')$ be the number of loops in the
graph $\gamma\bullet\gamma'$ and let $\gamma\circ\gamma'$ be the Brauer diagram obtained
by deleting these loops. The following pictures give two examples of the
multiplication $\gamma\circ\gamma'$ of diagrams.
$$\psset{dotsep=0.5mm}\BrauerDiagram{&&&}{\\&&&&\\&&&&%
     \ncline{1,1}{2,2}\ncline{1,2}{2,1}\ncline{1,3}{2,3}\ncline{1,4}{2,4}
     \ncline{3,1}{4,1}\ncline{3,2}{3,3}\ncline{4,2}{4,3}\ncline{3,4}{4,4}
     \ncline[linestyle=dotted]{2,1}{3,1}\ncline[linestyle=dotted]{2,2}{3,2}
     \ncline[linestyle=dotted]{2,3}{3,3}\ncline[linestyle=dotted]{2,4}{3,4} }
   \quad = \quad
   \BrauerDiagram{&&&}%
       { \ncarc[arcangle=50]{1,1}{1,3}\ncline{1,2}{2,1}\ncline{2,2}{2,3}\ncline{1,4}{2,4} }
\qquad\text{and}\qquad
  \BrauerDiagram{&&&}{\\&&&&\\&&&&%
     \ncline{1,1}{2,1}\ncline{1,2}{1,3}\ncline{2,2}{2,3}\ncline{1,4}{2,4}%
     \ncline{3,1}{4,1}\ncline{3,2}{3,3}\ncline{4,2}{4,3}\ncline{3,4}{4,4}%
     \ncline[linestyle=dotted]{2,1}{3,1}\ncline[linestyle=dotted]{2,2}{3,2}
     \ncline[linestyle=dotted]{2,3}{3,3}\ncline[linestyle=dotted]{2,4}{3,4} }
   \quad = \quad
  \BrauerDiagram{&&&}%
       { \ncline{1,1}{2,1}\ncline{1,2}{1,3}\ncline{2,2}{2,3}\ncline{1,4}{2,4} }
$$
In the first example $\gamma=\gamma(1,2)$,
$\gamma'=\gamma_2$ and $\ell(\gamma,\gamma')=0$. In the second example
$\gamma=\gamma'=\gamma_2$ and $\ell(\gamma,\gamma')=1$.

Recall that $R$ is a commutative ring.

\begin{Defn}[Brauer~\cite{BrauerAlg}]
Suppose that $\omega\in R$. The \textsf{Brauer algebra}
$\B_n(\omega)$, with parameter $\omega$, is the
$R$--algebra which is free as an $R$--module with basis
$\set{b_\gamma|\gamma\in\BDiag(n)}$ and with multiplication determined
by
$$b_\gamma
b_{\gamma'}=\omega^{\ell(\gamma,\gamma')}b_{\gamma\circ\gamma'},$$
for $\gamma,\gamma'\in\BDiag(n)$.
\end{Defn}

It is easy to see that $\B_n(\omega)$ is an associative algebra with
identity $b_{\gamma_e}$. We abuse notation and sometimes write
$1=b_{\gamma_e}$.

The second example above indicates that $e_i^2=\omega e_i$, for
$1\le i<n$. Similarly, $s_i^2=1$, for $1\le i<n$.

Let $\Sym_n$ be the symmetric group on $n$ letters. To each
permutation $w\in\Sym_n$ we associate the Brauer diagram $\gamma(w)$
which has edges $\set{\{i,\overline{w(i)}\}|\text{for } 1\le i\le n}$.
Notice that if $w=(i,i+1)$ then this is consistent with the notation
introduced above for the elements $s_i=b_{\gamma(i,i+1)}\in\B_n(\omega)$.

The diagrams $\set{\gamma(w)|w\in\Sym_n}$ are precisely the Brauer
diagrams which do not have any horizontal edges. It is easy to see
that the map $w\mapsto b_{\gamma(w)}$ induces an algebra embedding of
the group ring $R\Sym_n$ of $\Sym_n$ into $\B_n(\omega)$. In this way,
$R\Sym_n$ can be considered as a subalgebra of $\B_n(\omega)$.

There is a well--known presentation of $\B_n(\omega)$, which we now
describe. See \cite{MW00:BMW} for example.

\begin{Prop}\label{Brauer presentation}
Suppose that $R$ is a commutative ring. The Brauer algebra $\B_n(\omega)$ is
generated by the elements $s_1,\dots,s_{n-1},e_1,\dots,e_{n-1}$
subject to the relations
\begin{xalignat*}{3}
  s_i^2 &=1, & e_i^2 &=\omega e_i, & s_ie_i &=e_is_i=e_i,\\
s_is_j&=s_js_i,& s_ie_j&=e_js_i,& e_ie_j&=e_je_i,\\
  s_ks_{k+1}s_k&=s_{k+1}s_ks_{k+1},& e_ke_{k+1}e_k&=e_k,
                            & e_{k+1}e_ke_{k+1}&=e_{k+1},\\
  s_ke_{k+1}e_k&=s_{k+1}e_k, & e_{k+1}e_ks_{k+1}&=e_{k+1}s_k,
\end{xalignat*}
where $1\le i,j<n$, with $|i-j|>1$, and $1\le k<n-1$
\end{Prop}

Let $s_{ij}=b_{\gamma(i,j)}$, and let $e_{ij}=b_{\gamma_{ij}}$ where
$\gamma_{ij}$ is the Brauer diagram with edges $\{i,j\}$, $\{\bar i,
\bar j\}$ and $\{k, \bar k\}$, for $k\ne i,j$.

\begin{Cor}[\protect{Nazarov~\cite[(2.2)]{Nazarov:brauer}}]\label{Brauer surjection}
    Suppose that $\omega\in R$ and let $\Omega=\set{\omega_a|a\ge0}$, where
$\omega_a=\omega\big(\frac{\omega-1}2\big)^a$, for $a\ge0$. Then there
is a surjective algebra homomorphism
$\pi\map{\Waff_n(\Omega)}{\B_n(\omega)}$ which is determined by
$$\pi(S_i)=s_i,\quad \pi(E_i)=e_i,\quad\text{and}\quad
         \pi(X_j)=\frac{\omega-1}2 +\sum_{k=1}^{j-1}(s_{kj}-e_{kj}),$$
for $1\le i<n$ and $1\le j\le n$. Moreover,
$\ker\pi=\langle X_1-(\tfrac{\omega-1}{2})\rangle$, so that
$$\Waff_n(\Omega)/\big\langle X_1-(\tfrac{\omega-1}2)\big\rangle\cong\B_n(\omega).$$
\end{Cor}

Notice, in particular, that $\pi(X_1)=\tfrac{\omega-1}2$. To prove this
result it is enough to show that the elements $\pi(X_j)$, for $1\le j\le
n$, satisfy the relations in $\Waff_n(\Omega)$. For these calculations
see \cite[Lemma~2.1 and Proposition~2.3]{Nazarov:brauer}.

Fix a Brauer diagram $\gamma\in\BDiag(n)$. By
Proposition~\ref{Brauer presentation} we can write $b_\gamma$ as a
word in the generators $s_1,\dots,s_{n-1},e_1,\dots,e_{n-1}$. Fix such
a word for $b_\gamma$ and let $B_\gamma\in\Waff_n(\Omega)$ be the
corresponding word in the generators
$S_1,\dots,S_{n-1},E_1,\dots,E_{n-1}$. Then $\pi(B_\gamma)=b_\gamma$.
\label{B gamma}

Given $\alpha, \beta\in \N_0^n$ and $\gamma\in\BDiag(n)$ write
\label{X alpha}
$$
X^\alpha B_\gamma X^\beta
 =X_1^{\alpha_1}\dots X_n^{\alpha_n}B_\gamma X_1^{\beta_1}\dots X_n^{\beta_n}.
$$
We want to use these monomials to give a basis of $\Waff_n(\Omega)$.
The anti--symmetry relations $E_i(X_i+X_{i+1})=0$, for $1\le i<n$,
show that the set of all monomials is not linearly independent. In
Theorem~\ref{Waff basis} below we will show that the following
monomials are linearly independent.

\begin{Defn}\label{regular}
Suppose that $\alpha, \beta\in \N_0^n$ and $\gamma\in\BDiag(n)$.  A
monomial $X^\alpha B_\gamma X^\beta$ in~$\Waff_n(\Omega)$ is
\textsf{regular} if
\begin{enumerate}
\item $\alpha_r=0$ whenever $r$ is the left
endpoint of a horizontal edge in the top row of~$\gamma$.
\item if $\beta_l\ne0$ then $l$ is the left endpoint of a horizontal
edge in the bottom row of~$\gamma$.
\end{enumerate}
\end{Defn}

We can view a regular monomial $X^\alpha B_\gamma X^\beta$ as a Brauer
diagram if we colour the horizontal and vertical edges with the
non--negative integers using $\alpha$ and $\beta$.

Following Corollary~\ref{constraint} we also make the following
definition. (Recall that we are assuming that $2$ is invertible in $R$.)

\begin{Defn}\label{admissibility}
Let $\Omega=\set{\omega_a\in R|a\ge0}$. Then $\Omega$ is \textsf{admissible} if
$$
\omega_{2a+1}=\frac12\Big\{-\omega_{2a}
         +\sum_{b=1}^{2a+1}(-1)^{b-1}\omega_{b-1}\omega_{2a+1-b}\Big\},
$$
for all $a\ge0$.
\end{Defn}

By Corollary~\ref{constraint} $E_1$ is a torsion element if $\Omega$
is not admissible.

\begin{Remark}\label{W_1 remark}
    Let $y$ be an indeterminate and consider the generating
    series $\widetilde W_1(y)=\sum_{a\ge0}\omega_ay^{-a}$. Then the
    condition for $\Omega$ to be admissible can be written as
    $$\(\widetilde W_1(y)+y-\half\)\(\widetilde W_1(-y)-y-\half\)
            =(\half-y)(\half+y).$$
    Similar generating functions play an important
    role in section~4.
\end{Remark}

\begin{Theorem}[\protect{Nazarov~\cite[Theorem~4.6]{Nazarov:brauer}}]
\label{Waff basis}
Suppose $R$ is a commutative ring in which $2$ is a unit and that
$\Omega=\set{\omega_a\in R|a\ge0}$ is admissible. Then
$\Waff_n(\Omega)$ is free as an $R$--module with basis
$\set{X^\alpha B_\gamma X^\beta|\alpha,\beta\in\N_0^n, \gamma\in\BDiag(n),
\text{ and $X^\alpha B_\gamma X^\beta$ is regular}}$.
\end{Theorem}

\begin{proof}[Sketch of proof]
We have defined the elements of $\Omega$ to be scalars, but
Nazarov~\cite{Nazarov:brauer} works with a larger algebra
$\widehat{\Waff_n}(\widehat\Omega)$ generated by elements
$S_i$, $E_i$, $X_j$, for $1\le i<n$ and $1\le j\le
n$, and $\widehat\Omega=\set{\widehat\omega_a|a\ge0}$ where these
generators satisfy the same relations as the corresponding generators
of $\Waff_n(\Omega)$ except that the elements of $\Omega$ are
central elements of $\widehat{\Waff_n}(\widehat\Omega)$, rather
than scalars.  Hence,
$\Waff_n(\Omega)\cong\widehat{\Waff_n}(\widehat\Omega)/I$, where $I$
is the two sided ideal of $\widehat{\Waff_n}(\widehat\Omega)$
generated by the elements $\set{\widehat\omega_a-\omega_a|a\ge0}$.

Nazarov puts a grading on $\widehat{\Waff_n}(\widehat\Omega)$ by
setting $\deg S_i=\deg E_i=\deg \widehat\omega_a=0$ and $\deg X_i=1$.
To prove the result it is enough to work with the associated graded
algebra $\mathop{gr}(\widehat{\Waff_n}(\widehat\Omega))$, where the
grading is that induced by the degree function.  The arguments of
Lemma~4.4 and Lemma~4.5 from \cite{Nazarov:brauer} go through without
change for an arbitrary ring, so
$\widehat{\Waff_n}(\widehat\Omega)$ is spanned by
$$\Set[60]X^\alpha B_\gamma X^\beta
   {\widehat\omega}_2^{h_2}{\widehat\omega}_4^{h_4}\dots|%
         $\alpha, \beta\in\N_0^n$, $\gamma\in\BDiag(n)$, $h_{2i}\ge0$,
      for $i\ge1$, with only finitely many $h_{2i}\ne0$|,$$
where the monomials $X^\alpha B_\gamma X^\beta$ are all regular (see
\cite[Theorem~4.6]{Nazarov:brauer}). This implies that the regular
monomials span $\Waff_n(\Omega)$ for any ring $R$.

To complete the proof we first consider the case where the elements of
$\Omega'$ are indeterminates over $\Z$ and we consider the affine Wenzl
algebras defined over the field $\C(\Omega')$ and over the ring
$\Z[\Omega']$. We write $\Waff_{R,n}(\Omega')=\Waff_n(\Omega')$ to
emphasize that $\Waff_n(\Omega')$ is defined over the ring $R$.

Using Nazarov's algebra $\widehat{\Waff_n}(\widehat\Omega')$ and arguing as
above, it follows from \cite[Lemma~4.8]{Nazarov:brauer} that the
set of regular monomials are linearly independent when $R=\C(\Omega')$.
By the last paragraph, the regular monomials span
$\Waff_{\Z[\Omega'],n}(\Omega')$. Using the natural map
$\Waff_{\Z[\Omega'],n}(\Omega')\longrightarrow\Waff_{\C(\Omega'),n}(\Omega')$
it follows that
$\Waff_{\Z[\Omega'],n}(\Omega')$ is free as a $\Z[\Omega']$--module and
has basis the set of regular monomials. Hence, by a specialization
argument, if $R$ is arbitrary ring $R$ and $\Omega\subseteq R$ then
$$\Waff_{R,q}(\Omega)
      \cong\Waff_{\Z[\Omega'],n}(\Omega')\otimes_{\Z[\Omega']}R,$$
where we consider $R$ as a $\Z[\Omega']$--module by letting
$\omega_a'\in\Omega'$ act on $R$ as multiplication by
$\omega_a\in\Omega$, for $a\ge0$.  Hence, $\Waff_{R,n}(\Omega)$ is
free as an $R$--module with basis the set of regular monomials as
claimed.
\end{proof}

We are now ready to define the cyclotomic Nazarov--Wenzl algebras.  We assume
henceforth that $\Omega$ is admissible.

\begin{Defn}\label{cyclo NZ}
Fix an integer $r\ge1$ and $\bu=(u_1,\dots,u_r)\in R^r$.
The \textsf{cyclotomic Nazarov--Wenzl algebra} $\W_{r,n}=\W_{r,n}(\bu)$ is the
$R$--algebra $\Waff_n(\Omega)/\langle(X_1-u_1)\dots(X_1-u_r)\rangle$.
\end{Defn}

We should write $\W_{r,n}(\bu,\Omega)$, however, in section~3 we will
restrict to the case where $\Omega$ is $\bu$--admissible
(Definition~\ref{u-admissible}), which implies that $\omega_a$ is
determined by~$\bu$, for $a\ge0$.  For this reason we omit $\Omega$
from the notation for $\W_{r,n}(\bu)$.

By Corollary~\ref{Brauer surjection} the Brauer algebras
$\B_n(\omega)$ are a special case of the cyclotomic Nazarov--Wenzl algebras
corresponding to $r=1$ and
$\Omega=\set{\omega(\tfrac{\omega-1}2)^a|a\ge0}$.

By definition there is a surjection
$\pi_{r,n}\map{\Waff_n(\Omega)}{\W_{r,n}(\bu)}$. Abusing notation, we
write $S_i=\pi_{r,n}(S_i)$, $E_i=\pi_{r,n}(E_i)$, $X_j=\pi_{r,n}(X_j)$, and
$B_\gamma=\pi_{r,n}(B_\gamma)$ for $1\le i<n$, $1\le j\le n$ and
$\gamma\in\BDiag(n)$.

Notice that because $(X_1-u_1)\dots(X_1-u_r)=0$ in $\W_{r,n}(\bu)$ the
cyclotomic Nazarov--Wenzl algebras have only $r$ unwrapping relations; that is,
we only need to impose the relations
$E_1X_1^aE_1=\omega_aE_1$, for $0\le a\le r-1$.

Every $\W_{r,n}(\bu)$--module can be considered as a
$\Waff_n(\Omega)$--module by inflation along the surjection
$\pi_{r,n}\map{\Waff_n(\Omega)}{\W_{r,n}(\bu)}$. In particular, every
irreducible $\W_{r,n}(\bu)$--module is also an irreducible
$\Waff_n(\Omega)$--module.  Conversely, it is not hard to see that
every irreducible $\Waff_n(\Omega)$--module~$M$ over an algebraically
closed field can be considered as an irreducible module for some
cyclotomic Nazarov--Wenzl algebra~$\W_{r,n}(\bu)$, where $\bu$ depends
on $\Omega$ and ~$M$.  At first sight this is not very useful because
almost all of the results in this paper require that $\Omega$ be
$\bu$--admissible (Definition~\ref{u-admissible}) and, in general,  it
seems unlikely that $\Omega$ will be $\bu$--admissible for all of the
parameters $\bu$ that arise in this way.  Nevertheless, this
observation and the theory of cellular algebras allows us to construct
all of the finite dimensional $\Waff_n(\Omega)$ modules over an
algebraically closed field when $\Omega$ is admissible; see
Theorem~\ref{finite irreds}.

%

For our first result about the cyclotomic Nazarov--Wenzl algebras we
prove the easy half of Theorem~\ref{RANK}. That is, we show that
$\W_{r,n}(\bu)$ is spanned by $r^n(2n-1)!!$ elements.

\begin{Defn}
Suppose that $\alpha, \beta\in \N_0^n$ and $\gamma\in\BDiag(n)$.
\begin{enumerate}
    \item The monomial $X^\alpha B_\gamma X^\beta$ in $\W_{r,n}(\bu)$
        is \textsf{regular} if $X^\alpha B_\gamma X^\beta$ is a
        regular monomial in $\Waff_n(\Omega)$.
    \item The monomial $X^\alpha B_\gamma X^\beta$ in $\W_{r,n}(\bu)$
        is \textsf{$r$--regular} if it is regular and
        $0\le \alpha_i, \beta_i<r$, for all $1\le i\le n$.
\end{enumerate}
\end{Defn}

\begin{Prop}\label{spanning set}
The cyclotomic Nazarov--Wenzl algebra $\W_{r,n}(\bu)$ is spanned by the
set of $r$--regular monomials $\{X^\alpha B_\gamma X^\beta\}$. In
particular, if $R$ is a field then
$$\dim_R\W_{r,n}(\bu)\le r^n(2n-1)!!.$$
\end{Prop}

\begin{proof} By Theorem~\ref{Waff basis}, and the definitions,
$\W_{r,n}(\bu)$ is spanned by the regular monomials in
$\W_{r,n}(\bu)$.  As in the proof of Theorem~\ref{Waff basis}, we put
a grading on $\W_{r,n}(\bu)$. Then in the associated graded algebra,
$\mathop{gr}\W_{r,n}(\bu)$, we have the relation $(X_i-u_1)\cdots
(X_i-u_r)=0$. We claim that the regular monomial $X^\alpha B_\gamma
X^\beta$ can be written as a linear combination of $r$-regular
monomials. If $X^\alpha B_\gamma X^\beta$ is an $r$--regular monomial
then there is nothing to prove so we may assume that $X^\alpha
B_\gamma X^\beta$ is not $r$--regular and, in particular, that
$|\alpha|+|\beta|>0$. Then, using the relation $(X_i-u_1)\cdots
(X_i-u_r)=0$ we can subtract a linear combination of $r$--regular
monomials from $X^\alpha B_\gamma X^\beta$ to obtain a linear
combination of regular elements of smaller degree.  The claim now
follows by induction.

Finally, a counting argument shows that the number of $r$--regular
monomials is equal to $r^n(2n-1)!!$. Therefore, if $R$ is a field then
$\dim_R\W_{r,n}(\bu)\le r^n(2n-1)!!$.
\end{proof}

The \textsf{degenerate (cyclotomic) Hecke algebra} $\H_{r,n}(\bu)$ of type
$G(r,1,n)$ is the unital associative $R$--algebra with generators
$T_1,\dots,T_{n-1},Y_1,\dots,Y_n$ and relations
\begin{alignat*}{3}
  (Y_1-u_1)\dots(Y_1-u_r)&=0, & T_i^2 &=1,\\
  T_iT_j &=T_jT_i, &Y_iY_k &=Y_kY_i, \\
  T_iY_i-Y_{i+1}T_i&=-1,&  Y_iT_i-T_iY_{i+1}&=-1, \\
  T_jT_{j+1}T_j &=T_{j+1}T_jT_{j+1},& T_iY_\ell=Y_\ell T_i,
\end{alignat*}\label{degenerate H}
for $1\le i<n$, $1\le j<n-1$ with $|i-j|>1$, $1\le k\le n$ and $\ell\ne i,i+1$.
Therefore there is a surjective algebra homomorphism
$\W_{r,n}(\bu)\longrightarrow\H_{r,n}(\bu)$ determined by $$S_i\mapsto
T_i,\quad E_i\mapsto 0, \quad\text{and}\quad X_j\mapsto Y_j,$$ for
$1\le i<n$ and $1\le j\le n$.  (In fact, a special case of
Proposition~\ref{degen isom} below shows that
$\H_{r,n}(\bu)\cong\W_{r,n}(\bu)/\langle E_1\rangle$.) Consequently,
every irreducible $\H_{r,n}(\bu)$--module can be considered as an
irreducible $\W_{r,n}(\bu)$--module via inflation.  These irreducible
modules are precisely the irreducibles upon which $E_i$ acts as zero.
We record this fact for future use.

\begin{Cor}\label{AK irreds}
Suppose that $R$ is a field and that $M$ is an irreducible
$\W_{r,n}(\bu)$--module which is annihilated by some $E_i$. Then $M$ is
an irreducible $\H_{r,n}(\bu)$--module.
\end{Cor}

\begin{proof} As $E_{i+1}=S_iS_{i+1}E_iS_{i+1}S_i$ and $S_j$ is
    invertible for all $j$, the two--sided
    ideal of $\W_{r,n}(\bu)$ generated by $E_1$ is the same as the
    two--sided ideal generated by $E_i$, for $1\le i<n$. The result
    now follows from the remarks
    above.
\end{proof}

Recall that the degenerate affine Hecke algebra is a finitely
generated module over its center (see, for example,
\cite{Klesh:book}), which is the ring of the symmetric polynomials in
$Y_1,\dots,Y_n$.  This fact, together with Dixmier's version of
Schur's lemma, implies that all of the irreducible modules of the
degenerate affine Hecke algebra are finite dimensional. 

By the last paragraph the power sum symmetric functions are central
elements of the degenerate affine Hecke algebra. In contrast, only the
power sums of odd degree are central in $\Waff_n(\Omega)$. Another
difference is that the affine Wenzl algebra is not finitely generated
over its center. To see this, we give an example of an infinite
dimensional irreducible $\Waff_2(\Omega)$-module.

\begin{Example}
Suppose that $\Omega$ is admissible and that 
$\tilde W_1(y)=\sum_{a\ge0}\omega_a y^{-a}$ is not a rational 
function in $y$. Consider
$V=\oplus_{n\ge0}Rv_n$. Define an action of $\Waff_2(\Omega)$ on $V$ by
$Ev_n=\omega_nv_0$, $X_1v_n=v_{n+1}$, $X_2v_n=-v_{n+1}$ and
$$
Sv_n=(-1)^nv_n- \varepsilon v_{n-1}+ \sum_{k=0}^{n-1}(-1)^k\omega_{n-k-1}v_k,
$$
where $\varepsilon=1$, if $n\equiv1\pmod 2$, and $\varepsilon=0$,
otherwise.  All of the defining relations except for the relation $S^2=1$
are easy to check.  As $S^2$ commutes with $X_1$, $S^2v_0=v_0$ and
$X_1v_n=v_{n+1}$, we have that $S^2$ acts as the identity on $V$.

Now we show that~$V$ is irreducible. Let $W$ be a 
$\Waff_2(\Omega)$-submodule of $V$.  Suppose that $EW=0$.
If $\sum c_nv_n\in W$ then $\sum c_n\omega_{n+k}=0$, for all 
$k\ge0$. As the vectors 
$\{(\omega_k,\omega_{k+1},\dots)\in R^\infty|k\ge0\}$ 
span an infinite dimensional subspace of $R^\infty$, 
we have $c_n=0$, for all $n\ge0$. Hence $W=0$. 
Thus, $W\ne0$ implies $EW\ne0$. Then $v_0\in W$ and $W=V$.
Therefore,$V$ is an irreducible $\Waff_2(\Omega)$-module as claimed.
\end{Example}

In light of this example, we restrict our attention to finite
dimensional $\Waff_n(\Omega)$--modules in what follows.

\section{Restrictions on $\Omega$ and the irreducible representations
of $\W_{r,2}$}

In this section we explicitly compute the (possible) irreducible
representations of the cyclotomic Nazarov--Wenzl algebras $\W_{r,2}(\bu)$. As a
consequence we find a set of conditions on the parameter set $\Omega$
which ensure that $\W_{r,2}(\bu)$ has dimension
$3r^2=r^n(2n-1)!!\mid_{n=2}$ when $R$ is a field. In the next section
we will see that these conditions on $\Omega$ are exactly what we need
for general $n$.

The cyclotomic Nazarov--Wenzl algebra $\W_{r,2}(\bu)$ is generated by
$S_1,E_1,X_1$ and $X_2$. Throughout this section we suppose that $R$
is an algebraically closed field and, for convenience, we set $S=S_1$
and $E=E_1$.

\begin{Prop}\label{E=0} Suppose that $M$ is an irreducible
$\W_{r,2}(\bu)$--module such that $EM=0$. Then either:
\begin{enumerate}
\item $M=Rm$ is one dimensional and the action of $\W_{r,2}(\bu)$
        is determined by
$$Sm=\epsilon m,\quad Em=0, \quad X_1m=u_im,\quad\text{and}\quad
               X_2m=(u_i+\epsilon)m,$$
where $\epsilon=\pm1$ and $1\le i\le r$. In particular, up to
isomorphism, there are at most $2r$ such representations.
\item $M$ is two dimensional and the action of $\W_{r,2}(\bu)$ is given by
    $$S\mapsto\tfrac1{u_i-u_j}\big(\begin{smallmatrix}-1&b\\c&1\end{smallmatrix}\big),\quad
      E\mapsto\big(\begin{smallmatrix}0&0\\0&0\end{smallmatrix}\big),\quad
      X_1\mapsto\big(\begin{smallmatrix}u_i&0\\0&u_j\end{smallmatrix}\big),
      \quad\text{and}\quad
      X_2\mapsto\big(\begin{smallmatrix}u_j&0\\0&u_i\end{smallmatrix}\big),$$
for some non--zero $b,c\in R$ such that $bc=(u_i-u_j)^2-1$, where $u_i\ne u_j$.
Up to isomorphism there are at most $\binom r2$ such representations.
\item $M$ is two dimensional and the action of $\W_{r,2}(\bu)$ is given by
$$S\mapsto\big(\begin{smallmatrix}0&1\\1&0\end{smallmatrix}\big),\quad
      E\mapsto\big(\begin{smallmatrix}0&0\\0&0\end{smallmatrix}\big),\quad
      X_1\mapsto\big(\begin{smallmatrix}u_i&-1\\0&u_i\end{smallmatrix}\big),
      \quad\text{and}\quad
      X_2\mapsto\big(\begin{smallmatrix}u_i&1\\0&u_i\end{smallmatrix}\big).$$
Up to isomorphism there are at most $r$ such representations.
\end{enumerate}
\end{Prop}

\begin{proof}
As noted in Corollary~\ref{AK irreds} $M$ is an irreducible
$\H_{r,2}(\bu)$--module. The result now follows from the
representation theory of $\H_{r,2}(\bu)$: choose a simultaneous
eigenvector $m$ of $R[Y_1,Y_2]$. Then, because
$\H_{r,2}(\bu)=R[Y_1,Y_2]+T_1R[Y_1,Y_2]$, if $M$ is not one
dimensional then it must be two dimensional.  If this is the case,
$\{m,Sm\}$ is a basis of $M$.  Further, if the eigenvalues for the
action $Y_1$ on $M$ are distinct, then we can simultaneously
diagonalize $Y_1$ and $Y_2$. All of our claims now follow.
\end{proof}

Note that since $\prod_{i=1}^r(X_1-u_i)$ acts as zero on $M$, case~(c)
can only arise if the~$u_i$ are not pairwise distinct. The irreducible
representations of $\W_{r,2}(\bu)$ upon which~$E$ acts non--trivially
take more effort to understand.

\begin{Prop}\label{Ene0}
Let $F$ be a field in which $2$ is invertible and that
$u_1,\dots,u_r$ are algebraically independent over $F$. Let
$R=F(u_1,\dots,u_r)$  and let $\W_{r,2}(\bu)$ be the cyclotomic Nazarov--Wenzl
algebra defined over $R$, where $\omega_0\ne0$.  Then $\W_{r,2}(\bu)$
has a unique irreducible module $M$ such that $EM\ne0$. Moreover, if
$d=\dim_R M$ then $d\le r$ and there exists a basis
$\{m_1,\dots,m_d\}$ of~$M$ and scalars
$\{v_1,\dots,v_d\}\subseteq\{u_1,\dots,u_r\}$, with $v_i\ne v_j$
when~$i\ne j$, such that for $1\le i\le d$ the following hold:
\begin{enumerate}
    \item $X_1 m_i=v_im_i$ and $X_2 m_i=-v_im_i$,
    \item $E m_i= \gamma_i (m_1+\dots+m_d)$ and
    \item $\displaystyle S m_i= \frac{\gamma_i-1}{2v_i}m_i
        +\sum_{j\ne i}\frac{\gamma_i}{v_i+v_j}m_j$,
\end{enumerate}\noindent%
where $\gamma_i=(2v_i-(-1)^d)\displaystyle
        \prod_{\substack{1\le j\le d\\j\ne i}}\dfrac{v_i+v_j}{v_i-v_j}$.
Moreover, $\omega_a=\sum_{j=1}^d v_j^a\gamma_j$, for all $a\ge0$; and,
in particular,
$$\omega_0=\begin{cases}
              2(v_1+\dots+v_d),&\text{if $d$ is even},\\
              2(v_1+\dots+v_d)+1,&\text{if $d$ is odd}.\\
          \end{cases}$$
Conversely, if $\omega_a=\sum_{j=1}^d v_j^a\gamma_j$, for all $a\ge0$,
then (a)--(c) define a $\W_{r,2}(\bu)$--module~$M$ with $EM\ne0$.
\end{Prop}

\begin{proof}
Suppose that $M$ is an irreducible $\W_{r,2}$--module such that
$EM\ne0$. Note that $M$ is finite dimensional.
Let $d=\dim_R M$. We first show that (a)--(c) hold.
%
%
Since $u_1,\dots,u_r$ are pairwise
distinct, we can fix a basis $\{m_1,\dots,m_d\}$ of $M$ consisting of
eigenvectors for $X_1$. Write $X_1m_i=v_i m_i$, for some $v_i\in
\{u_1, \cdots, u_r\}$.

Set $f:=\frac1{\omega_0}E$. This is a non--zero idempotent and
$fM\ne0$ since $EM\ne0$.

Fix an element $0\ne m\in fM$. Then $Em=\omega_0 m$ and $Sm=m$ (since
$SE=E$).  As $0=(X_1+X_2)Em=(X_1+X_2)\omega_0m$, we have
$(X_1+X_2)m=0$. However, $X_1+X_2$ is central in $\W_{r,2}$, so
$X_1+X_2$ acts as a scalar on $M$ by Schur's lemma. Hence,
$X_2m_i=-X_1m_i=-v_i m_i$, for $i=1,\dots,d$, proving (a).

We claim that $\{m,X_1m,\dots,X_1^{d-1}m\}$ is a basis
of $M$. To see this, for any~$a\ge0$ let $M_a$ be the $R$--submodule
of $M$ spanned by $\{m,X_1m,\dots,X_1^am\}$. Notice that $M_a$ is
closed under left multiplication by $E$ since if $k\ge0$ then
$$EX_1^km=EX_1^kfm=\frac1{\omega_0}EX_1^kEm
         =\frac{\omega_k}{\omega_0}Em
         =\omega_k m.$$
Also, by Lemma~\ref{SX^a},
\begin{align*}
    SX_1^am&=\big(X_2^aS+\sum_{b=1}^aX_2^{b-1}(E-1)X_1^{a-b}\big)m\\
      &=X_2^am+\sum_{b=1}^a\big(X_2^{b-1}EX_1^{a-b}E\frac1{\omega_0}m
-X_1^{a-b}X_2^{b-1}m\big)\\
      &=X_2^am+\sum_{b=1}^a\big(\frac{\omega_{a-b}}{\omega_0}X_2^{b-1}Em
-X_1^{a-b}X_2^{b-1}m\big)\\
      &=(-1)^a X_1^am+
\sum_{b=1}^a\big((-1)^{b-1}\omega_{a-b}X_1^{b-1}m-(-1)^{b-1}X_1^{a-1}m\big)
      \end{align*}
So, $M_a$ is closed under multiplication by $S$. Choose $a\ge0$ to be
minimal such that $\{m,X_1m,\dots,X_1^{a+1}m\}$ is not linearly
independent. Since $X_1^{a+1}m\in M_a$, $M_a$ is closed under
multiplication by $X_1$.  Hence, $M_a=M$ since $M$ is irreducible. By
counting dimensions, $M=M_{d-1}$, proving the claim.

Next we show that $EM=Rm$.
Suppose that $m'=\sum_{i=0}^{d-1}c_iX_1^im\in EM$.
Then
$$
m'=\frac1{\omega_0}Em'=\frac1{\omega_0}\sum_{i=0}^{d-1}c_iEX_1^im
             =\frac1{\omega_0^2}\sum_{i=0}^{d-1}c_iEX_1^iEm
             =\frac1{\omega_0^2}\Big(\sum_{i=0}^{d-1}c_i\omega_i\Big)m,
$$
since $Ea=\omega_0a$ whenever $a\in EM$. Hence, $EM=Rm$, as claimed.

Recall that we have fixed a basis $\{m_1,\dots,m_d\}$ of $M$.
Write $m=\sum_{i=1}^dr_im_i$, for some $r_i\in R$. Suppose that $r_i=0$ for
some $i$. Then
$$\prod_{\substack{1\le j\le d\\j\ne i}}(X_1-v_j)\cdot m=0.$$
This contradicts the linear independence of
$\{m,X_1m,\dots,X_1^{d-1}m\}$; hence, $r_i\ne0$ for $i=1,\dots,d$.
By rescaling the $m_i$, if necessary, we can and do
assume that $m=m_1+\dots+m_d$ in the following.

By the argument of the last paragraph, all of the eigenvalues
$\{v_1,\dots,v_d\}$ of $m$ must be distinct. This also shows that
$d=\dim M\le r$ and that $\{v_1,\dots,v_d\}$ are algebraically
independent (since we are assuming that $u_1,\dots,u_r$ are
algebraically independent).  In particular, $v_i$ and $v_i+v_j$, for
$i\ne j$, are invertible. So the formula in part (c) makes sense.

As $EM=Rm$, we can define elements $\gamma_i\in R$ by
$$
Em_i=\gamma_im=\gamma_i(m_1+\dots+m_d),\qquad\text{for } i=1,\dots,d.
$$
Write $Sm_i=\sum_{j=1}^d c_j^{(i)}m_j$.
Then $X_1Sm_i-SX_2m_i=(E-1)m_i$ reads
$$
\sum_{j=1}^d c_j^{(i)}v_jm_j+v_i(\sum_{j=1}^d c_j^{(i)}m_j)=
\gamma_i(m_1+\cdots+m_d)-m_i.
$$
Thus, $(v_i+v_j)c_j^{(i)}=\gamma_i-\delta_{ij}$ and we have
$$Sm_i=\frac{\gamma_i-1}{2v_i}m_i+\sum_{j\ne i}\frac{\gamma_i}{v_i+v_j}m_j.$$
This proves (c).

Next we prove the formula for $\gamma_i$ given in (b).
Since $E=SE$ we find that
$$\gamma_i\sum_{j=1}^d m_j=Em_i=SEm_i
       =\gamma_i\sum_{j=1}^d\Big\{\frac{\gamma_j-1}{2v_j}
            +\sum_{k\ne j}\frac{\gamma_k}{v_j+v_k}\Big\}m_j,$$
for $i=1,\dots,r$. Note that some $\gamma_i$ is non--zero, since
$EM\ne 0$. Thus, comparing the coefficient of $m_j$ on both sides
shows that
$$\sum_{k=1}^d\frac{\gamma_k}{v_j+v_k}=1+\frac1{2v_j},$$
for $j=1,\dots,d$.

We claim that
$\displaystyle\det\big(\frac1{v_i+v_j}\big)_{1\le i,j\le d}
  =\Big(\prod_{i=1}^d 2v_i\Big)^{-1}\prod_{i>j}\big(\frac{v_i-v_j}{v_i+v_j}\big)^2$.
To see this, observe that
$$ \big(\prod_{i=1}^d 2v_i)\prod_{i>j}(v_i+v_j)^2
         \det\big(\frac1{v_i+v_j}\big)_{1\le i,j\le d}
$$
is a symmetric polynomial in $v_1,\dots, v_d$ which is divisible
by $v_i-v_j$ for $i\ne j$. This shows that
this determinant is a constant multiple of $\prod_{i>j}\big(v_i-v_j\big)^2$.
To determine the constant, we multiply
$\det\big(\frac1{v_i+v_j}\big)_{1\le i,j\le d}$ by $v_n$,
set $v_n=\infty$ and use induction.

By the last paragraph, the matrix
$\big(\frac1{v_i+v_j}\big)_{1\le i,j\le d}$ is invertible, so
$\gamma_1,\dots,\gamma_d$ are uniquely determined. Hence, to prove the
formula for $\gamma_i$ it suffices to show that
$$
\sum_{k=1}^d \frac{2v_k-(-1)^d}{v_j+v_k}\prod_{i\neq k} \frac{v_k+v_i}{v_k-v_i}
=1+\frac{1}{2v_j},
$$
for $1\le j\le d$. Let $f(z)=\frac{2z-(-1)^d}{2z(z+v_j)}\prod_{i=1}^d
\frac{z+v_i}{z-v_i}$ and view $f(z)$ as an element of the function
field of the projective line defined over $F(v_1,\dots,v_d)$. Then,
the left hand side can be interpreted as the sum
$\sum_{k=1}^d\Res_{z=v_k} f(z)dz$, where $\Res_{z=v}f(z)dz$ is the
residue of~$f(z)$ at $v$, if $v\ne\infty$, and it is the residue of
$-\frac 1{z^2}f(\frac 1z)$ at $0$, if $v=\infty$.  Thus, the residue
theorem for complete non--singular curves implies that
$$ \sum_{k=1}^d \frac{2v_k-(-1)^d}{v_j+v_k}\prod_{i\neq k}
        \frac{v_k+v_i}{v_k-v_i} = -\Big(\Res_{z=\infty} f(z)dz+\Res_{z=0}
        f(z)dz\Big)=1+\frac{1}{2v_j}, $$
as required. Hence, we have shown that, for $1\le j\le d$,
$$\gamma_j=(2v_j-(-1)^d) \prod_{k\neq j} \frac{v_j+v_k}{v_j-v_k}\,,$$
so (b) is proved. (For a combinatorial proof see
Proposition~\ref{identity}(a) below.)

Now, since $Em=\omega_0m$ and $m=\sum_{i=1}^d m_i$, we have that
$\omega_0=\sum_{i=1}^d \gamma_i$. Similarly, we have that
$\omega_a=\sum_{j=1}^m v_j^a \gamma_j$ because
\begin{align*}
    \omega_am&=\frac{\omega_a}{\omega_0}Em=\frac1{\omega_0}EX_1^aEm
              =EX_1^am\\
             &=\sum_{i=1}^d EX_1^am_i
              =\sum_{i=1}^d v_i^aEm_i
              =\Big(\sum_{i=1}^d v_i^a\gamma_i\Big)m.
\end{align*}

We now show that
$$\omega_0=\sum_{i=1}^d\gamma_i
     =\begin{cases}
        2(v_1+\dots+v_d),&\text{if $d$ is even},\\
        2(v_1+\dots+v_d)+1,&\text{if $d$ is odd}.\\
       \end{cases}
$$
To evaluate $\sum_{i=1}^d \gamma_i$, we consider
$g(z)=\frac{2z-(-1)^d}{2z}\prod_{i=1}^d \frac{z+v_i}{z-v_i}$
and interpret the sum as $\sum_{i=1}^d \Res_{z=v_i} g(z)dz$.
Then the residue theorem gives the desired formula for $\omega_0$.

We next show that $M$ is uniquely determined, up to isomorphism.
Suppose that $\W_{r,n}(\bu)$ has another irreducible module of
dimension $d'$ upon which $e$ acts non--trivially. Then, by the
argument above,
$$\omega_0=\begin{cases}
         2(v_1'+\dots+v'_{d'}),&\text{if $d'$ is even},\\
         2(v_1'+\dots+v'_{d'})+1,&\text{if $d'$ is odd},\\
       \end{cases}
$$
for some $v_1',\dots,v_{d'}'\subseteq\{u_1,\dots,u_r\}$. As we are
assuming that $u_1,\dots,u_r$ are algebraically independent, this
forces $d'=d$ and $v_i'=v_{(i)\sigma}$, for some
$\sigma\in\Sym_d$ and $1\le i\le d$. Hence, by (a)--(c),
$M\cong M'$ as required.

Finally, it remains to verify that (a)--(c) define a representation
of $\W_{r,2}(\bu)$ whenever $\omega_a=\sum_{i=1}^dv_i^a\gamma_i$, for
$a\ge0$ and $\gamma_i$ as above. It is easy to check that the action
respects the relations
$E(X_1+X_2)=0=(X_1+X_2)E$, $EX_1^aE=\omega_aE$ and
$X_1S-SX_2=E-1=SX_1-X_2S$. That $SE=E$ and $ES=E$ on $M$,
follows from the identity $\sum_{k=1}^d
\frac{\gamma_k}{v_j+v_k}=1+\frac1{2v_j}$ proved above. We now prove
that $S^2=1$.  Observe that $S^2$ commutes with $X_1$ when acting on
$M$. As the $v_i$
are pairwise distinct, we have $S^2m_i=c_im_i$, for some $c_i\in R$.
Explicit computation shows that $c_i=\frac{1-2\gamma_i}{4v_i^2}+
\gamma_i\sum_{j=1}^d \frac{\gamma_j}{(v_i+v_j)^2}$.  Computing
the residues of $h(z)dz$, when
$h(z)=\frac{2z-(-1)^d}{2z(z+v_i)^2}\prod_{k=1}^d \frac{z+v_i}{z-v_i}$,
proves that $c_i=1$, for $1\le i\le d$.
\end{proof}

\begin{Remark}
The action of $X_1$ on an irreducible $\W_{r,2}(\bu)$-module
is not semisimple in general. For example, let $\Omega$ be given by
$\omega_0=1$, $\omega_1=0$ and
$$\omega_{a+2}=\frac1{2}\omega_{a+1}-\frac1{16}\omega_a,$$
for $a\ge0$. For $r=2$ we set
$(u_1,u_2)=(\frac1{4},\frac1{4})$ and for
$r=3$ set $(u_1,u_2,u_3)=(\frac1{4},\frac1{4},-\frac1{2})$.
Then $\W_{r,2}(\bu)$ has a two dimensional irreducible module upon
which the generators act as follows:
$$
E\mapsto\Big(\begin{smallmatrix}1&0\\0&0\end{smallmatrix}\Big),\quad
S\mapsto\Big(\begin{smallmatrix}1&0\\0&-1\end{smallmatrix}\Big),\quad
X_1\mapsto\Big(\begin{smallmatrix}0&\frac14\\-\frac14&\frac12\end{smallmatrix}\Big)
\quad\text{and}\quad
X_2\mapsto\Big(\begin{smallmatrix}0&-\frac14\\\frac14&-\frac12\end{smallmatrix}\Big).
$$
Further, $X_1-\frac1{4}\ne0$ and $(X_1-\frac1{4})^2=0$.
\end{Remark}

\begin{Theorem}\label{admissible condition}
Let $F$ be a field in which $2$ is invertible and that $u_1,\dots,u_r$
are algebraically independent over $F$. Let $R=F(u_1,\dots,u_r)$ and
suppose that $\W_{r,2}(\bu)$ is a split semisimple $R$--algebra and
that $\omega_0\ne0$. Then $\W_{r,2}(\bu)$ has dimension
$3r^2=r^n(2n-1)!!\mid_{n=2}$ if and only if $\W_{r,2}(\bu)$ has an
irreducible representation of dimension~$r$ upon which $E$ acts
non--trivially.
\end{Theorem}

\begin{proof}
We have constructed all the irreducible $\W_{r,2}$-modules in
Propositions~\ref{E=0} and \ref{Ene0} above. Under our assumptions,
Proposition~\ref{E=0} implies that $\W_{r,2}(\bu)$ has (a) $2r$
pairwise non--isomorphic one dimensional representations and (b)
$\binom r2$ pairwise non--isomorphic two dimensional representations.
Note that case (c) from Proposition~\ref{E=0} does not occur since
$u_1,\dots,u_r$ are pairwise distinct. Further, Proposition~\ref{Ene0}
implies that $\W_{r,2}(\bu)$ has a unique irreducible representation
$M$ such that $EM\ne0$ and, moreover, if $d=\dim M$ then $1\le d\le r$.
Hence, by the Wedderburn--Artin theorem we have
$$\dim\W_{r,2}(\bu)=2r+4\tbinom r2+d^2=2r^2+d^2,$$
so that $\dim\W_{r,2}(\bu)=3r^2$ if and only if $r=d$. The result follows.
\end{proof}

Note that by Proposition~\ref{Ene0}, under the conditions of the theorem, $\W_{r,2}(\bu)$ has a
(unique) representation of dimension~$r$ upon which $E$ acts
non--trivially if and only if $\omega_a=\sum_{j=1}^r u_j^a\gamma_j$,
for $a\ge0$, where
$\displaystyle\gamma_i=(2u_i-(-1)^r)\prod_{\substack{1\le j\le r\\j\ne i}}
         \dfrac{u_i+u_j}{u_i-u_j}$.

Recall that Schur's $q$--functions $q_a=q_a(\mathbf{x})$ in the
indeterminates $\mathbf{x}=(x_1,\dots,x_r)$~\cite[p.~250]{Macdonald}
are defined by the equation
$$\prod_{i=1}^r\frac{1+x_iy}{1-x_iy}=\sum_{a\ge0}q_a(\mathbf{x})y^a.$$
Note that $q_a(\mathbf{x})$ is a polynomial in $\mathbf{x}$, for all
$a\ge0$.

\begin{Lemma}\label{omegas}
Assume that $R$ is an integral domain and that $2$ is invertible in
$R$. Suppose that $\bu\in R^r$, with $u_i- u_j\ne0$ whenever $i\ne j$.
Let $F$ be the quotient field of $R$ and for $a\ge0$ define
$$\omega_a=\sum_{i=1}^r\(2u_i-(-1)^r\)u_i^a
    \prod_{\substack{1\le j\le r\\j\ne i}}\dfrac{u_i+u_j}{u_i-u_j}\in F,$$
as in Theorem~\ref{admissible condition}.
Then $\omega_a=q_{a+1}(\bu)-\half (-1)^r q_a(\bu)+\half\delta_{a0}$. In
particular, $\omega_a\in R$.
\end{Lemma}

\begin{proof}If $a=0$ then the result follows from
    Proposition~\ref{Ene0}, so we can assume that $a>0$.
    Let $f(z)=\half z^{a-1}(2z-(-1)^r)\prod_{i=1}^r\frac{z+u_i}{z-u_i}$.
    Then $\omega_a$ can be interpreted as
    $\sum_{i=1}^r\Res_{z=u_i}f(z)\,dz=-\Res_{z=\infty}f(z)\,dz$.
    Calculating the residue of $f(z)dz$ at $z=\infty$ now shows that
    $\omega_a=q_{a+1}(\bu)-\half (-1)^r q_a(\bu)+\half\delta_{a0}$.
    (See \cite[(2.9), p.~209]{Macdonald} for a more direct proof.)
    Hence, $\omega_a\in R$ since $q_b(\mathbf{x})\in R[\mathbf x]$, for
    $b\ge0$.
\end{proof}

We want the cyclotomic Nazarov--Wenzl  algebras to be ``cyclotomic''
generalizations of the Brauer algebras. In particular, we want them to
be free $R$--modules of rank $r^n(2n-1)!!$. Theorem~\ref{admissible
condition} gives sufficient conditions on
$\Omega=\set{\omega_a|a\ge0}$ for $\W_{r,2}(\bu)$ to have dimension
$r^n(2n-1)!!$ when $R$ is an algebraically closed field and $n=2$.
Consequently, in our study of $\W_{r,n}(\bu)$ we will require that
$\Omega$ have the following property.

\begin{Defn}
\label{u-admissible}
Let $\Omega=\set{\omega_a|a\ge0}\subseteq R$ and suppose that
$\bu\in R^r$. Then $\Omega$ is
\textsf{$\bu$--admissible} if
    $\omega_a=q_{a+1}(\bu)-\half (-1)^r q_a(\bu)+\half\delta_{a0},$
for $a\ge0$.
\end{Defn}

\begin{Remark}
Let $R=\Z[\bu]$ where $u_1,\dots,u_r$ are indeterminates.  Assume that
each $\omega_a$, for $a\ge0$, is a polynomial in $\bu$ and that
$\omega_0\ne0$. Then, by Theorem~\ref{admissible condition} and
Theorems~\ref{generic} and ~\ref{W cellular} below, $\Omega$ is
$\bu$-admissible if and only if
\begin{enumerate}
\item $\W_{r,n}(\Omega)\otimes_{\Z[\bu]}\Q(\bu)$ is
semisimple, and,
\item $\W_{r,n}(\Omega)$ is a free $R$-module of rank $r^n(2n-1)!!$,
\end{enumerate}%
\noindent for all $n\ge0$.
\end{Remark}

Recall from Remark \ref{W_1 remark} that 
$\displaystyle\widetilde W_1(y)=\sum_{a\ge0}\omega_ay^{-a}$, where $y$ is an
indeterminate.

\begin{Lemma}\label{W_1 expansion}
Suppose that $\bu\in R^r$.  Then $\Omega$ is $\bu$--admissible if and only if
$$\widetilde W_1(y)+y-\half
              =(y-\half(-1)^r)\prod_{i=1}^r\frac{y+u_i}{y-u_i}.$$
\end{Lemma}

\begin{proof} By definition, $\u$--admissibility is equivalent to the
identity
$$\widetilde W_1(y)=
     \half+\sum_{a\ge0}\(q_{a+1}(\bu)-\half (-1)^r q_a(\bu)\)y^{-a}.$$
Now expand this equation using the definition of the Schur $q$--functions.
\end{proof}

\begin{Cor}\label{admiss}
    Suppose that $\Omega$ is $\bu$--admissible. Then $\Omega$ is admissible.
\end{Cor}

\begin{proof}First suppose that $\mathbf{x}=(x_1,\dots,x_r)$ are
    algebraically independent and let $\Omega=\set{\omega_a|a\ge0}$,
    where $\omega_a=q_{a+1}(\mathbf x)-\half (-1)^r
    q_a(\mathbf{x})+\half\delta_{a0},$ for $a\ge0$.  Then $\Omega$ is
    $\mathbf{x}$--admissible by definition and hence admissible by
    Corollary~\ref{constraint} and Proposition~\ref{Ene0}. Therefore,
    by the definition of admissibility we have the following
    polynomial identity in $x_1,\dots,x_r$
    $$
    \omega_{2a+1}=\frac12\Big\{-\omega_{2a}
         +\sum_{b=1}^{2a+1}(-1)^{b-1}\omega_{b-1}\omega_{2a+1-b}\Big\}.
   $$
   The general case now follows by specializing
   $x_i=u_i$, for $1\le i\le r$.

   For a second proof, note that if $\Omega$ is $\bu$--admissible then
   $$\(\widetilde W_1(y)+y-\half\)\(\widetilde W_1(-y)-y-\half\)
   =(\half-y)(\half+y),$$
   by Lemma~\ref{W_1 expansion}. Hence, $\Omega$
   is admissible by Remark~\ref{W_1 remark}.
\end{proof}

\begin{Defn}\label{rational}
The parameter set $\Omega$ is \textsf{rational} if there exists a
$k>0$ and $a_1,\dots,a_k\in R$ such that $\Omega$ satisfies the linear
recursion
$$\omega_{i+k}+a_1\omega_{i+k-1}+\dots+a_k\omega_i=0,$$
for all $i\gg0$.
\end{Defn}

Equivalently, $\Omega$ is rational if $\Omega$ is admissible and
$\tilde W_1(y)$ is a rational function. See Lemma~\ref{rational E}
for another characterization of rationality.

Rationality allows us to give a partial converse to
Corollary~\ref{admiss}.

\begin{Prop}\label{rational irreds}
Suppose that $R$ is an algebraically closed field and that $\Omega$ is
rational. Then every finite dimensional irreducible
$\Waff_n(\Omega)$--module can be considered as an irreducible module
for some cyclotomic Nazarov--Wenzl algebra $\W_{r,n}(\bu)$ with
$\Omega$ being $\bu$-admissible.

In particular, if $\Omega$ is rational then $\Omega$ is
$\bu$--admissible for some~$\bu$.
\end{Prop}

\begin{proof}
As $\Omega$ is rational, $\tilde W_1(y)$ is a rational function and we may 
write
$$
\frac{\tilde W_1(y)+y-\frac1{2}}{y+\frac1{2}}
  =\frac{\prod_i\ (y-\alpha_i)^{n_i}}{\prod_j\ (y-\beta_j)^{m_j}},
$$
for some non--negative integers $n_i$ and $m_j$ and with the
$\alpha_i, \beta_j\in R$ being pairwise distinct. Using
Remark~\ref{W_1 remark} it follows easily that
$$
\tilde W_1(y)+y-\frac1{2}=(y+\frac1{2})\prod_{i=1}^s \frac{y+c_i}{y-c_i},
$$
for some $c_i\in R$ and some $s\ge0$.

Now suppose that $M$ is a finite dimensional irreducible
$\Waff_n(\Omega)$--module and let
$(X_1-\lambda_1)\cdots(X_1-\lambda_d)$ be the characteristic
polynomial for the action of~$X_1$ on~$M$. Set
$$\bu=
   \begin{cases}
    (c_1,\dots,c_s,\lambda_1,\dots,\lambda_d,-\lambda_1,\dots,-\lambda_d),
          &\text{if $s$ is odd},\\ 
    (c_1,\dots,c_s,\lambda_1,\dots,\lambda_d,-\lambda_1,\dots,-\lambda_d,0),
          &\text{if $s$ is even}.
\end{cases}$$
Put $r=s+2d$ if $s$ is odd and $r=s+2d+1$ if $s$ is even.
Then $M$ is an irreducible $\W_{r,n}(\bu)$-module and
$\Omega$ is $\bu$-admissible.
\end{proof}

We will improve on this result by showing that we can construct all of
the irreducible modules for the affine Wenzl algebras in
Theorem~\ref{finite irreds} below.

\section{The seminormal representations of $\W_{r,n}(\bu)$}

\def\s{\mathfrak t}
\def\ts{\tilde\s}
\def\t{\mathfrak u}
\def\u{\mathfrak w}
\def\v{\mathfrak v}

In this section, we will give an explicit description of the
irreducible representations of $\W_{r,n}(\bu)$ in the special case
when $R$ is an field of characteristic greater than~$2n$ and when the
parameters $\bu$ satisfy some rather technical assumptions; see
Theorem~\ref{seminormal}.

The semisimple irreducible representations of the Brauer algebra
$\B_n(\omega)$ are labelled by partitions of $n-2m$, where $0\le
m\le\floor{n2}$, and a basis of the representation indexed by the
partition $\lambda$ is indexed by the set of updown
$\lambda$--tableaux. Analogously, we might expect that the semisimple
irreducible representations of $\W_{r,n}(\bu)$ should be indexed by
the multipartitions of $n-2m$, with the bases of these modules being
indexed by the updown $\lambda$--tableaux, where $\lambda$ is a
multipartition. We will see that this is the case. We begin by
defining these combinatorial objects.

Recall that a \textsf{partition} of $m$ is a sequence of weakly
decreasing non--negative integers $\tau=(\tau_1,\tau_2,\dots)$ such that
$|\tau|:=\tau_1+\tau_2+\cdots=m$. Similarly, an \textsf{$r$--multipartition}
of $m$, or more simply a multipartition, is an ordered $r$--tuple
$\lambda=(\lambda^{(1)},\dots,\lambda^{(r)})$ of partitions
$\lambda^{(s)}$, with $|\lambda|:=|\lambda^{(1)}|+\dots+|\lambda^{(r)}|=m$. If $\lambda$ is a multipartition of $m$ then we write $\lambda\vdash m$.

If $\lambda$ and $\mu$ are two multipartitions we say that $\mu$ is
obtained from $\lambda$ by \textsf{adding} a box if there exists a
pair $(i,s)$ such that $\mu^{(s)}_i=\lambda^{(s)}_i+1$ and
$\mu^{(t)}_j=\lambda^{(t)}_j$ for $(j,t)\ne(i,s)$. In this situation
we will also say that $\lambda$ is obtained from $\mu$ by
\textsf{removing} a box and we write $\lambda\subset\mu$ and
$\mu\setminus\lambda=(i,\lambda^{(s)}_i,s)$. We will also say that the
triple $(i,\lambda^{(s)}_i,s)$ is an \textsf{addable} node of $\lambda$
and a \textsf{removable} node of $\mu$. Note that $|\mu|=|\lambda|+1$.

Fix an integer $m$ with $0\le m\le\floor{n2}$ and let $\lambda$ be a
multipartition of $n-2m$.  An \textsf{$n$--updown
$\lambda$--tableau}, or more simply an updown $\lambda$--tableau, is a
sequence $\t=(\t_1,\t_2,\dots,\t_n)$ of multipartitions where
$\t_n=\lambda$ and the multipartition $\t_i$ is obtained from
$\t_{i-1}$ by either \textit{adding} or \textit{removing} a box, for
$i=1,\dots,n$, where we set $\t_0$ equal to the empty
multipartition $\emptyset$. Let $\UPD_n(\lambda)$ be the set of updown
$\lambda$--tableaux of $n$.  Note that $\lambda$ is a multipartition
of $n-2m$ and each element of $\UPD_n(\lambda)$ is an $n$--tuple of
multipartitions, so the $n$ is necessary in this notation.

In the special case when $\lambda$ is a multipartition of $n$ (so $m=0$),
there is a natural bijection between the set of $n$--updown
$\lambda$--tableaux and the set of standard $\lambda$--tableaux
in the sense of \cite{DJM:cyc}. This is the origin of the terminology of
updown $\lambda$--tableaux. If $\lambda$ is a multipartition of $n$ we set
$\Std(\lambda)=\UPD_n(\lambda)$ and refer to the elements of
$\Std(\lambda)$ as standard $\lambda$--tableaux.\label{tableaux}

\begin{Defn}\label{tilde}
    Suppose $1\le k\le n$. Define an equivalence relation $\simk$ on
    $\UPD_n(\lambda)$ by declaring that $\t\simk \s$ if $\t_j=\s_j$ 
    whenever $1\le j\le n$ and $j\neq k$, for $\s,\t\in\UPD_n(\lambda)$.
\end{Defn}

The following result is an immediate consequence of Definition~\ref{tilde}.

\begin{Lemma} Suppose $\s\in \UPD_n(\lambda)$ with
$\s_{k-1}=\s_{k+1}$. Then there is a bijection between the set of
all addable and removable nodes of $\s_{k-1}$ and the set of
$\t\in \UPD_n(\lambda)$ with $\t\simk \s$.
\end{Lemma}

Let $\lambda$ be a multipartition and suppose that $\t$ is an
$n$--updown $\lambda$--tableaux. For $k=2,\dots,n$ the mutipartitions
$\t_k$ and $\t_{k-1}$ differ by exactly one box; so either
$\t_k\subset\t_{k-1}$ or $\t_{k-1}\subset\t_k$. We define the
\textsf{content} of $k$ in $\t$ to be the scalar $c_\t(k)\in R$ given by
$$c_\t(k)=\begin{cases} \label{content}
     j-i+u_s,  &\text{if }\t_k{\setminus}\t_{k-1}=(i,j,s),\\ 
     i-j-u_s, &\text{if }\t_{k-1}{\setminus}\t_{k}=(i,j,s). 
\end{cases}$$ 
More generally, if $\alpha=(i,j,s)$ is an addable node of $\lambda$ we
define $c(\alpha)=u_s+j-i$ and if $\alpha$ is a removable node of
$\lambda$ we set $c(\alpha)=-(u_s+j-i)$.

The key property of contents that we need to construct the seminormal
representations is the following. Note that we are not (yet) assuming
that $R$ is a field.

\begin{Defn}\label{W-generic}
    The parameters $\bu=(u_1,\dots,u_r)$ are \textbf{generic} for
    $\W_{r,n}(\bu)$ if whenever there exists $d\in\Z$ such that either
    $u_i\pm u_j=d\cdot1_R$ and $i\ne j$,
    or $2u_i=d\cdot1_R$ then $|d|\ge2n$.
\end{Defn}

For example, $\bu$ is generic for $\W_{r,n}(\bu)$ if $u_1,\dots,u_r$
are algebraically independent over a subfield of~$R$.

\begin{Lemma}\label{generic u}
    Suppose that the parameters $\bu$ are generic for $\W_{r,n}(\bu)$
    and that $\Char R>2n$. Let $\lambda$ be a multipartition of
    $n-2m$, where $0\le m\le\floor{n2}$, and suppose that
    $\s,\t\in\UPD_n(\lambda)$. Then
\begin{enumerate}
\item $\s=\t$ if and only if $c_\s(k)=c_\t(k)$, for $k=1,\dots,n$;
\item if $1\le k<n$ then $c_\s(k)-c_\s(k+1)\ne0$; and,
\item if $\s_{k-1}=\s_{k+1}$ then $c_\s(k)\pm c_\t(k)\ne0$, whenever
    $\t\simk\s$ and $\t\ne\s$.
\item $2c_{\s}(k)\pm1\ne0$, for $1\le k\le n$. 
\end{enumerate}
\end{Lemma}

\begin{proof}Part (a) follows by induction on $n$. The key point is that our
assumptions imply that the contents of the addable and removable nodes in $\lambda$
are distinct so a $\lambda$--tableau $\s$ is uniquely determined by the sequence
of contents $c_\s(k)$, for $k=1,\dots,n$. The same argument proves parts (b),
(c) and (d).
\end{proof}

Until further notice we fix an integer $m$ with $0\le m\le \floor{n2}$ and
we fix a multipartition $\lambda$ of $n-2m$.

Motivated by \cite{Nazarov:brauer}, we introduce the following
rational functions in an indeterminate $y$. These functions will
play a key role in the construction of seminormal representations
of $\W_{r,n}(\bu)$.

\begin{Defn}\label{rational-W}
Suppose that $\s\in \UPD_n(\lambda)$.  For $1\le k\le n$, define
rational functions $W_k(y, \s)$ by
$$ W_k(y,\s)=\frac12-y+\big(y-\frac12(-1)^r\big)
         \prod_{\alpha}\frac{y+c(\alpha)}{y-c(\alpha)},$$
where $\alpha$ runs over the addable and removable nodes of the
multipartition $\s_{k-1}$.
\end{Defn}

The rational functions $W_k(y,\s)$ are related to the
combinatorics above by the following result. If $f(y)$ is a
rational function and $\alpha\in R$ then we write
$\Res_{y=\alpha}f(y)$ for the residue of $f(y)$ at $y=\alpha$.

\begin{Lemma}\label{w-residue}
    Suppose that $\bu$ is generic and $\Char R>2n$. Let
    $\s\in \UPD_n(\lambda)$ and $1\le k\le n$. Then
    $$ \frac{W_k(y, \s)}{y}=\sum_\alpha
         \Big(\Res_{y=c(\alpha)}\frac{W_k(y,\s)}{y}\Big)\cdot\frac 1{y-c(\alpha)},$$
    where $\alpha$ runs over the addable and removable nodes of
    $\s_{k-1}$.
\end{Lemma}

\begin{proof} As the $c(\alpha)$ are pairwise distinct, we can certainly write
$${W_k(y,\s)\over y}
    =a+\frac by +\sum_{ \alpha} \Big(\Res_{y=c(\alpha)}\frac{W_k(y,\s)}y\Big)
             \cdot {1\over y-c(\alpha)},$$
for some $a,b\in R$, where $\alpha$ runs over the addable and
removable nodes of $\s_{k-1}$.  Now, $a=\tfrac{W_k(y,
\s)}{y}\mid_{y=\infty}=0$. Let $c$ be the number of addable and
removable nodes of $\s_{k-1}$. Since a partition always has an odd
number of addable and removable nodes, we have that $(-1)^c=(-1)^r$.
Therefore,
$$b=\Res_{y=0}\frac{W_k(y,\s)}{y}=\frac12\big(1-(-1)^c(-1)^r\big)=0,$$
as we needed to show.
\end{proof}

We are now ready to define the matrices which make up the seminormal
form.

\begin{Defn}\label{e}
    Let $\lambda$ be a multipartition and $k$ an integer with
    $1\le k\le n$. Suppose that
    $\s$ and $\t$ are updown $\lambda$--tableaux in $\UPD_n(\lambda)$
    such that $\s_{k-1}=\s_{k+1}$. Then we define the scalars
    $e_{\s\t}(k)\in R$ by
$$\label{ew}
e_{\s\t}(k)=\begin{cases}\displaystyle
    \Res_{y=c_\s(k)} \frac{W_k(y,\s)}{y},&\text{if }\s=\t,\\[10pt]
    \sqrt{e_{\s\s}(k)}\sqrt{e_{\t\t}(k)},&
              \text{if }\s\ne\t \text{ and }\t\simk \s,\\[10pt]
    0,&\text{otherwise}.
\end{cases}
$$
(In (\ref{be}) below we will fix the choice of square roots
$\sqrt{e_{\s\s}(k)}$, for $\s\in\UPD_n(\lambda)$ and $1\le k\le n$.)
\end{Defn}

Note that when $c_\s(k)\ne0$ then
$\displaystyle e_{\s\s}(k)=\Res_{y=c_\s(k)}\frac{W_k(y,\s)+y-\frac12}{y}.$

We remark that if $\s_{k-1}\ne\s_{k+1}$ then the definition of
$e_{\s\s}(k)$ still makes sense, however, we do not define
$e_{\s\s}(k)$ in this case as we will not need it (see
Theorem~\ref{seminormal} below).

It follows from  Definition~\ref{rational-W} that
\begin{equation}\label{ek}
    e_{\s\s}(k)=\big(2c_\s(k)-(-1)^r\big)
             \prod_{\alpha}\frac{c_\s(k)+c(\alpha)}{c_\s(k)-c(\alpha)}
\end{equation}
where $\alpha$ runs over all addable and removable nodes of $\s_{k-1}$
with $c(\alpha)\neq c_\s(k)$.  Note that Lemma~\ref{generic u} now
implies that:
$$ \text{If $\t\simk\s$ then $e_{\s\t}(k)\ne0$, for $1\le k<n$.}$$
This will be used many times below without mention. Also observe that
Lemma~\ref{w-residue} can be restated as
\begin{equation}\label{w-res}
    \frac{W_k(y,\s)}{y}=\sum_{\t\simk\s}\frac{e_{\t\t}(k)}{y-c_\t(k)}.
\end{equation}

Given two partitions $\s$ and $\t$ write $\s\ominus\t=\alpha$ if
either $\t\subset\s$ and $\s\setminus\t=\alpha$, or $\s\subset\t$ and
$\t\setminus\s=\alpha$.

\begin{Defn}\label{ab}
    Let $\s\in\UPD_n(\lambda)$ and suppose that $\s_{k-1}\ne\s_{k+1}$, for some $k$
    with $1\le k<n$.
\begin{enumerate}
\item We define
$$
a_{\s}(k)=\frac1{c_\s(k+1)-c_\s(k)}\quad\text{ and }\quad
        b_{\s}(k)={\sqrt{1-a_{\s}(k)^2}}.
$$
(We fix the choice of square root for $b_\s(k)$ in (\ref{be})
below. Note that $c_\s(k+1)-c_\s(k)\ne0$ by Lemma~\ref{generic u}(b).)
\item If $\s_k\ominus\s_{k-1}$ and
$\s_{k+1}\ominus\s_k$ are in different rows and in different columns
then we define $S_k\s$ to be the updown $\lambda$--tableau
$$S_k\s=(\s_1, \cdots,
\s_{k-1}, \t_{k}, \s_{k+1}, \cdots, \s_{n})$$
where $\t_{k}$ is the multipartition which is uniquely determined by
the conditions $\t_k\ominus\s_{k+1}=\s_{k-1}\ominus \s_{k}$ and
$\s_{k-1}\ominus\t_{k}=\s_{k}\ominus \s_{k+1}$.
If the nodes $\s_k\ominus\s_{k-1}$ and $\s_{k+1}\ominus\s_k$ are both in the
same row, or both in the same column, then $S_k\s$ is not defined.
\end{enumerate}
\end{Defn}

We remark that if $\s_{k-1}=\s_{k+1}$ then the definitions of
$a_\s(k)$ and $b_\s(k)$ both make sense, however, we do not define
them in this case as we will never need them (see
Theorem~\ref{seminormal} below). Moreover, the condition
$\s_{k-1}\ne\s_{k+1}$ is crucial in proving Lemma~\ref{x-equal}(b)
below. (In fact, if we drop this condition then
Lemma~\ref{x-equal}(b) is not correct.)

We leave the following Lemma as an exercise to help the reader
familiarize themselves with the definitions.

\begin{Lemma}\label {x-equal}
    Suppose that $\s\in \UPD_n(\lambda)$ and $1\le k<n$. Then:
    \begin{enumerate}
    \item If $S_k\s$ is defined then $c_\s(k)=c_{S_k\s}(k+1)$
        and $c_\s(k+1)=c_{S_k\s}(k)$; consequently, $a_{S_k\s}(k)=-a_\s(k)$.
    \item If $S_k\s$ is not defined then $a_{\s}(k)=\pm1$ and $b_{\s}(k)=0$.
    \end{enumerate}
\end{Lemma}

Finally, if $\s_{k-1}=\s_{k+1}$ and $\t\simk\s$, where $1\le k<n$, we set
$$s_{\s\t}(k)=\frac{e_{\s\t}(k)-\delta_{\s\t}}{c_\s(k)+c_\t(k)}.$$
Note that $c_\s(k)+c_\t(k)\ne0$ by Lemma~\ref{generic u}.

We will assume that we have chosen the square roots in the definitions
of $b_\s(k)$ and $e_{\s\t}(k)$ so that the following equalities hold.

\begin{Assumption}[Root conditions]\label{be}
We assume that the ring $R$ is large enough so that
$\sqrt{e_{\s\s}(k)}\in R$ and $b_\s(k)=\sqrt{1-a_\s(k)^2}\in R$, for
all $\s,\t\in\UPD_n(\lambda)$ and $1\le k<n$, and that the following
equalities hold:
\begin{enumerate}
  \item If $\s_{k-1}\ne\s_{k+1}$ and $S_k\s$ is defined
        then $b_{S_k\s}(k)=b_\s(k)$.
  \item If $\s_{k-1}\ne\s_{k+1}$ and $\s\siml\t$, where $|k-l|>1$,
      then $b_\s(k)=b_\t(k)$.
  \item If $\s_{k-1}\ne\s_{k+1}$, $\s_k\ne\s_{k+2}$ and $S_k\s$
    and $S_{k+1}\s$ are both defined then $b_{S_{k+1}\s}(k)=b_{S_k\s}(k+1)$.
  \item If $\s_{k-1}=\s_{k+1}$ and $\s_k=\s_{k+2}$ then
      $\sqrt{e_{\s\s}(k)}\sqrt{e_{\s\s}(k+1)}=1$.
  \item If $\s_{k-1}=\s_{k+1}$, $\t_{k-1}=\t_{k+1}$ and
      $e_{\s\s}(k)=e_{\t\t}(k)$ then
      $\sqrt{e_{\s\s}(k)}=\sqrt{e_{\t\t}(k)}$.
  \item If $\s_{k-1}=\s_{k+1}$, $\s_k=\s_{k+2}$ and
      $\t\simkk\s$, $\u\simk\s$ with $S_k\t$ and $S_{k+1}\u$ both
      defined and $S_k\t=S_{k+1}\u$ then
      $b_\t(k)\sqrt{e_{\t\t}(k+1)}=b_\u(k+1)\sqrt{e_{\u\u}(k)}$.
\end{enumerate}
\end{Assumption}

In Lemma~\ref{be real} below we show that if $R=\R$ then it is
possible to choose $\bu$ so that  the Root Condition is satisfied.

Assuming (\ref{be}) we can now give the formulas for the seminormal
representations of $\W_{r,n}(\bu)$.

\begin{Theorem} \label{seminormal}
Suppose that $R$ is a field such that $\Char R>2n$ and that the
root conditions $($\ref{be}$)$ hold in~$R$.  Assume that $\bu$ is
generic for $\W_{r,n}(\bu)$.  Let $\Delta(\lambda)$ be the~$R$--vector
space with basis $\set{v_\s|\s\in\UPD_n(\lambda)}$. Then
$\Delta(\lambda)$ becomes a $\W_{r,n}(\bu)$--module via
\begin{align*}
\bullet\quad
   S_k v_{\s}&=\begin{cases}
         \displaystyle \sum_{\t\simk\s} s_{\s\t}(k)v_\t,&
               \text{if }\s_{k-1}= \s_{k+1},\\[15pt]
     a_{\s}(k) v_{\s}+b_{\s}(k)v_{S_k\s}, &\text{if }\s_{k-1}\neq \s_{k+1},
      \end{cases}\\
\bullet\quad
    E_k v_{\s}&=\begin{cases}
        \displaystyle\sum_{\t\simk\s} e_{\s\t}(k)v_\t,&\text{if }\s_{k-1}=\s_{k+1}\\[15pt]
    0, & \text{if }\s_{k-1}\neq \s_{k+1},\\
      \end{cases}\\
\bullet\quad
     X_j v_\s&=c_\s(j)v_\s,
 \end{align*}
 for $1\le k<n$ and $1\le j\le n$ and where we set $v_{S_k\s}=0$ if
 $S_k\s$ is not defined.
\end{Theorem}

\begin{Defn}\label{seminormal rep}
We call $\Delta(\lambda)$ a \textsf{ seminormal representation}
of $\W_{r,n}(\bu)$.
\end{Defn}

We note that the action of the operators $E_k$ and $S_k$ on
$\Delta(\lambda)$, with respect to the basis
$\set{v_\s|\s\in\UPD_n(\lambda)}$, are given by symmetric matrices,
for $0\le k<n$.

For the remainder of this section we assume that $R$ is an
algebraically closed field with $\Char R>2n$ and that the
parameters $\bu$ are generic for $\W_{r,n}(\bu)$ and satisfy
(\ref{be}). The proof of Theorem~\ref{seminormal} will occupy the
rest of this section. Our strategy is to use the rational functions
$W_k(\s,k)$ to verify that the action that we have just defined of
$\W_{r,n}(\bu)$ on $\Delta(\lambda)$ respects all of the relations in
$\W_{r,n}(\bu)$.

Throughout this section it will be convenient to work with formal
(infinite) linear combinations of elements of $\Delta(\lambda)$ and
$\W_{r,n}(\bu)$; alternatively, the reader may prefer to think that we
have extended our coefficient ring from $R$ to $R((y^{-1}))$, where
$y$ is an indeterminate over $R$. In fact, at times we will need to
work with formal series involving more than one indeterminate.

If $A$ is an algebra we let $Z(A)$ be its center.

\begin{Lemma}\label{tilde W}
    Suppose $k\ge0$ and that $a\ge0$. Then there exist elements
    $\omega^{(a)}_k$ in $Z\(\W_{r,k-1}(\bu)\)\cap R[X_1,\dots,X_{k-1}]$
    such that
    $$E_kX_k^aE_k=\omega^{(a)}_kE_k,$$
    and the degree of $\omega_k^{(a)}$, as a polynomial in
    $X_1,\dots,X_{k-1}$, is less than or equal to~$a$.
    Moreover, the generating series
    $\widetilde W_k(y)=\sum_{a\ge0}\omega_k^{(a)}y^{-a}$ satisfies
    $$\widetilde W_{k+1}(y)=-y+\half+
     \frac{(y+X_k)^2-1}{(y-X_k)^2-1}\frac{(y-X_k)^2}{(y+X_k)^2}
     \(\widetilde W_k(y)+y-\half\).$$
\end{Lemma}

\begin{proof}Observe that
    $\sum_{a\ge0}E_kX_k^aE_ky^{-a}=E_k\frac{y}{y-X_k}E_k$, so to prove
    the Lemma it is enough to argue by induction on $k$ to show that
    $E_k\frac1{y-X_k}E_k=\frac1y\widetilde W_k(y)E_k$, where
    $\widetilde W_k(y)$ and its coefficients are as above.

    If $k=1$ then there is nothing to prove. Assume that $k>1$.
    Starting with the identity
    $$S_k\frac1{y-X_k}=\frac1{y-X_{k+1}}S_k+\frac1{y+X_k}E_k\frac1{y-X_k}
                         -\frac1{(y-X_k)(y-X_{k+1})}$$
    Nazarov~\cite[Prop.~4.2]{Nazarov:brauer} proves that
    $E_{k+1}\frac1{y-X_{k+1}}E_{k+1}=\frac1y\widetilde W_{k+1}(y)E_{k+1}$,
    where $\widetilde W_{k+1}(y)$ satisfies the recurrence relation
    above. Nazarov assumes that he is working over the complex field
    (so, $R=\C)$, however, his arguments are valid over an arbitrary
    ring. Nazarov also proves that if $R=\C$ then the coefficients of
    $\widetilde W_k(y)$ are central in $\W_{r,k-1}(\bu)$. We
    modify Nazarov's arguments to establish centrality for
    fields of positive characteristic.

    By induction we may assume that the coefficients of $\widetilde
    W_k(y)$ commute with $E_1,\dots,E_{k-2}$ and $S_1,\dots,S_{k-2}$,
    so it is enough to show that the coefficients of $\widetilde
    W_{k+1}(y)$  commute with $E_{k-1}$ and $S_{k-1}$.
    Since $k\ge2$ we can write
     $$\frac{\widetilde W_{k+1}(y)+y-\half}{\widetilde W_{k-1}(y)+y-\half}
     =\frac{\mathcal X}{\mathcal Y}
     :=\frac{(y+X_k)^2-1}{(y-X_k)^2-1}\frac{(y-X_k)^2}{(y+X_k)^2}
       \frac{(y+X_{k-1})^2-1}{(y-X_{k-1})^2-1}\frac{(y-X_{k-1})^2}{(y+X_{k-1})^2}
     $$
     As $E_{k-1}$ and $S_{k-1}$ commute with $\W_{r,k-2}(\bu)$ it is enough
     to show that
     $E_{k-1}\frac{\mathcal X}{\mathcal Y}=\frac{\mathcal X}{\mathcal Y}E_{k-1}$
     and $S_{k-1}\frac{\mathcal X}{\mathcal Y}=\frac{\mathcal X}{\mathcal Y}S_{k-1}$.
     Now, $E_{k-1}\frac{\mathcal X}{\mathcal Y}=\frac{\mathcal X}{\mathcal Y}E_{k-1}$
     if and only if ${\mathcal Y}E_{k-1}\mathcal X=\mathcal XE_{k-1}\mathcal Y$, and
     this follows easily using relation \ref{Waff relations}(i).

     To prove that $S_{k-1}$ commutes with
     $\frac{\mathcal X}{\mathcal Y}$ let
$$\sum_{m\ge0}a_mz^m=\frac{(1+X_{k-1}z)(1+X_kz)}{(1-X_{k-1}z)(1-X_kz)},$$
     where $z=-y^{-1}$ or $z=(y\pm1)^{-1}$.
     Then $a_0=1$, $a_1=2(X_{k-1}+X_k)$, $a_2=2(X_{k-1}+X_k)^2$ and
     $$a_m=(X_{k-1}+X_k)a_{m-1}-X_{k-1}X_ka_{m-2}, \quad\text{for } m\ge3.$$
     Consequently, if $m\ge1$ then $a_m=(X_{k-1}+X_k)f_m(X_{k-1},X_k)$, for some
     $f_m\in R[X_{k-1},X_k]$. Now, relation \ref{Waff relations}(e)
     implies that $S_{k-1}$ and $X_{k-1}+X_k$ commute. Therefore, by induction,
     \begin{align*}
         S_{k-1}a_m
         &=(X_{k-1}+X_k)S_{k-1}a_{m-1}-(X_{k-1}X_kS_{k-1}+E_{k-1}X_k-X_kE_{k-1})a_{m-2}\\
           &=(X_{k-1}+X_k)a_{m-1}S_{k-1}-X_{k-1}X_ka_{m-2}S_{k-1}\\
           &=a_mS_{k-1}
       \end{align*}
       as required.

     Finally, it follows from the recurrence relation that
     $\omega_k^{(a)}\in R[X_1,\dots,X_{k-1}]$ and that
     $\omega_k^{(a)}$ has total degree at most $a$ as a polynomial in
     $X_1,\dots,X_{k-1}$.
\end{proof}

\begin{Remark}To prove that the $\omega_k^{(a)}\in Z\(\W_{r,k-1}(\bu)\)$ Nazarov
uses the identity
    $$\exp\Bigg(\sum_{a\ge0}2\big(X_{k-1}^{2a+1}+X_k^{2a+1}\big)\frac{y^{-2a-1}}{2a+1}\Bigg)
                 =\frac{(y+X_{k-1})(y+X_k)}{(y-X_{k-1})(y-X_k)}.$$
                 However, this formula is only valid in characteristic zero.
\end{Remark}

By Lemma~\ref{tilde W}, we have
$$
\widetilde W_k(y)+y-\half
      =\Big(\widetilde W_1(y)+y-\half\Big)
      \prod_{i=1}^{k-1}\frac{(y+X_i)^2-1}{(y-X_i)^2-1}
      \cdot\frac{(y-X_i)^2}{(y+X_i)^2}.
$$
As the right hand side acts on each $v_{\s}$ as multiplication by
a scalar we can define $\widetilde W_k(y,\s)\in R((y^{-1}))$ by
$\widetilde W_k(y)v_\s=\widetilde W_k(y,\s)v_\s$.

The next Proposition gives a representation theoretic interpretation
of the rational functions $W_k(y,\s)$ which were introduced in
Definition~\ref{rational-W}. 

\begin{Prop}\label{W=tilde W}
    Suppose that $\s\in \UPD_n(\lambda)$ and that $1\le k\le n$.
Then $$W_k(y,\s)=\widetilde W_k(y,\s).$$
\end{Prop}

\begin{proof}
As $\Omega$ is $\bu$-admissible, by Lemma~\ref{W_1 expansion} we have
$$
\widetilde W_1(y, \s)+y-\frac12
        =\big(y+\frac12(-1)^{r+1}\big)\prod_{t=1}^r {y+u_t\over y-u_t}.$$
Consequently, we can rewrite the definition of $\widetilde W_k(y,\s)$ as
$$
    \widetilde W_k(y, \s)+y-\frac12
       = (y-\frac12(-1)^r)\cdot\prod_{t=1}^r \frac{(y+u_t)}{(y-u_t)}
         \prod_{i=1}^{k-1} \frac{(y+c_\s (i))^2-1}{(y-c_\s (i))^2-1}
         \cdot\frac{(y-c_\s (i))^2}{(y+c_\s (i))^2}.
$$
If $c_\s(i)=-c_\s(j)$, for some $1\le i, j\le k-1$ with $i\neq j$, then
$$ \frac{(y+c_\s (i))^2-1}{(y-c_\s (i))^2-1}
          \cdot\frac{(y-c_\s (i))^2}{(y+c_\s (i))^2}
          \cdot\frac{(y+c_\s (j))^2-1}{(y-c_\s (j))^2-1}
          \cdot\frac{(y-c_\s (j))^2}{(y+c_\s (j))^2}
        =1.$$
Hence, in computing $\widetilde W_k(y, \s)$
we can assume that
$\s=(\s_1,\dots,\s_m,\dots,\s_{k-1},\dots,\s_n)$
where $m=|\s_{k-1}|$,
$\s_m=\s_{k-1}$ and $c_\s(i)+c_\s(i+1)=0$ for $m<i<k-1$ with $i-m$ odd
(so $\s_{i+1}$ is obtained by adding a box to
$\s_i$, for $1\le i<m$,
and $\s_i=\s_{k-1}$ for $m\le i\le k-1$ with
$i-m$ even). Let
$\s_{k-1}=(\mu^{(1)},\mu^{(2)},\ldots,\mu^{(r)})$.
Fix $t$ with $1\le t\le r$ and, abusing notation, write
$\beta\in\mu_k^{(t)}$ to indicate
that $\beta=(k,j,t)$ is a node in row $k$ of $\mu^{(t)}$.
Let $p_1=(k,1,t)$, $p_2=(k,\mu_k^{(t)},t)$,
$p_3=(k,\mu_k^{(t)}+1,t)$ and $p_4=(k+1,1,t)$.
Then
$$\begin{array}{l}
\displaystyle\prod_{\beta\in\mu_k^{(t)}}
  \frac{(y+c'(\beta))^2-1}{(y-c'(\beta))^2-1}
  \cdot\frac{(y-c'(\beta))^2}{(y+c'(\beta))^2}\\
\displaystyle\qquad
  =\prod_{\beta\in\mu_k^{(t)}}\frac{y-c'(\beta)}{y-(c'(\beta)+1)}
    \frac{y-c'(\beta)}{y-(c'(\beta)-1)}
    \frac{y+(c'(\beta)+1)}{y+c'(\beta)}
    \frac{y+(c'(\beta)-1)}{y+c'(\beta)}\\
\displaystyle\qquad
  =\frac{y-c'(p_1)}{y-c'(p_3)}
    \frac{y-c'(p_2)}{y-c'(p_4)}
    \frac{y+c'(p_3)}{y+c'(p_1)}
    \frac{y+c'(p_4)}{y+c'(p_2)}\\
\displaystyle\qquad
  =\frac{y-c'(p_1)}{y+c'(p_1)}
    \frac{y-c'(p_2)}{y+c'(p_2)}
    \frac{y+c'(p_3)}{y-c'(p_3)}
    \frac{y+c'(p_4)}{y-c'(p_4)}\,,
\end{array}$$
where for $\beta=(a,b,t)$ we write
$c'(\beta)=b-a+u_t$.
Taking the product over all~$k$ shows that
$$
\frac{(y+u_t)}{(y-u_t)}\prod_{\beta\in\mu^{(t)}}
\frac{(y+c(\beta))^2-1}{(y-c(\beta))^2-1}
\cdot\frac{(y-c(\beta))^2}{(y+c(\beta))^2}
=\prod_\alpha\frac{y+c(\alpha)}{y-c(\alpha)},
$$
where, in the first product, every node is considered to be an addable
node and, in the second product, $\alpha$ runs over the addable and
removable nodes of $\mu^{(t)}$. Hence,
$$ \widetilde W_k(y,\s)+y-\frac12=\big(y-\frac12(-1)^r\big)
         \prod_{\alpha}\frac{y+c(\alpha)}{y-c(\alpha)},$$
where $\alpha$ runs over the addable and removable nodes
of~$\s_{k-1}=(\mu^{(1)},\ldots,\mu^{(r)})$.
\end{proof}

\begin{Cor}\label{exe=oe}
    Suppose that $\s\in \UPD_n(\lambda)$ and that $1\le k<n$ and $a\ge0$. Then
    $E_kX_k^aE_kv_\s=\omega_k^{(a)} E_k v_\s$.
\end{Cor}
\begin{proof}
    If $\s_{k-1}\ne\s_{k+1}$ then
$E_k X_k^i E_k v_{\s}=0=\omega_k^{(i)} E_k v_{\s}$, so we may assume that
$\s_{k-1}=\s_{k+1}$.
 Now, by definition,
 $e_{\s\t}(k)=\sqrt{e_{\s\s}(k)}\sqrt{e_{\t\t}(k)}$. So
\begin{align*}
    E_k {y\over y-X_k} E_k v_{\s}
      &=E_k \sum_{\t\overset k \sim \s}\frac{y}{y-c_\t(k)}e_{\s\t}(k) v_{\t}
         =\sum_{\u\simk \t} \sum_{\t\overset k \sim \s}
                  \frac{y}{y-c_\t(k)}e_{\t\u}(k) e_{\s\t}(k) v_{\u}\\
      &= \sum_{\u\simk \s}\Bigg(\sum_{\t\overset k \sim \s}
           \frac{y}{ y-c_\t(k) }e_{\t\t}(k)\Bigg)e_{\s\u}(k) v_{\u}\\
      &=W_k(y, \s) E_k v_{\s} =\widetilde W_k(y, \s) E_kv_{\s},
\end{align*}
by Proposition~\ref{W=tilde W}. By Lemma~\ref{tilde W},
$\omega_k^{(a)}\in R[X_1,\dots,X_{k-1}]$, so
$\omega_k^{(a)}v_\s=\omega_k^{(a)}v_\t$ whenever $\s\simk\t$.
Therefore,
\begin{align*}
    E_k {y\over y-X_k} E_k v_{\s}
    &=\sum_{\t\simk\s}e_{\s\t}(k)\widetilde W_k(y, \s)v_\t
     =\sum_{\t\simk\s}e_{\s\t}(k)\widetilde W_k(y, \t)v_\t\\
    &=\widetilde W_k(y)\sum_{\t\simk\s}e_{\s\t}(k)v_\t
     =\widetilde W_k(y)E_k v_{\s}.
\end{align*}
Comparing the coefficient of $y^{-a}$, for $a\ge0$, on both sides of the last
equation proves the Corollary.
\end{proof}

\begin{Lemma}\label{ees} Suppose that $\s\in \UPD_n(\lambda)$ with
$\s_{k-1}=\s_{k+1}$ and $\s_k=\s_{k+2}$. Then
$e_{\s\s}(k)e_{\s\s}(k+1)=1$.
\end{Lemma}

\begin{proof} The recursion formula of
Lemma~\ref{tilde W} and Proposition \ref{W=tilde W} show that
$$W_{k+1}(y,\s)+y-\frac12=
\big(W_k(y,\s)+y-\frac12\big)
{(y-c_\s(k))^2\over (y+c_\s(k))^2}{(y+c_\s(k))^2-1\over (y-c_\s(k))^2-1},$$
and, by definition,
$$W_k(y,\s)+y-\frac12=\big(y-\frac12(-1)^r \big)
\prod_{\t\simk\s} {y+c_\t (k)\over y-c_\t(k)}.$$ Thus,
$$\begin{aligned}
\frac{W_{k+1}(y, \s)+y-\half}y = & \Big(1-{1\over 2y}(-1)^r \Big)
{y-c_\s(k)\over y+c_\s(k)}{(y+c_\s(k))^2-1\over (y-c_\s(k))^2-1} \\
&\qquad\times \prod_{\t\simk\s, \t\neq \s} {y+c_\t (k)\over
y-c_\t(k)}.\\ \end{aligned} $$ Taking  residues at
$y=-c_\s(k)=c_\s(k+1)$ on both sides of this equation, we have
$$\begin{aligned} e_{\s\s}(k+1) &={2c_\s(k)+(-1)^r\over 4c_\s(k)^2-1}
\prod_{\t\simk\s, \t\ne\s} {c_\s(k)-c_\t(k)\over c_\s(k)+c_\t(k)}\\ &=
{1\over 2c_\s(k)-(-1)^r}\prod_{\t\simk \s, \t\ne\s}
{c_\s(k)-c_\t(k)\over c_\s(k)+c_\t(k)}={1\over e_{\s\s}(k)}.\\
\end{aligned} $$ where the last equality uses (\ref{ek}).
\end{proof}

We remark that the condition $\s_k=\s_{k+2}$ is needed in
Lemma~\ref{ees} only because $e_{\s\s}(k+1)$ is not defined without
this assumption.

\begin{Lemma}\label{be equality}
Fix an integer $k$ with $1\le k<n-1$ and suppose that
$\s,\t,\u\in\UPD_n(\lambda)$ are updown $\lambda$--tableaux such that
$\s_{k-1}=\s_{k+1}$, $\s_k=\s_{k+2}$, $\t\simkk\s$, $\u\simk\s$ and
that $S_k\t$ and $S_{k+1}\u$ are both defined with
$S_k\t=S_{k+1}\u$. Then $b_\t(k)^2e_{\t\t}(k+1)=b_\u(k+1)^2e_{\u\u}(k)$.
\end{Lemma}

\begin{proof}
Let $\sigma=\s_k\ominus\s_{k-1}$ and
$\tau=\t_{k+1}\ominus\s_k$. $S_k\t=S_{k+1}\u$ implies
$\tau=\u_k\ominus\s_{k-1}$. Then, by (\ref{ek}),
\begin{align*}
e_{\u\u}(k)
  &=(2c(\tau)-(-1)^r)
      \prod_\alpha\frac{c(\tau)+c(\alpha)}{c(\tau)-c(\alpha)},\\
\intertext{where $\alpha$ runs over the addable and removable nodes
  of $\s_{k-1}=\u_{k-1}$ with $c(\alpha)\ne c(\tau)$ and, similarly,}
e_{\t\t}(k+1)
  &=(2c(\tau)-(-1)^r)
      \prod_\alpha\frac{c(\tau)+c(\alpha)}{c(\tau)-c(\alpha)},
\end{align*}
where $\alpha$ runs over all addable and removable nodes of
$\s_k=\t_k$ with $c(\alpha)\ne c(\tau)$.  We have
$e_{\u\u}(k)=\Res_{y=c(\tau)} \frac{W_k(y,\s)+y-\frac12}{y}$
and
$e_{\t\t}(k+1)=\Res_{y=c(\tau)} \frac{W_{k+1}(y,\s)+y-\frac12}{y}$.
Further, by Lemma~\ref{tilde W} and Proposition \ref{W=tilde W}, we have
$$W_{k+1}(y,\s)+y-\frac12=\big(W_k(y,\s)+y-\frac12\big)
{(y+c(\sigma))^2-1\over (y-c(\sigma))^2-1}{(y-c(\sigma))^2\over
(y+c(\sigma))^2}.$$
It follows that
$$
\frac{e_{\t\t}(k+1)}{e_{\u\u}(k)}
   =\frac{(c(\sigma)+c(\tau))^2-1}{(c(\sigma)+c(\tau))^2}
      \frac{(c(\tau)-c(\sigma))^2}{(c(\tau)-c(\sigma))^2-1}
   =\frac{b_{\u}(k+1)^2}{b_{\t}(k)^2},$$
where the last equality follows from the definitions because
$(c_{\t}(k),c_{\t}(k+1),c_{\t}(k+2))=(c(\sigma),c(\tau),-c(\tau))$
and
$(c_{\u}(k),c_{\u}(k+1),c_{\u}(k+2))=(c(\tau),-c(\tau),c(\sigma))$.
\end{proof}

The following combinatorial identities will be used in the proof of
Theorem~\ref{seminormal}.

\begin{Prop}\label{identity}  Suppose that $\s, \t'\in \UPD_n(\lambda)$ with
$\s_{k-1}=\s_{k+1}$, $\s_k\neq \s_{k+2}$,  $ \t'\simk \s$
and $\t'\neq \s$, where $1\le k<n-1$. Let $\ts\in \UPD_n(\lambda)$ be the
updown tableau which is uniquely determined by the conditions
$\ts\simk \s$ and $\ts_k=\s_{k+2}$. Then the following
identities hold:
\begin{enumerate}
    \item$\displaystyle\sum_{\t\simk\s} {e_{\t\t}(k)\over
c_\s(k)+c_\t(k)} =1+{1\over 2c_\s(k)},$
\item
$\displaystyle\sum_{\t\simk\s}\frac{e_{\t\t}(k)}{(c_\s(k)+c_\t(k))^2}
         =\Big(1-\frac1{4c_\s(k)^2}\Big)\frac1{e_{\s\s}(k)} +
               \frac1{2c_\s(k)^2}$
\item
$\displaystyle\sum_{\t\simk\s}
  {e_{\t\t}(k)\over (c_\s(k)+c_\t(k)) (c_\t(k)+c_{\t'}(k))}={1\over
2c_\s(k)c_{\t'}(k)},$
\end{enumerate}
\end{Prop}

\begin{proof}
It follows from (\ref{w-res}) and Definition~\ref{ew} that
$$ \frac{W_k(y,\s)}y=\sum_{\t\simk\s}\frac{e_{\t\t}(k)}{y-c_\t(k)}.$$
Evaluating both sides at $y=-c_\s(k)$ and using (\ref{rational-W}) gives
(a).

By Proposition~\ref{W=tilde W} and Corollary~\ref{exe=oe} we have
$$ E_k {1\over (y-X_k) (v-X_k)}E_kv_{\s}
       ={1\over v-y}\Big( {W_k(y,\s)\over y}-{W_k(v,\s)\over v}\Big) E_k
v_{\s}.$$
Comparing the coefficients of $v_{\s}$ on both sides of this equation we
obtain
$$
  \sum_{\t\simk \s}\frac{e_{\t\t}(k)}{(y-c_\t(k))(v-c_\t(k))}
    =\frac1{v-y}\Big\{ \frac{W_k(y,\s)}{y}-\frac{W_k(v,\s)}{v} \Big\}.
$$
Setting $y=-c_\s(k)$ we obtain
\begin{align*}
\sum_{\t\simk \s}\frac{e_{\t\t}(k)}{(c_\s(k)+c_\t(k))(v-c_\t(k))}
   &=\frac1{v+c_\s(k)}\Bigg\{\Big(\frac{W_k(v,\s)}{v}+1-\frac1{2v}\Big)
   +\Big(\frac1{2c_\s(k)}+\frac1{2v}\Big)\Bigg\}\\
&=\frac{2v-(-1)^r}{2v(v+c_\s(k))}\prod_{\t\simk\s}\frac{v+c_\t(k)}{v-c_\t(k)
}
            +\frac{1}{2c_\s(k)v}.
\end{align*}
Setting $v=-c_{\t'}(k)$ gives (c). Now we set $v=-c_\s(k)$.
Then it gives
$$\sum_{\t\simk\s}\frac{e_{\t\t}(k)}{(c_\s(k)+c_\t(k))^2}
         =\frac{2c_\s(k)+(-1)^r}{4c_\s(k)^2}
            \prod_{\substack{\t\simk\s\\\t\neq\s}}
    \frac{ c_\s(k)-c_\t(k)}{c_\s(k)+c_\t(k)}+\frac1{2c_\s(k)^2}.$$
On the other hand, multiplying the reciprocal of (\ref{ek}) by
$\(1-\frac1{4c_\s(k)^2}\)$ gives
$$
\Big(1-\frac1{4c_\s(k)^2}\Big)\frac1{e_{\s\s}(k)}
=\frac{2c_\s(k)+(-1)^r}{4c_\s(k)^2}
 \prod_{\t\simk \s,\t\neq\s}\frac{c_\s(k)-c_\t(k)}{c_\s(k)+c_\t(k)}.
$$
Combining these two equations gives (b).
\end{proof}

We are now ready to start checking that the action of $\W_{r,n}(\bu)$
on $\Delta(\lambda)$ respects the relations of $\W_{r,n}(\bu)$.
We break the proof into several lemmas and propositions.

\begin{Lemma}\label{start relations}
Suppose  $\s\in\UPD_n(\lambda)$. Then
\begin{enumerate}
    \item $E_i^2v_\s=\omega_0 E_iv_\s$, for $1\le i<n$.
    \item $E_1X_1^aE_1v_\s=\omega_aE_1v_\s$, for $a>0$.
    \item $(X_1-u_1)(X_1-u_2)\cdots (X_1-u_r)v_{\s}=0$.
    \item $X_iX_j v_{\s}=X_jX_i v_{\s}$ for $1\le i, j\le n$.
    \item $E_i(X_i+X_{i+1}) v_{\s}=(X_i+X_{i+1}) E_i v_{\s}=0$,
           $1\le i\le n-1$.
    \item $(S_iX_i-X_{i+1}S_i)v_\s=(E_i-1)v_{\s}
               =(X_iS_i-S_{i}X_{i+1})v_{\s}$, for $1\le i\le n-1$
    \item $E_kE_l v_{\s}=E_lE_k v_{\s}$ if $|k-l|>1$.
    \item $E_k X_l v_{\s}=X_lE_k v_{\s}$ if $l\neq k, k+1$.
    \item $S_kX_lv_{\s}=X_lS_kv_{\s}$ if  $l\neq k, k+1$.
\end{enumerate}\end{Lemma}

\begin{proof} As $\omega_0=\omega_1^{(0)}$ and
$\omega_a=\omega_1^{(a)}$ by Lemma~\ref{tilde W}, parts (a)
    and (b) have already been proved in Corollary~\ref{exe=oe}.  Parts
    (c)--(f) follow directly from the definitions of the actions.  If
    $|k-l|>1$ then (\ref{be})(e) shows that
    (g) holds.  Assume now that $l\ne k,k+1$. If $\s_{k-1}\neq
    \s_{k+1}$ then $c_{S_k\s}(l)=c_\s(l)$.  If $\t\simk \s$ then
    $c_\t(l)=c_\s(l)$.  Combining the last two statements forces (h)
    and (i) to be true.
\end{proof}

\begin{Lemma} Suppose  $\s\in \UPD_n(\lambda)$. Then
$E_kE_{k\pm1} E_kv_{\s}=E_k v_{\s}$.
\end{Lemma}

\begin{proof} We only prove that $E_kE_{k+1} E_kv_{\s}=E_k v_{\s}$, since the argument for
    the case $E_kE_{k-1} E_kv_{\s}=E_k v_{\s}$ is almost identical.
  
    We may assume $\s_{k-1}=\s_{k+1}$ since, otherwise, $E_kE_{k+1}E_k
    v_{\s}=0=E_k v_{\s}$. Let $\ts$ be the unique $n$--updown
    tableau such that
    $\ts\simk\s$ and $\ts_k=\s_{k+2}$. We have
$$
E_kE_{k+1}E_k v_{\s}
  = e_{\s\ts}(k) e_{\ts\ts}(k+1) \sum_{\t\simk\ts}
        e_{\ts\t}(k) v_{\t}
  =e_{\ts\ts}(k)e_{\ts\ts}(k+1)
         \sum_{\t\simk\s}  e_{\s\t}(k) v_\t.
$$
Hence, $E_kE_{k+1} E_kv_{\s}=E_k v_{\s}$ by Lemma~\ref{ees}.
\end{proof}

It remains to check relations (a), (b)(i), (b)(ii), (d)(i) and~(g)
from Definition~\ref{Waff relations}.

\begin{Lemma}\label{s-square} Suppose that $\s\in \UPD_n(\lambda)$.
Then $S_k^2 v_{\s}=v_{\s}$.
\end{Lemma}

\begin{proof}
\Case{1. $\s_{k-1}\neq \s_{k+1}$}
If $S_k\s$ is not defined then $a_{\s}(k)\in\{-1, 1\}$ and
$b_{\s}(k)=0$, which implies $S_k^2 v_{\s}=v_{\s}$. If
$S_k\s\in\UPD_n(\lambda)$ then by the choice of the square roots in
(\ref{be})(a) we have
$$
S_k^2 v_{\s}=\Big(a_{\s}(k)^2+ b_{\s}(k)b_{S_k\s}(k)\Big) v_{\s}
               +\Big(a_{\s}(k)+a_{S_k\s}(k)\Big)b_{\s}(k) v_{S_k\s}
            =v_{\s}.
$$
\Case{2.  $\s_{k-1}=\s_{k+1}$}
We have
$S_k^2v_\s=\sum_{\t\simk\s}
        \Big(\sum_{\v\simk\s}s_{\s\v}(k)s_{\v\t}(k)\Big)v_\t$.
So, the coefficient of $v_\s$ in~$S_k^2\v_s$ is
$$\sum_{\t\simk \s} s_{\s\t}(k) s_{\t\s}(k)
     =\sum_{\t\simk\s}\frac{e_{\s\s}(k)e_{\t\t}(k)}{(c_\s(k)+c_\t(k))^2}
              -\frac{e_{\s\s}(k)}{2c_\s(k)^2}+\frac1{4c_\s(k)^2}=1,$$
where the last equality follows by rearranging
Proposition~\ref{identity}(b).  If $\t\simk\s$ and $\t\ne\s$  then
the coefficient of $v_\t$ in
$S_k^2v_\s$ is
\begin{align*}
\sum_{\v\simk \s} s_{\s\v}(k)s_{\v\t}(k)
   &=\sum_{\substack{\v\sim \s\\\s\ne\v\ne\t}}
      \frac{e_{\s\v}(k) e_{\v\t}(k)}%
           {(c_\s(k)+c_\v(k))(c_\v(k)+c_{\t}(k))}\\
    &\qquad
  +\frac{(e_{\s\s}(k)-1)e_{\s\t}(k)}{2c_\s (k) (c_\s(k)+c_{\t}(k))}
  +\frac{(e_{\t\t}(k)-1)e_{\s\t}(k)}{2c_{\t} (k) (c_\s(k)+c_{\t}(k))}\\
   &= e_{\s\t}(k)\Big( \sum_{\v\simk \s}
      \frac{e_{\v\v}(k)}{(c_\s(k)+c_\v(k))(c_\v(k)+c_{\t}(k))}
      -\frac1{2c_\s(k)c_{\t}(k)}\Big)\\
   &=0
\end{align*}
by Proposition~\ref{identity}(c). Therefore, $S_k^2 v_\s=v_{\s}$.
\end{proof}

The next two Propositions prove that the action of $\W_{r,n}(\bu)$ on
$\Delta(\lambda)$ respects the tangle relations \ref{Waff relations}(g).

\begin{Prop}\label{se}
 For any $\s\in \UPD_n(\lambda)$, $E_kS_k v_{\s}=E_k
v_{\s}=S_kE_k v_{\s}$.
\end{Prop}
\begin{proof}
Suppose that $\s_{k-1}\neq \s_{k+1}$. Then either $S_k\s$ is not
defined, or $(S_k\s)_{k-1}\ne(S_k\s)_{k+1}$. In either case, we have
$E_kS_k v_\s=E_k v_\s=S_kE_k v_\s=0$.
Suppose  $\s_{k-1}=\s_{k+1}$. Then
$$S_kE_kv_\s=\sum_{\t\simk\s}e_{\s\t}(k)S_kv_\t
          =\sum_{\t'\simk\t}\sum_{\t\simk\s}s_{\t\t'}(k)e_{\s\t}(k)v_{\t'}.$$
By Proposition~\ref{identity}(a), we have
$$\begin{aligned} \sum_{\t\overset{k}\sim \s}
    e_{\s\t}(k)s_{\t\t'}(k)& =\sum_{\t\simk \s, \t\neq
\t'} {e_{\s\t}(k) e_{\t\t'}(k) \over c_\t(k)+c_{\t'}(k)} +
e_{\s\t'}(k) {e_{\t'\t'}(k)-1\over 2 c_{\t'} (k)}\\
&=e_{\s\t'}(k)\Big(\sum_{\t\simk \s} {e_{\t\t}(k) \over
c_\t(k)+ c_{\t'}(k)}- \frac1{2c_{\t'}(k)}\Big)\\
& =e_{\s\t'}(k).\\
\end{aligned}
$$
Hence, $S_kE_k v_{\s}=E_k v_{\s}$. One can prove that
 $E_k S_k v_{\s}=E_k v_{\s}$ similarly.
\end{proof}

\begin{Prop}\label{see} Suppose that $\s\in \UPD_n(\lambda)$. Then
\begin{enumerate}
 \item $S_{k}E_{k+1}E_{k} v_{\s}=S_{k+1} E_kv_{\s}$.
 \item$E_{k+1}E_kS_{k+1}v_{\s} =E_{k+1} S_{k}v_{\s}$.
\end{enumerate}
\end{Prop}

\begin{proof}
(a) We may assume that $\s_{k-1}=\s_{k+1}$ since otherwise
$S_{k}E_{k+1}E_{k} v_{\s}=S_{k+1} E_kv_{\s}=0$. Let
$\ts\in\UPD_n(\lambda)$ be the unique updown tableau such that
$\ts\simk\s$ and $\ts_k=\s_{k+2}$. We have
\begin{align*}
S_kE_{k+1}E_k v_{\s}
   &= e_{\s\ts}(k)e_{\ts\ts}(k+1)\Big(
   s_{\ts\ts}(k)v_{\ts}
     +\sum_{\t\simk\ts, \t\ne\ts}s_{\ts\t}(k) v_{\t}\Big)\\
   &\qquad+ \sum_{\t\simkk\ts, \t\ne\ts}
     e_{\s\ts}(k)e_{\ts\t}(k+1)  \Big(a_{\t}(k) v_{\t} + b_{\t}(k)
           v_{S_k\t}\Big).\\
\intertext{Observe that if $S_k\t$ is defined, for $\t$ in the second sum,
then $(S_k\t')_k\ne \s_{k+2}$ and $\u=S_{k+1}S_k\t$ is also
defined. Further, we have $\u\simk\ts$ and $\u\ne\s$. Similarly,}
S_{k+1} E_kv_{\s}
    &= e_{\s\ts}(k)\Big(s_{\ts\ts}(k+1)v_{\ts}
          +\sum_{\t\simkk\ts, \t\ne\ts}
          s_{\ts\t}(k+1)v_{\t}\Big)\\
    &\qquad +\sum_{\t\simk\ts, \t\ne\ts} e_{\s\t}(k)
\Big(a_\t(k+1)v_{\t}+b_\t(k+1) v_{S_{k+1} \t}\Big).
\end{align*}
We now compare the coefficients of $v_{\t}$ in $S_kE_{k+1}E_k v_\s$
and $S_{k+1} E_k v_\s$. First, observe that
$e_{\ts\ts}(k)e_{\ts\ts}(k+1)=1$ by Lemma~\ref{ees}.

\Case{1. $\t=\ts$}
Since $c_{\ts}(k)=-c_{\ts}(k+1)$, the definitions and the
remarks above show that the coefficient of $v_\t$ in $S_kE_{k+1}E_k v_{\s}$
is equal to
$$
e_{\s\ts}(k) e_{\ts\ts}(k+1)s_{\ts\ts}(k)
  =e_{\s\ts}(k) {1-e_{\ts\ts}(k+1)\over 2 c_{\ts} (k)}
  =e_{\s\ts}(k)s_{\ts\ts}(k+1),
$$
which is the coefficient of $v_\t$ in $S_{k+1} E_kv_{\s}$.

\Case{2.  $\t\simk\ts$ and $\t\ne\ts$}
Now, $c_{\ts}(k)=c_\t(k+2)$ and $c_\t(k+1)=-c_\t(k)$, so the
coefficient of $v_\t$ in $S_kE_{k+1}E_k v_{\s}$ is
$$
e_{\s\ts}(k) e_{\ts\ts}(k+1)s_{\ts\t}(k)
  =\frac{e_{\s\t}(k)}{c_{\ts}(k)+c_\t(k)}
  =e_{\s\t}(k) a_\t(k+1),
$$
which is the coefficient of $v_\t$ in $S_{k+1} E_kv_{\s}$.

\Case{3. $\t\simkk\ts$ and $\t\ne\ts$}
Since $c_\t(k)=-c_{\ts}(k+1)$, the
coefficient of $v_\t$ in $S_kE_{k+1}E_k v_{\s}$ is
$$a_{\t}(k)e_{\ts\t}(k+1)e_{s\ts}(k)
=\frac{e_{\ts\t}(k+1)e_{s\ts}(k)}{c_\t(k+1)+c_{\ts}(k+1)}
  =e_{\s\ts}(k)s_{\ts\t}(k+1),$$
which is the coefficient of $v_\t$ in $S_{k+1} E_kv_{\s}$.

Now suppose that $S_k\t$ is defined and let $\u=S_{k+1}S_k\t$ be as
above. Then the coefficient of $v_{S_k\t}$ in $S_kE_{k+1}E_k v_{\s}$
is
\begin{align*}
e_{\s\ts}(k)e_{\ts\t}(k+1)b_\t(k)
  &=\sqrt{e_{\s\s}(k)}\sqrt{e_{\t\t}(k+1)}b_\t(k)\\
  &=\sqrt{e_{\s\s}(k)}\sqrt{e_{\u\u}(k)}b_\u(k+1)\\
  &=e_{\s\u}(k)b_\u(k+1),
\end{align*}
where the second equality comes from (\ref{be})(f). As
$S_k\t=S_{k+1}\u$ this is the coefficient of $v_{S_k\t}$ in
$S_{k+1} E_kv_{\s}$. This completes the proof of~(a).

\medskip\noindent(b) We let the reader work out the expansions of $E_{k+1}E_kS_{k+1}v_{\s}$
and $E_{k+1} S_{k}v_{\s}$. To show that these two expressions are
equal there are four cases to consider.

\Case{1. $\s_k=\s_{k+2}$ and $\s_{k-1}=\s_{k+1}$}
We have
$$\begin{aligned} E_{k+1}E_{k} S_{k+1} v_{\s}&=E_{k+1} e_{\s\s}(k) s_{\s\s}(k+1) v_{\s}
={1-e_{\s\s}(k)\over  2c_\s(k+1)}E_{k+1}v_{\s}\\              &
=s_{\s\s}(k) E_{k+1}v_{\s}=E_{k+1}S_{k} v_{\s}.
\end{aligned}
$$

\Case{2. $\s_k\neq \s_{k+2}$ and $\s_{k-1}=\s_{k+1}$}
Define $\ts\in\UPD_n(\lambda)$ to be the unique updown tableau
such that $\ts\simk\s$ and $\ts_k=\s_{k+2}$. Then
$\ts\ne\s$ and
$$E_{k+1}E_{k} S_{k+1} v_{\s}
  = a_{\s}(k+1)e_{\s\ts}(k)E_{k+1} v_{\ts}
  =s_{\s\ts}(k)E_{k+1} v_{\ts}=E_{k+1}S_{k}v_{s},
$$
where the second equality uses the facts that $c_\s(k+1)=-c_\s(k)$,
$c_\s(k+2)=c_{\ts}(k)$ and
$(S_{k+1}\s)_{k-1}\ne(S_{k+1}\s)_{k+1}$.

\Case{3.  $\s_k=\s_{k+2}$ and $\s_{k-1}\neq\s_{k+1}$}
Define $\ts\in\UPD_n(\lambda)$ to be the unique updown tableau
such that $\ts\simkk\s$ and $\ts_{k+1}=\s_{k-1}$.
Then
\begin{align*}
    E_{k+1}E_{k} S_{k+1} v_{\s}
       &= s_{\s\ts}(k+1) e_{\ts\ts}(k) E_{k+1}v_{\ts}\\
       &=\frac{e_{\s\ts}(k+1)e_{\ts\ts}(k)}{c_\s(k+1)+c_{\ts}(k+1)}
       \sum_{\t\simkk\ts} e_{\ts\t}(k+1)v_\t\\
       &=a_\s(k)\sum_{\ts\simkk\t}e_{\s\t}(k+1)v_\t
        = E_{k+1}S_kv_\s,
\end{align*}
where we have used the facts that $c_{\ts}(k+1)=-c_\s(k)$ and
$(S_k\s)_k\ne(S_k\s)_{k+2}$.

\Case{4.  $\s_k\neq \s_{k+2}$ and $\s_{k-1}\neq \s_{k+1}$}
First observe that because of our assumptions we have
$E_{k+1}E_kS_{k+1}v_\s=b_\s(k+1)E_{k+1}E_kv_{S_{k+1}\s}$ and
$E_{k+1}S_kv_\s=b_\s(k)E_{k+1}v_{S_k\s}$. If
$(S_{k+1}\s)_{k-1}\ne(S_{k+1}\s)_{k+1}$ then we also have
$(S_{k}\s)_{k}\ne(S_{k}\s)_{k+2}$ so that
$E_{k+1}S_{k}S_{k+1}v_{\s}=0=E_{k+1}S_{k}v_{\s}$.

Suppose now that $(S_{k+1}\s)_{k-1}=(S_{k+1}\s)_{k+1}$ and let
$\ts\in\UPD_n(\lambda)$ be the unique updown tableau such that
$\ts\simk S_{k+1}\s$ and $\ts_{k}=\s_{k+2}$. Set
$\t=S_k\s$ and $\u=S_{k+1}\s$ and observe that the assumptions of
(\ref{be})(f) hold, so that
$b_\t(k)\sqrt{e_{\t\t}(k+1)}=b_\u(k+1)\sqrt{e_{\u\u}(k)}$. As
$b_\s(k)=b_\t(k)$ and $b_\s(k+1)=b_\u(k+1)$, the reader should now have
no difficulty in using (\ref{be})(d), together with the fact that
$\t'\simkk\ts$ if and only if $\t'\simkk S_k\s$, to show that
\begin{align*}
E_{k+1}E_kS_{k+1}v_{\s}
  &=b_{\s}(k+1)\sum_{\t'\simkk\ts}
      e_{\ts,S_{k+1}\s}(k)e_{\ts,\t'}(k+1)v_{\t'}\\
  &=b_{\s}(k)\sum_{\t'\simkk S_k\s}e_{S_k\s,\t'}(k+1)v_{\t'}
   =E_{k+1}S_kv_{\s}.
\end{align*}
\end{proof}

The next Proposition shows that the action of $\W_{r,n}(\bu)$ on
$\Delta(\lambda)$ respects the two relations
\ref{Waff relations}(b)(i) and \ref{Waff relations}(d)(i).

\begin{Prop}
Suppose that $\s\in \UPD_n(\lambda)$ and that $|k-l|>1$. Then:
\begin{enumerate}
    \item $S_kS_lv_{\s}=S_lS_kv_{\s}$.
    \item $S_kE_lv_{\s}=E_lS_kv_{\s}$.
\end{enumerate}
\end{Prop}

\begin{proof}
We prove only (a) as the proof of part (b) is similar to, but easier
than~(a).

First suppose that $\s_{k-1}= \s_{k+1}$ and $\s_{l-1}=\s_{l+1}$. Then
$$ S_kS_lv_\s=\sum_{\t\siml\s,\,\u\simk\t} s_{\s\t}(l)s_{\t\u}(k)v_\u.$$
Now for each pair of updown tableaux $(\u,\t)$ with
$\u\simk\t\siml\s$ there is a unique updown tableau $\t'$ such that
$\u\siml\t'\simk\s$; more precisely, $\t'_k=\u_k$ and $\t'_a=\s_a$ for
$a\ne l$. Notice that $\delta_{\t\u}=\delta_{\s\t'}$ and
$\delta_{\s\t}=\delta_{\u\t'}$. Therefore,
\begin{align*}
s_{\s\t}(l)s_{\t\u}(k)
    &=\frac{\sqrt{e_{\s\s}(l)}\sqrt{e_{\t\t}(l)}-\delta_{\s\t}}
           {c_\s(l)+c_\t(l)}
      \frac{\sqrt{e_{\t\t}(k)}\sqrt{e_{\u\u}(k)}-\delta_{\t\u}}
           {c_\t(k)+c_\u(k)}\\
    &=\frac{\sqrt{e_{\t'\t'}(l)}\sqrt{e_{\u\u}(l)}-\delta_{\u\t'}}
           {c_{\t'}(l)+c_\u(l)}
      \frac{\sqrt{e_{\s\s}(k)}\sqrt{e_{\t'\t'}(k)}-\delta_{\t'\s}}
           {c_\s(k)+c_{\t'}(k)}\\
    &=s_{\t'\u}(l)s_{\s\t'}(k),
\end{align*}
where the second equality uses (\ref{ek}) and (\ref{be})(e).
Hence,
$$ S_kS_lv_\s=\sum_{\t\siml\s,\,\u\simk\t} s_{\s\t}(l)s_{\t\u}(k)v_\u
       =\sum_{\t'\simk\s,\,\u\siml\t} s_{\s\t'}(k)s_{\t'\u}(l)v_\u
       =S_lS_kv_\s,$$
as required.

Assume now that $\s_{k-1}\ne\s_{k+1}$ and $\s_{l-1}=\s_{l+1}$. Then
\begin{align*}
S_kS_lv_\s
   &= \sum_{\t\siml\s} s_{\s\t}(l)\big(a_\t(k)v_\t+b_{\t}(k)v_{S_k\t}\big)\\
   &= a_\s(k)\sum_{\t\siml\s} s_{\s\t}(l)v_\t
                  +b_{\s}(k)\sum_{\t\siml\s}s_{\s\t}(l)v_{S_k\t}\\
   &= a_\s(k)\sum_{\t\siml\s} s_{\s\t}(l)v_\t
                  +b_{\s}(k)\sum_{\t'\siml S_k\s}s_{S_k\s,\t'}(l)v_{\t'}
    =S_lS_kv_\s.
\end{align*}
Interchanging $k$ and $l$ covers the case when
$\s_{k-1}=\s_{k+1}$ and $\s_{l-1}\ne\s_{l+1}$

Finally, consider the case when $\s_{k-1}\ne\s_{k+1}$ and $\s_{l-1}\ne\s_{l+1}$. Then
\begin{align*}
    S_kS_lv_\s &=a_\s(k)a_\s(l)v_\s + a_{S_l\s}(k)b_\s(l)v_{S_l\s}
                  +b_\s(k)a_\s(l)v_{S_k\s}+b_{S_l\s}(k)b_\s(l)v_{S_kS_l\s}\\
               &=a_\s(l)a_\s(k)v_\s + a_\s(k)b_\s(l)v_{S_l\s}
                  +a_{S_k\s}(l)b_\s(k)v_{S_k\s} +b_{S_k\s}(l)b_\s(k)v_{S_lS_k\s},
\end{align*}
since $a_{S_l\s}(k)=a_\s(k)$ and $a_\s(l)=a_{S_k\s}(l)$, by definition,
and $b_{S_l\s}(k)=b_\s(k)$ and $b_{S_k\s}(l)=b_\s(l)$ by (\ref{be})(b).
Hence, $S_kS_lv_\s=S_lS_kv_\s$ if $\s_{k-1}\ne\s_{k+1}$ and
$\s_{l-1}\ne\s_{l+1}$. This completes the proof of (a).
\end{proof}

Finally, we prove that the action of $\W_{r,n}(\bu)$ on
$\Delta(\lambda)$ respects the braid relations of length three.

\begin{Lemma}\label{sts} Suppose that $\s\in \UPD_n(\lambda)$ with
$\s_{k-1}\neq \s_{k+1}$ and $\s_{k}\neq \s_{k+2}$, where $1\le k<n-1$. Then
$S_kS_{k+1}S_kv_{\s}=S_{k+1}S_kS_{k+1} v_{\s}$.
\end{Lemma}

\begin{proof} We consider two cases.

\Case{1. $S_k\s$ is not defined, or $S_k\s$ is defined and
              $(S_k\s)_{k}\neq (S_k\s)_{k+2}$}
First suppose that $S_k\s$ is defined.  If $S_{k+1}\s$ is defined then
$(S_{k+1}\s)_{k-1}\ne(S_{k+1}\s)_{k+1}$, and if $S_{k+1}S_k\s$ is
defined then $(S_{k+1}S_k\s)_{k-1}\ne(S_{k+1}S_k\s)_{k+1}$ because
$\s_k\ne\s_{k+2}$.  Thus we have
$$\begin{array}{l}
S_kS_{k+1}S_k v_{\s}
  =\(a_\s(k)^2a_\s(k+1)+b_{\s}(k)a_{S_k\s}(k+1)b_{S_k\s}(k)\) v_{\s}
\\[4pt]\quad
+\(a_{\s}(k)a_{\s}(k+1)b_{\s}(k)+b_\s(k)a_{S_k\s}(k+1)a_{S_k\s}(k)\)v_{S_k\s}
\\[4pt]\quad
     + a_{\s}(k)b_\s(k+1)a_{S_{k+1}\s}(k) v_{S_{k+1}\s}
     + a_{\s}(k)b_\s(k+1)b_{S_{k+1}\s}(k) v_{S_kS_{k+1}\s}
\\[4pt]\quad
   + b_\s(k)b_{S_{k}\s}(k+1)a_{S_{k+1}S_k\s}(k)v_{S_{k+1}S_k\s}
   + b_\s(k)b_{S_{k}\s}(k+1)b_{S_{k+1}S_k\s}(k)v_{S_kS_{k+1}S_k\s}
\end{array}$$
Now, $\s_{k-1}\ne\s_{k+1}$, or if $S_kS_{k+1}\s$ is
defined, then $(S_kS_{k+1}\s)_k\ne(S_kS_{k+1}\s)_{k+2}$. Therefore, we have
$$\begin{array}{l}
S_{k+1}S_kS_{k+1} v_{\s}
  =\(a_\s(k+1)^2a_\s(k)+b_{\s}(k+1)a_{S_{k+1}\s}(k)b_{S_{k+1}\s}(k+1)\)v_{\s}
\\[4pt]\quad
+\(a_{\s}(k+1)a_{\s}(k)b_{\s}(k+1)
  +b_\s(k+1)a_{S_{k+1}\s}(k)a_{S_{k+1}\s}(k+1)\)v_{S_{k+1}\s}
\\[4pt]\quad
     + a_{\s}(k+1)b_\s(k)a_{S_k\s}(k+1) v_{S_k\s}
     + a_{\s}(k+1)b_\s(k)b_{S_k\s}(k+1) v_{S_{k+1}S_k\s}
\\[4pt]\quad
   + b_\s(k+1)b_{S_{k+1}\s}(k)a_{S_kS_{k+1}\s}(k+1)v_{S_kS_{k+1}\s}
\\[4pt]\quad
   + b_\s(k+1)b_{S_{k+1}\s}(k)b_{S_kS_{k+1}\s}(k+1)v_{S_{k+1}S_kS_{k+1}\s}
\end{array}$$
Now, $b_{S_k\s}(k)=b_{\s}(k)$ and $b_{S_{k+1}\s}(k+1)=b_{\s}(k+1)$ by
(\ref{be})(a). So, in order to check that the coefficients of
$v_{\s}$ are equal in the last two equations we have to show that
$$
a_{\s}(k)^2a_{\s}(k+1)+a_{S_k\s}(k+1)(1-a_{\s}(k)^2)
=a_{\s}(k)a_{\s}(k+1)^2+a_{S_{k+1}\s}(k+1)(1-a_{\s}(k+1)^2);
$$
however, this is just a special case of the easy identity
$$
\frac1{(b-a)^2(c-b)}+\frac1{c-a}\Big(1-\frac1{(b-a)^2}\Big)=
\frac1{(b-a)(c-b)^2}+\frac1{c-a}\Big(1-\frac1{(c-b)^2}\Big).
$$
To see that the coefficients of $v_{S_k\s}$ and $v_{S_{k+1}\s}$
are equal amounts to the following easily checked identities
\begin{align*}
a_{S_k\s}(k)a_{S_k\s}(k+1)+a_{\s}(k)a_{\s}(k+1)&=a_{\s}(k+1)a_{S_k\s}(k+1),\\
a_{S_{k+1}\s}(k)a_{S_{k+1}\s}(k+1)+a_{\s}(k)a_{\s}(k+1)
        &=a_{\s}(k)a_{S_{k+1}\s}(k).
\end{align*}
For the coefficients of $v_{S_{k+1}S_k\s}$ and $v_{S_kS_{k+1}\s}$,
note that
$a_{S_{k+1}S_k\s}(k)=a_{\s}(k+1)$ and $a_{S_kS_{k+1}\s}(k+1)=a_{\s}(k)$.
Finally, three applications of (\ref{be})(c) shows that the
coefficients in $v_{S_kS_{k+1}S_k\s}=v_{S_{k+1}S_kS_{k+1}\s}$ are
equal in both equations.

If $S_k\s$ is not defined then $a_\s(k)=\pm1$ and $b_\s(k)=0$ by
Lemma~\ref{x-equal}(b). Hence, the argument above is still valid if we
set $b_\s(k)=0$.

\Case{2.  $S_k\s$ is defined and $(S_k\s)_{k}=(S_k\s)_{k+2}$}
If $S_{k+1}\s$ is defined then $(S_{k+1}\s)_{k-1}=(S_{k+1}\s)_{k+1}$.
Let $\ts$ be the unique updown tableau such that $\ts\simkk S_k\s$ and
$\ts_{k+1}=\s_{k-1}$. Observe that
if $\t\simkk\ts$ and $\t\ne\ts$ then $\t_{k-1}\ne\t_{k+1}$. Therefore,
$$\begin{aligned}
S_k S_{k+1}&S_k v_\s
  = a_{\s}(k)^2 a_\s(k+1)v_{\s} +a_{\s}(k)a_\s(k+1)b_\s(k) v_{S_k\s}\\
  & + a_{\s}(k)b_\s(k+1)\sum_{\t\simk\ts} s_{S_{k+1}\s,\t}(k)v_{\t}
          +b_\s(k)\sum_{\t\simk\ts} s_{S_k\s,\ts}(k+1) s_{\ts\t}(k)v_{\t}\\
          &+ \sum_{\substack{\t\simkk S_{k}\s\\\t\neq \ts}}
     b_{\s}(k)s_{S_k\s,\t}(k+1)\Big(a_{\t}(k) v_{\t} +b_{\t}(k)v_{S_k\t}\Big),
 \end{aligned}$$
Similarly,
$$\begin{aligned}
S_{k+1} & S_{k}S_{k+1} v_{\s}= a_\s(k+1)^2 a_\s(k)v_{\s}
           +a_\s(k+1)a_\s(k)b_\s(k+1) v_{S_{k+1}\s}\\
   & + a_\s(k+1)b_\s(k)\sum_{\t\simkk\ts} s_{S_{k}\s,\t}(k+1)v_{\t}
            +b_\s(k+1)\sum_{\t\simkk\ts} s_{S_{k+1}\s,\ts}(k)
s_{\ts\t}(k+1)v_{\t}\\
&+ \sum_{\substack{\t\simk S_{k+1}\s\\\t\neq\ts}}
b_\s(k+1)s_{S_{k+1}\s,\t}(k)\(a_\t(k+1) v_{\t} +b_\t(k+1)v_{S_{k+1}\t}\).
\end{aligned}$$
We now compare each of the coefficients in the last two displayed equations.

First we consider the coefficient of $v_\s$.  To show that the
coefficients of $v_{\s}$ are equal in the two expressions above, we
have to prove that
$$\begin{array}{l}
a_{\s}(k)^2 a_\s(k+1)+b_{\s}(k)s_{S_k\s,S_k\s}(k+1)b_{S_k\s}(k)\\[4pt]
 \hspace*{20mm}=  a_\s(k+1)^2
a_\s(k)+b_\s(k+1)s_{S_{k+1}\s,S_{k+1}\s}(k)b_{S_{k+1}\s}(k+1).
\end{array}$$
Now, $b_{\s}(k)=b_{S_k\s}(k)$ and
$b_\s(k+1)=b_{S_{k+1}\s}(k+1)$ by (\ref{be})(a). So, the last
identity is equivalent to
$$\begin{array}{l}
a_{\s}(k)^2 a_\s(k+1)
   +\dfrac{e_{S_k\s,S_k\s}(k+1)-1}{2c_{S_k\s}(k+1)}b_{S_k\s}(k)^2\\[4pt]
 \hspace*{20mm}=  a_\s(k+1)^2a_\s(k)
  +\dfrac{e_{S_{k+1}\s,S_{k+1}\s}(k)-1}{2c_{S_{k+1}\s}(k)}b_{S_{k+1}\s}(k+1)^2.
\end{array}$$
This equation is easily verified using the definitions and
Lemma~\ref{be equality}. Hence, the coefficients of $v_{\s}$ in
$S_kS_{k+1}S_k v_\s$ and $S_{k+1}S_kS_{k+1}v_\s$ are equal.

Now consider the coefficient of $v_{S_k\s}$ in both equations. Since
$a_{S_k\s}(k)-a_\s(k+1)=2c_{S_k\s}(k+1)/
      {(c_{\ts}(k)+c_{S_{k+1}\s}(k))(c_{\ts}(k)+c_{S_k\s}(k))}$,
we see that
\begin{gather*}
s_{S_k\s,S_k\s}(k+1)(a_{S_k\s}(k)-a_{\s}(k+1))b_{\s}(k)
           +a_{\s}(k)a_{\s}(k+1)b_{\s}(k)\\
\begin{aligned}
  &=e_{S_k\s,S_k\s}(k+1)a_{\s}(k)a_{\s}(k+1)b_{\s}(k)\\
  &=\frac{b_{S_k\s}(k)e_{S_k\s,S_k\s}(k+1)}
      {(c_{\ts}(k)+c_{S_{k+1}\s}(k))(c_{\ts}(k+1)+c_{S_k\s}(k+1))}\\
  &=\frac{b_{S_{k+1}\s}(k+1)\sqrt{e_{S_{k+1}\s,S_{k+1}\s}(k)}
             \sqrt{e_{S_k\s,S_k\s}(k+1)}}
    {(c_{\ts}(k)+c_{S_{k+1}\s}(k))(c_{\ts}(k+1)+c_{S_k\s}(k+1))},\\
  &=b_{\s}(k+1)s_{S_{k+1}\s,\ts}(k)s_{\ts,S_k\s}(k+1).
\end{aligned}
\end{gather*}
where the second last equality uses (\ref{be})(f). Consequently,
$$\begin{array}{l}
  a_{\s}(k)a_\s(k+1)b_{\s}(k)
  +b_{\s}(k)s_{S_k\s,S_k\s}(k+1)a_{S_k\s}(k)\\[4pt]
   \hspace*{15mm} =  a_\s(k+1) b_{\s}(k) s_{S_k\s, S_k\s}(k+1)
     + b_\s(k+1)s_{S_{k+1}\s,\ts}(k) s_{\ts, S_k\s}(k+1)
\end{array}$$
Hence, the coefficients of $v_{S_k\s}$ in $S_kS_{k+1}S_kv_\s$ and
$S_{k+1}S_kS_{k+1}v_\s$ are equal. A similar argument shows that
$$\begin{array}{l}
a_{\s}(k)b_\s(k+1)s_{S_{k+1}\s,S_{k+1}\s}(k)
    +b_{\s}(k)s_{S_k\s, \ts}(k+1)s_{\ts,S_{k+1}\s}(k)\\[4pt]
\hspace*{20mm}= a_\s(k+1)a_\s(k)b_\s(k+1)
    +b_\s(k+1)s_{S_{k+1}\s,S_{k+1}\s}(k)a_{S_{k+1}\s}(k+1).
\end{array}$$
This proves that
the coefficient of $v_{S_{k+1}\s}$ in $S_kS_{k+1}S_kv_\s$ and in
$S_{k+1}S_kS_{k+1}v_\s$ are equal.

Now consider the coefficient of $v_\t$ where $\t\simk\ts$ and
$\t\not\in\{\ts, S_{k+1}\s\}$. This time
$$
a_{\t}(k+1)-a_{\s}(k)=
\dfrac{c_{S_{k+1}\s}(k)+c_{\t}(k)}
{(c_{S_k\s}(k+1)+c_{\ts}(k+1))(c_{\ts}(k)+c_{\t}(k))}.
$$
An argument similar to that for $v_{S_k\s}$ now shows that
$$
b_{\s}(k) s_{S_k\s,\ts}(k+1) s_{\ts\t}(k)
     =b_\s(k+1)s_{S_{k+1}\s,\t}(k) \( a_\t(k+1)-a_{\s}(k)\).$$
Therefore, the coefficients of $v_{\t}$ for such $\t$ in
$S_kS_{k+1}S_k v_\s$ and $S_{k+1}S_kS_{k+1}v_\s$ are equal.

Another variation of this argument shows that if $\t\simkk S_k\s$ and
$\t\not\in\{\ts,S_k\s\}$ then the coefficients of $v_\t$ in
$S_kS_{k+1}S_kv_\s$ and $S_{k+1}S_kS_{k+1}v_\s$ are both equal.

Next, we suppose that $S_k\t$ is defined and we compare the
coefficients of $v_{S_k\t}$ in $S_kS_{k+1}S_kv_\s$ and
$S_{k+1}S_kS_{k+1}v_\s$, when $\t\simkk S_k\s$ and
$\t\notin\{\ts,S_k\s\}$.  As $S_k\s$ is defined,  $\u=S_{k+1}S_k\t$ is
defined and $\u\simk S_{k+1}\s$ with $\u\not\in\{\ts,S_{k+1}\s\}$.
Conversely, if $S_{k+1}\u$ is defined for such $\u$ then
$\t=S_kS_{k+1}\u$ is defined.  Applying (\ref{be})(f) twice, we
have
$$
b_{\s}(k)b_{\t}(k)\sqrt{e_{S_k\s,S_k\s}(k+1)}\sqrt{e_{\t\t}(k+1)}
  =b_{\s}(k+1)b_{\u}(k+1)\sqrt{e_{S_{k+1}\s,S_{k+1}\s}(k)}
         \sqrt{e_{\u\u}(k)}.
$$
Consequently, because
$c_{S_k\s}(k+1)+c_\t(k+1)=c_\u(k)+c_{S_{k+1}\s}(k)$, we have
$$
b_{\s}(k) s_{S_k\s, \t}(k+1)b_{\t}(k)
    =b_\s(k+1) s_{S_{k+1}\s,\u}(k) b_\u(k+1).
$$
That is, the coefficients of $v_{S_k\t}$ in
$S_{k+1}S_kS_{k+1} v_{\s}$ and $S_{k} S_{k+1} S_{k} v_{\s}$ are equal.

It remains to compare the coefficients of $v_{\ts}$ in the two
equations. To show that these two coefficients are equal we have to
prove that
$$\begin{array}{l}
a_\s(k)b_\s(k+1)s_{S_{k+1}\s,\ts}(k)
    +b_{\s}(k)s_{S_k\s,\ts}(k+1)s_{\ts\ts}(k)\\[4pt]
   \hspace*{20mm} =a_\s(k+1)b_\s(k) s_{S_k\s, \ts}(k+1)
          +b_\s(k+1) s_{S_{k+1}\s,\ts}(k) s_{\ts\ts}(k+1).
\end{array}$$
First note that, by the definitions and (\ref{be})(a),
\begin{align*}
b_{\s}(k+1)s_{S_{k+1}\s,\ts}(k)
   &=\frac{b_{\s}(k+1)\sqrt{e_{S_{k+1}\s,S_{k+1}\s}(k)}
               \sqrt{e_{\ts,\ts}(k)}} {c_{S_{k+1}\s}(k)+c_{\ts}(k)}\\
   &=\frac{b_{\s}(k)\sqrt{e_{S_{k+1}\s,S_{k+1}\s}(k)}
               \sqrt{e_{\ts,\ts}(k+1)}} {c_{S_{k+1}\s}(k)+c_{\ts}(k)}\\
   &= b_{\s}(k)s_{S_k\s,\ts}(k+1)
   \frac{c_{S_k\s}(k+1)+c_{\ts}(k+1)}{c_{S_{k+1}\s}(k)+c_{\ts}(k)}
    e_{\ts\ts}(k).
\end{align*}
So, it is enough to show that
$$
\(c_{S_k\s}(k+1)+c_{\ts}(k+1)\)e_{\ts\ts}(k)
       \(a_{\s}(k)-s_{\ts\ts}(k+1)\)
  =\(c_{S_{k+1}\s}(k)+c_{\ts}(k)\)\(a_{\s}(k+1)-s_{\ts\ts}(k)\);
$$
however, this follows from Lemma~\ref{ees}. Hence,
the coefficients of $ v_{\ts}$ in
$S_{k+1}S_kS_{k+1} v_{\s}$ and $S_{k} S_{k+1} S_{k} v_{\s}$ are
equal.

This completes the proof of Lemma~\ref{sts}.
\end{proof}

\begin{Lemma}\label{sts rel2} Suppose that $\s\in \UPD_{n}(\lambda)$
and that either $\s_{k-1}=\s_{k+1}$ and $\s_k\neq \s_{k+2}$, or
$\s_{k-1}\neq \s_{k+1}$ and $\s_k= \s_{k+2}$, for $1\le k<n-1$.
Then $S_kS_{k+1}S_k v_{\s}=S_{k+1}S_kS_{k+1} v_{\s}$.
\end{Lemma}

\begin{proof}There are again two cases to consider.

\Case{1. $S_{k+1}\s$ is defined}
Suppose first that $\s_{k-1}=\s_{k+1}$ and $\s_k\neq
\s_{k+2}$. Then $\t=S_{k+1}\s\in \UPD(\lambda)$ is well-defined.
Furthermore, $\t_k\neq \t_{k+2}$ and $\t_{k-1}\neq \t_{k+1}$, so
$S_kS_{k+1}S_k v_{\t}=S_{k+1}S_kS_{k+1} v_{\t}$ by Lemma~\ref{sts}.
Now, $S_{k+1}v_\t=a_\t(k+1)v_\t+b_\t(k+1)v_\s$ and $b_\t(k+1)\ne0$.
Therefore
\begin{align*}
S_kS_{k+1}S_kv_\s
   &=\frac1{b_\t(k+1)}S_kS_{k+1}S_k\Big(S_{k+1}v_\t-a_\t(k+1)v_\t\Big)\\
   &=\frac1{b_\t(k+1)}\Big(S_k(S_{k+1}S_kS_{k+1})v_\t-a_\t(k+1)(S_kS_{k+1}S_k)v_\t\Big)\\
   &=\frac1{b_\t(k+1)}\Big(S_k(S_kS_{k+1}S_k)v_\t
                  -a_\t(k+1)(S_{k+1}S_kS_{k+1})v_\t\Big)\\
\intertext{by Lemma~\ref{sts}. Hence, using Lemma~\ref{s-square} twice,}
S_kS_{k+1}S_kv_\s
   &=\frac1{b_\t(k+1)}\Big(S_{k+1}S_kv_\t-a_\t(k+1)(S_{k+1}S_kS_{k+1})v_\t\Big)\\
   &=\frac1{b_\t(k+1)}\Big(S_{k+1}S_k(S_{k+1}S_{k+1})v_\t
            -a_\t(k+1)(S_{k+1}S_kS_{k+1})v_\t\Big)\\
   &=\frac1{b_\t(k+1)}(S_{k+1}S_kS_{k+1})\Big(S_{k+1}v_\t
                  -a_\t(k+1)v_\t\Big)\\
   &=(S_{k+1}S_kS_{k+1})v_\s
\end{align*}
as required.

The case when $\s_{k-1}\neq \s_{k+1}$ and $\s_k= \s_{k+2}$ can be
proved similarly.

\Case{2. $S_{k+1}\s$ is not defined}
This is  equivalent to saying
that the two nodes $\s_{k+2}\ominus\s_{k+1}$ and $\s_{k+1}\ominus
\s_k$ are in the same row or in the same column. Therefore, either
$\s_k\subset \s_{k+1}\subset \s_{k+2}$ or $\s_k\supset
\s_{k+1}\supset \s_{k+2}$. Note that in either case
$\s_{k-1}=\s_{k+1}$, so we have
$$ E_k v_\s=\sum_{\substack{\t\siml\s\\\t\neq \s}}
        e_{\s\t}(k) v_\t +e_{\s\s}(k) v_\s.$$
By Proposition~\ref{see} and Proposition~\ref{se},
$S_{k}S_{k+1}S_kE_k v_\s=S_kS_{k+1}E_kv_\s=E_{k+1}E_k v_\s$ and
$S_{k+1}S_{k}S_{k+1}E_kv_\s=S_{k+1}E_{k+1}E_k v_\s=E_{k+1}E_kv_\s$.

Suppose that $\t\simk\s$ and $\t\ne\s$. Then $S_{k+1}\t$ is
well-defined and $\t_{k-1}=\t_{k+1}$---indeed, the two boxes
$\s_{k+2}\ominus\s_{k+1}$ and $\s_{k+1}\ominus \t_k$ belong to
different rows and columns. Hence, by Case~1,
$S_{k+1}S_{k}S_{k+1} v_\t=S_kS_{k+1}S_{k} v_\t$. Consequently,
$S_{k+1}S_{k}S_{k+1}e_{\s\s}(k)v_\s=S_k S_{k+1}S_{k}e_{\s\s}(k)v_\s$.
Cancelling the non-zero factor $e_{\s\s}(k)$ shows that
$S_{k}S_{k+1}S_k v_\s=S_{k+1}S_kS_{k+1} v_\s$.
\end{proof}

\begin{Prop}\label{end relations} Suppose that $1\le k<n-1$
    and $\s\in \UPD_{n}(\lambda)$. Then
    $S_kS_{k+1}S_k v_{\s}=S_{k+1}S_kS_{k+1} v_{\s}$.
\end{Prop}

\begin{proof} By Lemma~\ref{sts} and Lemma~\ref{sts rel2} it only remains
    to consider the case when $\s_{k-1}=\s_{k+1}$ and $\s_k= \s_{k+2}$.
By Lemma~\ref{s-square}, Proposition~\ref{se} and Proposition~\ref{see}(a), we have
\begin{align*}
    S_{k+1}S_kS_{k+1}E_kv_\s
      &= S_{k+1}S_k\cdot S_kE_{k+1}E_kv_\s
       = S_{k+1}E_{k+1}E_kv_\s
       = E_{k+1}E_kv_\s,\\
\intertext{on the one hand. Similarly, we also have}
    S_kS_{k+1}S_kE_kv_\s
      &= S_kS_{k+1}E_kv_\s
       = S_k\cdot S_kE_{k+1}E_kv_\s
       = E_{k+1}E_kv_\s,
\end{align*}
Therefore, recalling the definition of $E_kv_\s$, we have
$$
\(S_{k+1}S_kS_{k+1}-S_kS_{k+1}S_k\)
   \Big(e_{\s\s}(k)v_\s+\sum_{\t\simk\s, \t\ne\s} e_{\s\t}(k)v_\t\Big)=0.$$
Now, if $\t\simk \s$ and $\t\neq \s$ then
$S_kS_{k+1}S_k v_{\t}=S_{k+1}S_k S_{k+1}v_{\t}$ by Lemma~\ref{sts rel2}.
Consequently, $S_kS_{k+1}S_k v_{\s}=S_{k+1}S_k S_{k+1}v_{\s}$ since
$e_{\s\s}(k)\ne0$. This completes the proof.
\end{proof}

\begin{proof}[Proof of Theorem~\ref{seminormal}]
    The results from Lemma~\ref{start relations} to
    Proposition~\ref{end relations} show that the action of the
    generators of $\W_{r,n}(\bu)$ on $\Delta(\lambda)$ respects all of
    the relations of $\W_{r,n}(\bu)$.  Hence, $\Delta(\lambda)$ is a
    $\W_{r,n}(\bu)$--module, as we wanted to show.
\end{proof}

\section{Irreducible representations and Theorem~A}

In this section we use the seminormal representations to show that the
cyclotomic Nazarov--Wenzl algebras are always free of rank
$r^n(2n-1)!!$.  Before we can do this we need to recall some
identities involving updown tableaux.

First, if $\lambda$ is a multipartition of $n-2m$ let
$f^{(n,\lambda)}$ be the number of $n$--updown $\lambda$--tableaux.
So, in particular, $f^{(|\lambda|,\lambda)}=\#\Std(\lambda)$ is the
number of standard $\lambda$--tableaux.
Sundaram~\cite[Lemma~8.7]{Sundaram:PhD} has given a combinatorial
bijection to show that if $\tau$ is a \textit{partition} (so $r=1$)
then the number of $n$--updown $\tau$--tableaux is equal to $\binom
n{|\tau|}(n-|\tau|-1)!!f^{(|\tau|,\tau)}$. Terada~\cite{Terada:brauer}
has given a geometric version of this bijection when $|\tau|=0$ and
$n$ is even.

\begin{Lemma}\label{updown}
Suppose that $0\le m\le\floor{n2}$ and that
$\lambda\in\Lambda_r^+(n-2m)$. Then
$$f^{(n,\lambda)}=r^m\binom n{2m}(2m-1)!!\#\Std(\lambda).$$
\end{Lemma}

\begin{proof} Using Sundaram's formula from above we have
\begin{align*}
    f^{(n,\lambda)}&=\sum_{\substack{n_1,\dots,n_r\\n_1+\dots+n_r=n\\
                    n_t-|\lambda^{(t)}|\in2\Z}}
         \binom n{n_1,\dots,n_r}
           \prod_{t=1}^r\binom{n_t}{|\lambda^{(t)}|}
           (n_t-|\lambda^{(t)}|-1)!!f^{(|\lambda^{(t)}|,\lambda^{(t)})}\\
           &=\sum_{\substack{n_1,\dots,n_r\\n_1+\dots+n_r=n\\
                      n_t-|\lambda^{(t)}|\in2\Z}}n!
  \prod_{t=1}^r\frac{(n_t-|\lambda^{(t)}|-1)!!f^{(|\lambda^{(t)}|,\lambda^{(t)})}}
              {(n_t-|\lambda^{(t)}|)!|\lambda^{(t)}|!}\\
  &=n!\prod_{t=1}^r\frac{f^{(|\lambda^{(t)}|,\lambda^{(t)})}}
                                   {|\lambda^{(t)}|!}
        \sum_{\substack{n_1,\dots,n_r\\n_1+\dots+n_r=n}}
          \prod_{t=1}^r\frac{(n_t-|\lambda^{(t)}|-1)!!}{(n_t-|\lambda^{(t)}|)!}\\
  &=\frac{n!}{(n-2m)!}\#\Std(\lambda)
        \sum_{\substack{a_1,\dots,a_r\\a_1+\dots+a_r=m}}
          \prod_{t=1}^r\frac{(2a_t-1)!!}{(2a_t)!},
\intertext{where the summation is now over
$a_t=\frac{n_t-|\lambda^{(t)}|}2$, for $1\le t\le r$. Hence}
    f^{(n,\lambda)}&=\frac{n!}{(n-2m)!}\#\Std(\lambda)
        2^{-m}\sum_{\substack{a_1,\dots,a_r\\a_1+\dots+a_r=m}}
              \prod_{t=1}^r\frac1{a_t!}\\
    &=\frac{n!}{(n-2m)!}\#\Std(\lambda)\frac{r^m}{2^m m!}
     =r^m\binom n{2m}(2m-1)!!\#\Std(\lambda).
\end{align*}
\end{proof}

It is well--known from the representation theory of the degenerate Hecke
algebras $\H_{r,k}$ that $\sum_\lambda\#\Std(\lambda)^2=r^kk!$, where
in the sum $\lambda\in\Lambda_r^+(k)$.

\begin{Cor}\label{(2n-1)!!}
Suppose that $n\ge1$ and $r\ge1$. Then
    $$\sum_{m=0}^{\floor{n2}}\sum_{\lambda\vdash n-2m}{f^{(n,\lambda)}}^2
                 =r^n(2n-1)!!.$$
\end{Cor}

\begin{proof}Using the Lemma we have
\begin{align*}
\sum_{m=0}^{\floor{n2}}\sum_{\lambda\vdash n-2m}{f^{(n,\lambda)}}^2
  &=\sum_{m=0}^{\floor{n2}}\sum_{\lambda\vdash n-2m}
             \Big\{r^m\binom n{2m}(2m-1)!!\#\Std(\lambda)\Big\}^2.\\
  &=\sum_{m=0}^{\floor{n2}}r^{2m}\binom n{2m}^2\big((2m-1)!!\big)^2
        \sum_{\lambda\vdash n-2m} \#\Std(\lambda)^2.\\
  &=\sum_{m=0}^{\floor{n2}}r^{2m}\binom n{2m}^2\big((2m-1)!!\big)^2
              r^{n-2m}(n-2m)!\\
  &=r^n\sum_{m=0}^{\floor{n2}}\binom n{2m}^2\big((2m-1)!!\big)^2(n-2m)!
\end{align*}
To complete the proof, notice that the sum on the right hand side does
not depend on $r$, so we can set $r=1$ and deduce the result from the
representation theory of the Brauer algebras.
\end{proof}

A representation theoretic proof of this result is given in
\cite{RuiYu} where it is obtained as a consequence of the branching
rules for the cyclotomic Brauer algebra. The cell modules of the
cyclotomic Brauer algebras are indexed by the multipartitions of
$n-2m$, for $0\le m\le\floor{n2}$. The branching
rule~\cite[Theorem~6.1]{RuiYu} shows that the dimension of the cell
module indexed by $\lambda$ is $f^{(n,\lambda)}$. On the other hand,
the cellular basis of the cyclotomic Brauer algebras constructed
in~\cite[Theorem~5.11]{RuiYu} contains $r^n(2n-1)!!$ elements.
Combining these two facts proves the result.

Given two multipartitions $\lambda$ and $\mu$ such that $\mu$ is
obtained by adding a box to $\lambda$ we write
$\lambda\rightarrow\mu$, or $\mu\leftarrow\lambda$.

\begin{Theorem}\label{generic}
    Suppose that $R$ is a field  with $\Char R>2n$ and the root
    conditions (Assumption~\ref{be}) hold in $R$. Assume that the
    parameters $u_1,\dots,u_r$ are generic for $\W_{r,n}(\bu)$ and
    that $\Omega$ is $\bu$--admissible. Then:
\begin{enumerate}
\item Suppose $n>1$. There is a $\W_{r,n-1}(\bu)$-module isomorphism
$$
\Delta(\lambda)\downarrow
   =\bigoplus_{\substack{\mu\\\mu\rightarrow\lambda}}\Delta(\mu)
      \quad\bigoplus\quad
      \bigoplus_{\substack{\nu\\\lambda\rightarrow\nu}} \Delta(\nu).
$$
where $\Delta(\lambda)\downarrow$ is $\Delta(\lambda)$ considered as a
$\W_{r,n-1}(\bu)$-module.
\item The seminormal representation $\Delta(\lambda)$ is an
        irreducible $\W_{r,n}(\bu)$--module for each multipartition
        $\lambda$ of $n-2m$, where $0\le m\le\floor{n2}$.
\item The set
       $\set{\Delta(\lambda)|\lambda\vdash n-2m, 0\le m\le\floor{n2}}$
       is a complete set of  irreducible
       $\W_{r,n}(\bu)$--modules.
\item $\W_{r,n}(\bu)$ is a split semisimple $R$--algebra of dimension
$r^n(2n-1)!!$.
     \end{enumerate}
\end{Theorem}

\begin{proof}
Part (a) follows if we define  $\Delta(\mu)$ to be  the vector subspace
spanned by $v_\t$ with $\t\in \UPD_n(\lambda)$ and $\t_{n-1}=\mu$.

Let $\mathscr X=\langle X_1,\dots,X_n\rangle$. Since $X_k
v_\s=c_\s(k) v_\s$, for all $\s\in\UPD_n(\lambda)$ and $1\le k\le n$,
the seminormal representation
$\Delta(\lambda)=\bigoplus_{\s\in\UPD_n(\lambda)}Rv_\s$ decomposes
into a direct sum of one dimensional submodules as an
$\mathscr X$--module. Further, by Lemma~\ref{generic u}(a), this
decomposition is multiplicity free. In particular,
$\Delta(\lambda)\cong\Delta(\mu)$ if and only if $\lambda=\mu$.
Further, if $M$ is a $\W_{r,n}(\bu)$-submodule of $\Delta(\lambda)$ then
$M$ is spanned by some subset of $\set{v_\s|\s\in\UPD_n(\lambda)}$.

To prove (b) we now argue by induction on $n$.  If $n=1$ then
$\Delta(\lambda)$ is one dimensional and hence irreducible, for all
$\lambda$. Suppose now that $n>1$ and let $M\subset \Delta(\lambda)$
be a non--zero $\W_{r,n}(\bu)$--submodule of $\Delta(\lambda)$. By the
remarks in the last paragraph, $M$ is spanned by a subset of
$\set{v_\s|\s\in\UPD_n(\lambda)}$. Therefore, if we consider $M$ as a
$\W_{r,n-1}(\bu)$--module then $M\supset\Delta(\mu)$, for some
multipartition $\mu$ which is obtained by adding or removing a node
from $\lambda$.

\Case{1. $|\lambda|=n$}
Since $|\lambda|=n$, The multipartition $\mu$ is obtained from
$\lambda$ by removing a node. If
$\lambda=((0),\dots,(0),(a^b),(0),\dots,(0))$ then
$\Delta(\lambda)\downarrow$ is irreducible as a $\W_{r,n-1}(\bu)$--module,
so there is nothing to prove. Suppose then that $\lambda$ is not of
this form and that $\nu$ is a different multipartition which is obtained
from $\lambda$ by removing a node. Let $\s\in\UPD_n(\lambda)$ be an
updown tableau such that $\s_{n-1}=\mu$ and $\mu\setminus\s_{n-2}=\lambda\setminus\nu$.
So $v_{\s}\in \Delta(\mu)\subset M$ and $(S_{n-1}\s)_{n-1}=\nu$. Now,
$$S_{n-1} v_\s=a_\s(n-1)v_\s+b_\s(n-1) v_{S_{n-1} \s}\in M,$$
and $b_\s(n-1)\ne0$ since $\lambda{\setminus}\mu$ and $\lambda{\setminus}\nu$ cannot be in
the same row or in the same column. Consequently, $v_{S_{n-1}\s}\in
M$. This implies that $\Delta(\nu)\subset M$ since $(S_{n-1}\s)_{n-1}=\nu$.
Therefore, $\sum_{\nu\rightarrow\lambda}\Delta(\nu)\subset M$, so $M=\Delta(\lambda)$ by part(a). Hence, $\Delta(\lambda)$ is irreducible as required.

\Case{2. $|\lambda|<n$}
Since $|\lambda|<n$, $\UPD_{n-2}(\lambda)$ is non--empty so we fix
$\t\in \UPD_{n-2}(\lambda)$. Let
$\s=(\t_1,\dots,\t_{n-2},\mu,\lambda)$, then $\s\in\UPD_n(\lambda)$
and $v_\s\in\Delta(\mu)\subset M$.  Then
$$E_{n-1} v_\s=\sum_{\u\overset{n-1}\sim\s}e_{\s\u}(n-1)v_\u\in M.$$
As $e_{\s\u}(n-1)\ne0$ whenever $\u\overset{n-1}\sim\s$, we have
$v_\u\in M$ for each term in this sum. If
$\nu\leftarrow\lambda$ or $\nu\rightarrow\lambda$ then
$\u=(\t_1,\dots,\t_{n-2},\nu,\lambda)\overset{n-1}\sim\s$, so
$\Delta(\nu)\subset M$.  Hence, $M=\Delta(\lambda)$ and
$\Delta(\lambda)$ is irreducible as claimed.  This completes the proof
of (b).

Finally, we prove (c) and (d). We have already seen that the
seminormal representations are pairwise non--isomorphic, so it remains
to show that every irreducible is isomorphic to $\Delta(\lambda)$ for
some $\lambda$. Let $\Rad\W_{r,n}(\bu)$ be the Jacobson radical
of~$\W_{r,n}(\bu)$. Then
\begin{align*}
 \dim_R \W_{r,n}(\bu)&\ge \dim_R(\W_{r,n}(\bu)/\Rad\W_{r,n}(\bu))
                  \ge \sum_{m=0}^{\lfloor n/2\rfloor}
      \sum_{\lambda\vdash n-2m} \(\dim_R\Delta(\lambda)\)^2.\\
\intertext{By construction,
$\dim\Delta(\lambda)=\#\UPD_n(\lambda)=f^{(n,\lambda)}$. So using
Corollary~\ref{(2n-1)!!}, and then Proposition~\ref{spanning set}, we have}
 \dim_R \W_{r,n}(\bu) &\ge r^n (2n-1)!!\ge \dim_R\W_{r,n}(\bu).
\end{align*}
Therefore, $\Rad\W_{r,n}(\bu)=0$, which forces
$\dim_{R} \W_{r,n}(\bu)=r^n (2n-1)!!$.  Now, parts (c) and (d) both follow
from the Wedderburn-Artin Theorem.
\end{proof}

Before establishing a strong version of Theorem~A, we show that the
Root conditions (Assumption~\ref{be}) can be satisfied when
$R=\R$.

\begin{Lemma}\label{be real}
Suppose that $R=\R$ and
we choose $u_i\in R$ in such a way that
\begin{enumerate}
\item
$|u_1|>\cdots>|u_r|\ge n$ and
$|u_i|-|u_{i+1}|\ge 2n$,
\item
$u_i<0$ if $i$ is even and $u_i>0$ if
$i$ is odd.
\end{enumerate}\noindent%
Suppose that $\s\in\UPD_n(\lambda)$ and $1\le k<n$. Then
$|a_{\s}(k)|\le1$, if
$\s_{k-1}\ne\s_{k+1}$, and
$e_{\s\s}(k)>0$, if $\s_{k-1}=\s_{k+1}$.
In particular, the Root Condition (\ref{be}) holds
if we choose non--negative square roots $\sqrt{b_{\s}(k)}\ge0$ and
$\sqrt{e_{\s\s}(k)}>0$.
\end{Lemma}

\begin{proof}
We start with the case $\s_{k-1}\ne\s_{k+1}$.  Let
$\alpha=\s_k\ominus\s_{k-1}$ and $\beta=\s_{k+1}\ominus\s_k$. Note
that $c(\alpha)+c(\beta)\ne0$.  Write $\alpha=(i,j,t)$ and
$\beta=(i',j',t')$.  If $t=t'$ and both nodes are addable, or both
nodes are removable, then $\alpha\ne\beta$.  Thus,
$c(\beta)-c(\alpha)$ is a nonzero integer and $|a_\s(k)|\le1$. If
$t=t'$ and only one of the nodes is addable (and the other is
removable), then
$$
\frac1{|a_\t(k)|}=|c(\alpha)-c(\beta)|=|2u_t+(j-i)+(j'-i')|
\ge 2|u_t|-2(n-1)\ge2.
$$
Hence, $|a_t(k)|\le1$ if $t=t'$. A similar argument shows that
$|a_\s(k)|\le1$ when $t\ne t'$.

Next we consider the case $\s_{k-1}=\s_{k+1}$.
Let $\alpha=\s_k\ominus\s_{k-1}$ and $\lambda=\s_{k-1}$.
Write $\alpha=(i,j,t)$. By (\ref{ek}) and because
$R=\R$, we have
$$
e_{\s\s}(k)=\big(2c(\alpha)-(-1)^r\big)
         \prod_{\beta}\frac{c(\alpha)+c(\beta)}{c(\alpha)-c(\beta)},
$$
where $\beta$ runs over all of the addable and removable nodes of
$\lambda$ with $\beta\neq\alpha$.

Suppose that $t$ is even. First we show that
$$\prod_{\beta\not\in\lambda^{(t)}}
 \frac{c(\alpha)+c(\beta)}{c(\alpha)-c(\beta)}<0.$$
Consider the contents of all of the addable and removable nodes of
$\lambda^{(t')}$, where $t'\ne t$.  If $t'$ is even then there are $l$
positive contents $|u_{t'}|+d_j$ with $|d_j|<n$, for $1\le j\le l$,
and $l+1$ negative contents $-|u_{t'}|-c_i$ with $|c_i|<n$, for
$1\le i\le l+1$.  Let~$\epsilon_{t'}$ be the sign of the product of
$\frac{c(\alpha)+c(\beta)} {c(\alpha)-c(\beta)}$ over all addable and
removable nodes~$\beta$ of $\lambda^{(t')}$.  Our aim is to show that
$$\prod_{t'\ne t}\epsilon_{t'}=-1.$$
By our assumptions, $\epsilon_{t'}$ is
equal to the sign of
$$
\frac{(-|u_t|+|u_{t'}|)^l}{(-|u_t|-|u_{t'}|)^l}
\frac{(-|u_t|-|u_{t'}|)^{l+1}}{(-|u_t|+|u_{t'}|)^{l+1}}
=\frac{|u_t|+|u_{t'}|}{|u_t|-|u_{t'}|}.
$$
Thus, $\epsilon_{t'}<0$ if and only if $t'<t$.  If $t'$ is odd then
there are $l+1$ positive contents $|u_{t'}|+c_i$ with $|c_i|<n$, for
$1\le i\le l+1$, and $l$ negative contents $-|u_{t'}|-d_j$ with
$|d_j|<n$, for $1\le j\le l$. Then, by the same argument,
$\epsilon_{t'}<0$ if and only if $t'<t$ again. Thus
$$\prod_{t'\ne t}\epsilon_{t'}=(-1)^{t-1}=-1.$$

Let $-|u_t|-c_i$, for $1\le i\le l+1$, be the contents of the addable
nodes of $\lambda^{(t)}$ and let $|u_t|+d_j$, for $1\le j\le l$, be
the contents of the removable nodes of $\lambda^{(t)}$. We may assume
that
$$c_1>d_1>\cdots>c_l>d_l>c_{l+1}.$$
Let $\epsilon_t$ be the sign of the product of
$\frac{c(\alpha)+c(\beta)}{c(\alpha)-c(\beta)}$, where $\beta$ runs
over all of the addable and removable nodes of $\lambda^{(t)}$ such
that $\beta\ne\alpha$.

If $c(\alpha)=-|u_t|-c_i$, for some $i$, then $\epsilon_t$ is equal to
the sign of $$\prod_{k\ne i}\frac{-2|u_t|-c_i-c_k}{c_k-c_i}
\prod_{k=1}^l\frac{d_k-c_i}{-2|u_t|-c_i-d_k},$$
so
$\epsilon_t=\frac{(-1)^l}{(-1)^{l+1-i}}
    \frac{(-1)^{l-i+1}}{(-1)^l}=1$.
As $2c(\alpha)-(-1)^r=-2|u_t|-2c_i\pm1<0$ and
$$\prod_{1\le t'\le r}\epsilon_{t'}=-1,$$
we have $e_{\s\s}(k)>0$.

If $c(\alpha)=|u_t|+d_j$, for some $j$, then
$\epsilon_t$ is equal to the sign of
$$\prod_{k=1}^{l+1}\frac{d_j-c_k}{2|u_t|+d_j+c_k}
\prod_{k\ne j}^l\frac{2|u_t|+d_j+d_k}{d_j-d_k},$$
so
$\epsilon_t=\frac{(-1)^j}{(-1)^{j-1}}=-1$.
As $2c(\alpha)-(-1)^r=2|u_t|+2d_j\pm1>0$ and
$$\prod_{1\le t'\le r}\epsilon_{t'}=-1,$$
we have $e_{\s\s}(k)>0$ again.

The case when $t$ is odd is handled similarly.
In this case, we have
$$\prod_{\beta\not\in\lambda^{(t)}}
 \frac{c(\alpha)+c(\beta)}{c(\alpha)-c(\beta)}>0,$$
because its sign is equal to $(-1)^{t-1}=1$.
Let $|u_t|+c_i$, for $1\le i\le l+1$, be the
contents of the addable nodes of $\lambda^{(t)}$ and let $-|u_t|-d_j$, for
$1\le j\le l$, be the contents of the removable nodes of $\lambda^{(t)}$
such that
$$c_1>d_1>\cdots>c_l>d_l>c_{l+1}.$$
If $c(\alpha)=|u_t|+c_i$, for some $i$, then
$\epsilon_t$ is equal to the sign of
$$\prod_{k\ne i}\frac{2|u_t|+c_i+c_k}{c_i-c_k}
\prod_{k=1}^l\frac{c_i-d_k}{2|u_t|+c_i+d_k},$$
so
$\epsilon_t=\frac{(-1)^{i-1}}{(-1)^{i-1}}=1$.
As $2c(\alpha)-(-1)^r>0$ we have $e_{\s\s}(k)>0$.

If $c(\alpha)=-|u_t|-d_j$, for some $j$, then
$\epsilon_t$ is equal to the sign of
$$\prod_{k=1}^{l+1}\frac{c_k-d_j}{-2|u_t|-d_j-c_k}
\prod_{k\ne j}^l\frac{-2|u_t|-d_j-d_k}{d_k-d_j},$$
so
$\epsilon_t=\frac{(-1)^{l-j+1}}{(-1)^{l+1}}
\frac{(-1)^{l-1}}{(-1)^{l-j}}=-1$.
As $2c(\alpha)-(-1)^r<0$
we have $e_{\s\s}(k)>0$ again.
\end{proof}

We can now prove a stronger version of Theorem~A.

\begin{Theorem}\label{W-basis}
    Suppose that $R$ is a commutative ring in which $2$ is
    invertible and that $\Omega$ is $\bu$--admissible. Then
    $\W_{r,n}(\bu)$ is free as an $R$--module with basis the set of
    $r$--regular monomials. Consequently, $\W_{r,n}(\bu)$ is free of
    rank $r^n(2n-1)!!$.
\end{Theorem}

\begin{proof}Recall that if $R$ is a ring in which $2$ is invertible
then $\W_{r,n}(\bu)$ is spanned by the set of $r$-regular monomials by
Proposition~\ref{spanning set}. For convenience, if $S$ is a ring and
$\bu_s\in S^r$ then we let $\W_S(\bu_S)$ be the cyclotomic
Nazarov--Wenzl algebra defined over~$S$ with parameters $\bu_S$.

First, we consider the special case when $R=\mathcal Z$, where
$\ZC=\Z[\frac12,\dot u_1,\dots,\dot u_r]$ and the $\dot u_i$ are
indeterminates over $\Z$. Let $\dot\bu=(\dot u_1,\dots,\dot u_r)$,
define $\dot\Omega$ in accordance with Definition~\ref{u-admissible}
and  consider the cyclotomic Nazarov--Wenzl algebra $\W_\ZC(\dot\bu)$.
As $\R$ is not finitely generated over $\Q$ we can find~$r$
algebraically independent transcendental real numbers
$u_1',\dots,u_r'\in\R$ which satisfy the hypotheses of Lemma~\ref{be
real}. Let $\ZC'=\Z[\frac12,u_1',\dots,u_r']$ and let
$\theta\map{\ZC}{\ZC'}$ be the $\ZC$--linear map determined by
$\theta(\dot u_i)=u_i'$, for $1\le i\le r$. Then $\theta$ is a ring
isomorphism. Let $\bu'=(u_1',\dots,u_r')$ and
$\Omega'=\set{\theta(\dot\omega_a)|a\ge0}$. Then $\Omega'$ is
$\bu'$--admissible and $\theta$ induces an isomorphism of
$\ZC$--algebras $\W_\ZC(\dot\bu)\cong\W_{\ZC'}(\bu')$, where the
inverse map is the homomorphism induced by $\theta^{-1}\map{\ZC'}\ZC$.

Now, by Lemma~\ref{be real} and Theorem~\ref{generic}(d),
$\W_\R(\bu')$ is an $\R$--algebra of dimension $r^n(2n-1)!!$. Hence
the set of $r$--regular monomials is an $\R$--basis of $\W_\R(\bu')$
since there are $r^n(2n-1)!!$ $r$-regular monomials. In particular,
the set of $r$-regular monomials is linearly independent over $\R$,
and hence linearly independent over~$\ZC'$.  Therefore,
$\W_{\ZC'}(\bu')$ is free as a $\ZC'$-module of rank $r^n(2n-1)!!$.
Hence, $\W_\ZC(\dot\bu)$ is free as a $\ZC$-module of rank $r^n(2n-1)!!$.

Now suppose that $R$ is an arbitrary commutative ring (in which $2$ is
invertible). Then we can consider $R$ as a $\ZC$--algebra by letting
$\dot u_i$ act on $R$ as multiplication by~$u_i$, for $1\le i\le r$.
Since $\W_\ZC(\dot\bu)$ is $\ZC$--free, the $R$--algebra
$\W_\ZC(\dot\bu)\otimes_\ZC R$ is free as an $R$--module of rank
$r^n(2n-1)!!$. As the generators of $\W_\ZC(\dot\bu)\otimes_\ZC R$
satisfy the relations of $\W_{r,n}(\bu)=\W_R(\bu)$ we have a
surjective homomorphism
$\W_{r,n}(\bu)\longrightarrow\W_\ZC(\dot\bu)\otimes_\ZC R$.
By Proposition~\ref{spanning set} this map must be an isomorphism, so
we are done.
\end{proof}

As an easy application of the Theorem we obtain the following useful
fact which we will use many times below without mention.

\begin{Prop}\label{parabolic}
Suppose that $R$ is a commutative ring in which $2$ is invertible
and that $\Omega$ is $\bu$--admissible.
\begin{enumerate}
\item For $1\le m\le n$, let $\W_{r,m}'(\bu)$ be the subalgebra of
$\W_{r,n}(\bu)$
generated by $\set{S_i,E_i,X_j|1\le i<m\text{ and }1\le j\le m}$.
Then $\W_{r,m}'(\bu)\cong\W_{r,m}(\bu)$.
\item The Brauer algebra $\B_n(\omega_0)$ is isomorphic to the
subalgebra of $\W_{r,n}(\bu)$ generated by
$\set{S_i,E_i|1\le i<n}$.
\end{enumerate}
\end{Prop}

\section{The degenerate Hecke algebras of type $G(r, 1, n)$}
\label{cellular}

\def\s{\mathfrak s}
\def\ts{\tilde\s}
\def\t{\mathfrak t}
\def\u{\mathfrak u}
\def\v{\mathfrak v}

Suppose $R$ is a commutative ring and let $\bu\in R^r$. Recall from
section~2 that $\H_{r,n}(\bu)$ is the degenerate Hecke algebra with
parameters $\bu$. In this section we give several results from the
representation theory of $\H_{r,n}=\H_{r,n}(\bu)$ which we will need in our
study of the cyclotomic Nazarov--Wenzl algebras. As the proofs of these results
are very similar to (and easier than) the proofs of the corresponding
results for the Ariki--Koike algebras we are very brief with the
details.

The following result is proved by
Kleshchev~\cite{Klesh:book}. We use the seminormal representations of
$\W_{r,n}(\bu)$ to give another proof.

Let $\Lambda_r^+(n)$ be the set of $r$--multipartitions of $n$. We
consider $\Lambda_r^+(n)$ as a partially ordered set under dominance
$\gedom$, where $\lambda\gedom\mu$ if
$$\sum_{t=1}^{s-1}|\lambda^{(t)}|+\sum_{j=1}^k\lambda^{(s)}_k
     \ge\sum_{t=1}^{s-1}|\mu^{(t)}|+\sum_{j=1}^k\mu^{(s)}_k,$$
for $1\le s\le r$ and all $k\ge0$. If $\lambda\gedom\mu$ and
$\lambda\ne\mu$ we sometimes write $\lambda\gdom\mu$.

\begin{Theorem}\label{H-basis}
    The degenerate Hecke algebra $\H_{r,n}(\bu)$ is free as an
    $R$--module of rank $r^n n!$.
\end{Theorem}

\begin{proof}
It is not difficult to see that for any ring $R$ set
$$ \set{Y_1^{k_1}Y_2^{k_2}\cdots Y_n^{k_n}T_w|%
            0\le k_i\le r-1, w\in \Sym_{n}}$$
spans $\H_{r,n}(\bf u)$ as an $R$-module. So we need to prove that
these elements are linearly independent.

We adopt the notation from the proof of Theorem~\ref{W-basis}. As in
the proof of that result, we first consider the case when $R=\ZC$,
where $\ZC=\Z[\frac12,\dot u_1,\dots,\dot u_r]$, and we choose $r$
algebraically independent transcendental real numbers $u_1',\dots,u_r'$
which satisfy the hypotheses of Lemma~\ref{be real}. Let
$\ZC'=\Z[\frac12,u_1',\dots,u_r']$. Then
$\ZC\cong\ZC'\hookrightarrow\R$ and we can ask whether
the degenerate Hecke algebra $\H_\R(\bu')$, defined
over $\R$ and with parameters $\bu'=(u_1',\dots,u_r')$, acts on the
seminormal representations of $\W_\R(\bu')$.
By definition, if $\lambda\in\Lambda_r^+(n)$ then $E_i\Delta(\lambda)=0$,
for $1\le i<n$.  Therefore, over $\R$, $\Delta(\lambda)$ can be
considered as an $\H_\R(\bu')$--module by Corollary~\ref{AK irreds}.
Hence, as in the proof of Theorem~\ref{W-basis},
$$ \dim_\R \H_{\R}(\bu')
      \ge \sum_{\lambda\in \Lambda_r^+(n)} (\dim_\R\Delta(\lambda))^2
      =r^n n!.$$
Consequently, by the opening paragraph of the proof, this set is a
basis of $\H_\R(\bu')$.  As in the proof of Theorem~\ref{W-basis} it
follows that $\H_\ZC(\dot\bu)$ is free as a $\ZC$--module of rank
$r^nn!$. The result for a general ring $R$ now follows by a
specialization argument.
\end{proof}

We remark that the definition of the seminormal representations of
$\W_{r,n}(\bu)$ required that~$R$ satisfy assumption (\ref{be}). It is
not hard to modify the definition of the seminormal representations of
$\H_{r,n}(\bu)$ so that the formulae do not involve any square roots
and so that they work over an arbitrary field (cf.~\cite{AK}). In
particular, this leads to a simplification of the last argument.

\begin{Defn}[Graham and Lehrer~\cite{GL}]\label{GL}
    Let $R$ be a commutative ring and $A$ an $R$--algebra.
    Fix a partially ordered set $\Lambda=(\Lambda,\gedom)$ and for each
    $\lambda\in\Lambda$ let $T(\lambda)$ be a finite set. Finally,
    fix $C^\lambda_{\s\t}\in A$ for all
    $\lambda\in\Lambda$ and $\s,\t\in T(\lambda)$.

    Then the triple $(\Lambda,T,C)$ is a \textsf{cell datum} for $A$ if:
    \begin{enumerate}
    \item $\set{C^\lambda_{\s\t}|\lambda\in\Lambda\text{ and }\s,\t\in
        T(\lambda)}$ is an $R$--basis for $A$;
    \item the $R$--linear map $*\map AA$ determined by
        $(C^\lambda_{\s\t})^*=C^\lambda_{\t\s}$, for all
        $\lambda\in\Lambda$ and all $\s,\t\in T(\lambda)$ is an
        anti--isomorphism of $A$;
    \item for all $\lambda\in\Lambda$, $\s\in T(\lambda)$ and $a\in A$
        there exist scalars $r_{\s\u}(a)\in R$ such that
        $$aC^\lambda_{\s\t}
            =\sum_{\u\in T(\lambda)}r_{\s\u}(a)C^\lambda_{\u\t}
                     \pmod{A^{\gdom\lambda}},$$
            where
    $A^{\gdom\lambda}=R\text{--span}%
      \set{C^\mu_{\u\v}|\mu\gdom\lambda\text{ and }\u,\v\in T(\mu)}$.
    \end{enumerate}
    \noindent An algebra $A$ is a \textsf{cellular algebra} if it has
    a cell datum and in this case we call
    $\set{C^\lambda_{\s\t}|\s,\t\in T(\lambda), \lambda\in\Lambda}$
    a \textsf{cellular basis} of $A$.
\end{Defn}

To show that $\H_{r,n}(\bu)$ is a cellular algebra we modify the
construction of the Murphy basis of the Ariki--Koike
algebras; see~\cite{DJM:cyc}. For any multipartition
$\lambda=(\lambda^{(1)}, \lambda^{(2)}, \cdots, \lambda^{(r)})$ we
define
$u_{\lambda}=u_{a_1,1}u_{a_2,2}\cdots u_{a_{r-1},{r-1}}$,
where
$u_{a,i}=(Y_1-u_{i+1})(Y_2-u_{i+1})\cdots (Y_a-u_{i+1})$
and $a_i=\sum_{j=1}^i |\lambda^{(j)}|$, $1\le i\le r-1$.  Let
$\Sym_{\lambda}$ be the Young subgroup
$\Sym_{\lambda^{(1)}}\times \Sym_{\lambda^{(2)}}
                  \times\cdots\times \Sym_{\lambda^{(r)}}$
of $\Sym_n$.  Let
$x_{\lambda}=\sum_{w\in \Sym_{\lambda}}T_w$ and define
$$
    m_{\s\t}=T_{d(\s)^{-1}}u_{\lambda} x_{\lambda}T_{d(\t)}
                  \quad\in\H_{r,n}(\bu),
$$
where $\s, \t$ are standard $\lambda$-tableaux.

\begin{Theorem}\label{H cellular}
The set
$\set{m_{\s\t}|\s,\t\in \Std(\lambda)\text{ and }\lambda\in \Lambda_r^+(n)}$
is a cellular basis of~$\H_{r,n}(\bu)$.
\end{Theorem}

\begin{proof}
    The proof of this result is similar to, but much easier than, the
    corresponding result for the  cyclotomic  Hecke algebras.  See
    \cite{DJM:cyc} for details.
\end{proof}

We next give a formula for the Gram determinant of the cell modules of
$\H_{r,n}(\bu)$. This requires some definitions.

\begin{Defn}\label{H-generic}
    The parameters $\bu=(u_1,\dots,u_r)$ are \textbf{generic} for
    $\H_{r,n}(\bu)$ if whenever there exists $d\in\Z$ such that
    $u_i-u_j=d\cdot1_R$ then $|d|\ge n$.
\end{Defn}

The following Lemma
is well-known (cf.~\cite[Lemma~3.12]{JM:cyc-Schaper}), and is easily
verified by induction on~$n$.

\begin{Lemma}\label{generic u for H}
    Suppose that the parameters $\bu$ are generic for $\H_{r,n}(\bu)$
    and that~$R$ is a field with $\Char R>n$. Let $\lambda$ and
    $\mu$ be multipartitions of $n$ and suppose that
    $\s\in\Std_n(\lambda)$ and $\t\in\Std_n(\mu)$. Then $\s=\t$ if and
    only if $c_\s(k)=c_\t(k)$, for $k=1,\dots,n$.
\end{Lemma}

As in the definition of a cellular basis, if
$\lambda\in\Lambda_r^+(n)$ then we let $\H_{r,n}^{\gdom\lambda}$ be
the free $R$--submodule $\H_{r,n}(\bu)$ with basis
$\set{m_{\s\t}|\s,\t\in\Std(\mu)\text{ for }\mu\gdom\lambda}$. It
follows directly from Definition~\ref{GL}(c) that
$\H_{r,n}^{\gdom\lambda}$ is a two--sided ideal of $\H_{r,n}(\bu)$.

\begin{Lemma}\label{triangular}
Suppose that $\lambda$ is a multipartition of $n$ and
that $\s,\t\in\Std_n(\lambda)$. Then
$$Y_km_{\s\t}=c_\s(k) m_{\s\t}
+\sum_{\substack{\u\in\Std_n(\lambda)\\\u\gdom\s}}r_{\u\t}m_{\u\t}
                  \pmod{\H_{r,n}^{\gdom\lambda}},$$
for some $r_{\u\t}\in R$.
\end{Lemma}

\begin{proof} If $r=1$ then this is a result of Murphy's~\cite{M:Nak}.
The general case can be deduced from this following the argument of
\cite[Prop.~3.7]{JM:cyc-Schaper}.
\end{proof}

We can now follow the arguments of \cite{M:gendeg} to construct a
``seminormal'' basis of $\H_{r,n}(\bu)$.

\begin{Defn}Suppose that $\lambda\in\Lambda_r^+(n)$.
\begin{enumerate}
\item For each $\t\in\Std(\lambda)$ let
$$F_\t=\prod_{k=1}^n\prod_{\substack{\mu\in\Lambda_r^+(n)\\\u\in\Std(\mu)\\
                           c_\u(k)\ne c_\t(k)}}
                           \frac{Y_k-c_\u(k)}{c_\t(k)-c_\u(k)}.$$
\item If $\s,\t\in\Std(\lambda)$ then let
           $f_{\s\t}=F_\s m_{\s\t} F_\t$.
\end{enumerate}
\end{Defn}

Using the last two results and the definitions it is not hard to show
that if~$\s,\t$ and $\u$ are standard tableaux then
$f_{\s\t}F_\u=\delta_{\t\u}f_{\s\t}$; see, for example,
\cite[Prop.~3.35]{M:ULect}. Hence, from Theorem~\ref{H cellular} and
Lemma~\ref{triangular} we obtain the following.

\begin{Prop}\label{orthog basis}
Suppose that $R$ is a field with $\Char R>n$ and that $\bu$ is
generic for $\H_{r,n}(\bu)$. Then $\set{f_{\s\t}|\s,\t\in
\Std(\lambda),\lambda\in \Lambda_r^+(n)}$ is a basis of
$\H_{r,n}(\bu)$. Moreover, for each standard tableau $\t$ there exists
a scalar $\gamma_\t\in R$ such that
$$f_{\s\t}f_{\u\v}= \delta_{\t\u}\gamma_\t f_{\s\v},$$
where $\s,\t\in\Std(\lambda)$, $\u,\v\in\Std(\mu)$,
and $\lambda,\mu\in\Lambda_r^+(n)$.
\end{Prop}

Notice, in particular, that the Proposition implies that
$\{f_{\s\t}\}$ is also a cellular basis of $\H_{r,n}(\bu)$.

Although we will not pursue this here, we remark that
$F_\t=\frac1{\gamma_\t}f_{\t\t}$ and that these elements give a
complete set of pairwise orthogonal primitive idempotents for
$\H_{r,n}(\bu)$. This can be proved by repeating the argument of
\cite[Theorem~2.15]{M:gendeg}

Suppose that $\lambda$ is a multipartition of $n$ and let $S(\lambda)$
be the associated \textsf{Specht module}, or cell module, of
$\H_{r,n}(\bu)$. Thus, $S(\lambda)$ is the free $R$--module with basis
$\set{m_\s|\s\in\Std(\lambda)}$ and where the action of
$\H_{r,n}(\bu)$ on $S(\lambda)$ is given by
$$a m_\s=\sum_{\u\in\Std(\lambda)}r_{\s\u}(a)m_\u,$$
where the scalars $r_{\s\u}(a)\in R$ are as in
Definition~\ref{GL}(c).

It follows directly from Definition~\ref{GL} that $S(\lambda)$ comes
equipped with a symmetric bilinear form $\langle\ ,\ \rangle$ which is
determined by
$$\langle m_\s,m_\t\rangle m_{\u\v}\equiv
        m_{\u\s}m_{\t\v}\pmod{\H_{r,n}^{\gdom\lambda}},$$
for $\s,\t,\u,\v\in\Std(\lambda)$. Let
$G(\lambda)
    =\det\(\langle m_\s,m_\t\rangle\)$,
for $\s,\t\in\Std(\lambda)$, be the Gram determinant of this form. So
$\G(\lambda)$ is well--defined up to a unit in~$R$.

\begin{Cor}Suppose that $R$ is a field with $\Char R>n$ and that
$\bu$ is generic for $\H_{r,n}(\bu)$. Let $\lambda$ be a
multipartition of $n$. Then
$$\G(\lambda)=\prod_{\t\in\Std(\lambda)}\gamma_\t.$$
\end{Cor}

\begin{proof} Fix $\t\in\Std(\lambda)$. Then Specht module
$S(\lambda)$ is isomorphic to the submodule of
$\H_{r,n}/\H_{r,n}^{\gdom\lambda}$ which is spanned by
$\set{m_{\s\t}+\H_{r,n}^{\gdom\lambda}|\s\in\Std(\lambda)}$, where
the isomorphism is given by
$\theta\map{\H_{r,n}/\H_{r,n}^{\gdom\lambda}}S(\lambda);
          m_{\s\t}+\H_{r,n}^{\gdom\lambda}\mapsto m_\s$.
Let $f_\s=\theta(f_{\s\t})$. Then $\set{f_\s|\s\in\Std(\lambda)}$ is a
basis of $S(\lambda)$ and the transition matrix between the two bases
$\{m_\s\}$ and $\{f_\s\}$ of $S(\lambda)$ is unitriangular by
Lemma~\ref{triangular}. Consequently,
$\G(\lambda)=\det\(\langle f_\s,f_\t\rangle\)$, where
$\s,\t\in\Std(\lambda)$. However, it follows from the multiplication
formulae in Proposition~\ref{orthog basis} that
$\langle f_\s,f_\t\rangle=\delta_{\s\t}\gamma_\t$; see the proof of
\cite[Theorem~2.11]{M:gendeg} for details. Hence the result.
\end{proof}

Consequently, in order to compute $\G(\lambda)$ it is sufficient to
determine $\gamma_t$, for all $\t\in\Std(\lambda)$. It is possible to
give an explicit closed formula for $\gamma_\t$
(cf.~\cite[(2.8)]{M:gendeg}), however, the following recurrence
relation is easier to check and sufficient for our purposes.

Given two standard $\lambda$--tableaux $\s$ and $\t$ write
$\s\gedom\t$ if $\s_k\gedom\t_k$, for $1\le k\le n$. Let $\t^\lambda$
be the unique standard $\lambda$--tableaux such that
$\t^\lambda\gedom\s$ for all $\s\in\Std(\lambda)$. If $\s\gedom\t$ and
$\s\ne\t$ then we write $\s\gdom\t$.

\begin{Lemma}\label{gamma recurrence}
Suppose that $R$ is generic for $\H_{r,n}(\bu)$ and that
$\Char R>n$. Let $\lambda$ be a multipartition of $n$.
\begin{enumerate}
\item $\displaystyle\gamma_{\t^\lambda}
   =\prod_{1\le t\le r}\prod_{i\ge1}(\lambda^{(t)}_i)!
    \cdot\prod_{1\le s<t\le r}
       \prod_{\substack{i,j\ge 1\\1\le j\le\lambda^{(s)}_i}}
          (j-i+u_s-u_t).$
\item Suppose that $\s,\t\in\Std(\lambda)$ such that $\s\gdom\t$
and $\s=S_k\t$, for some $k$. Then
$\gamma_\t=\frac{(c_\s(k)-c_\t(k)+1)(c_\s(k)-c_\t(k)-1)}%
                 {(c_\s(k)-c_\t(k))^2}\gamma_\s$.
\end{enumerate}
\end{Lemma}

\begin{proof} Part (a) follows easily by induction on $n$. Part~(b)
follows using arguments similar to
\cite[Cor.~3.14 and Prop.~3.19]{JM:cyc-Schaper}
\end{proof}

We remark that the arguments of ~\cite[3.30-3.37]{JM:cyc-Schaper} can
now be adapted to give a closed formula for $\G(\lambda)$. The final
result is that
$$\G(\lambda)
    =\prod_{\nu\in\Lambda_r^+(n)}g_{\lambda\nu}^{|\Std(\lambda)|},$$
where $g_{\lambda\nu}$ is a product of terms of the form
$\(c_{\t^\lambda}(k)-c_{\t^\nu}(l)\)^{\pm1}$, where these terms are
determined in exactly the same way as in
\cite[Defn~3.36]{JM:cyc-Schaper}. As we do not need the precise
formula we leave the details to the interested reader.

\begin{Theorem}\label{H semisimple}
Suppose that $R$ is a field and that $\bu\in R^r$. Then
$\H_{r,n}(\bu)$ is (split) semisimple if and only if $\Char R>n$ and
$\bu$ is generic for $\H_{r,n}(\bu)$.
\end{Theorem}

\begin{proof} First, note that because $\H_{r,n}(\bu)$ is cellular, it
is semisimple if and only if it is split semisimple; see, for example,
\cite[Cor~2.21]{M:ULect}.

Next, suppose that $\Char R>n$ and that $\bu$ is generic for
$\H_{r,n}(\bu)$. Then $\G(\lambda)\ne0$ for all
$\lambda\in\Lambda_r^+(n)$ by Lemma~\ref{gamma recurrence}.
Consequently, for each $\lambda\in\Lambda_r^+(n)$ the Specht module
$S(\lambda)$ is irreducible. Hence, by \cite[Cor~2.21]{M:ULect} again,
$\H_{r,n}(\bu)$ is semisimple.

To prove the converse, let
$\lambda=((n),(0),\dots,(0))\in\Lambda_r^+(n)$ and set
$m_\lambda=m_{\t^\lambda\t^\lambda}$; more explicitly,
$$m_\lambda=\sum_{w\in\Sym_n}T_w
                \cdot\prod_{t=2}^r\prod_{k=1}^n(Y_k-u_t).$$
It is easy to see that $T_\sigma m_\lambda=m_\lambda=m_\lambda T_\sigma$,
for any $\sigma\in\Sym_n$. It also follows from Lemma~\ref{triangular} that
$Y_km_\lambda=c_{\t^\lambda}(k)m_\lambda=m_\lambda Y_k$. Hence,
$\H_{r,n}(\bu)m_\lambda\H_{r,n}(\bu)=Rm_\lambda$ and
$$m_\lambda^2
     =n!\prod_{t=2}^r\prod_{d=0}^{n-1}(u_1+d-u_t)\cdot m_\lambda.$$
If $\Char R\le n$ then $n!=0$ in $R$ so that $Rm_\lambda$ is a
nilpotent ideal in $\H_{r,n}(\bu)$, so $\H_{r,n}(\bu)$ is not
semisimple. On the other hand if $\bu$ is not generic for
$\H_{r,n}(\bu)$ then $u_i-u_j=d1_R$ for some $i\ne j$ and some
$d\in\Z$ with $|d|<n$. By renumbering $u_1,\dots,u_r$, if necessary,
we see that $Rm_\lambda$ is a nilpotent ideal. Hence, if either
$\Char R\le n$, or if $\bu$ is not generic for $\H_{r,n}(\bu)$, then
$\H_{r,n}(\bu)$ is not semisimple.
\end{proof}

\section{A cellular basis of $\W_{r,n}(\bu)$}

Throughout this section we assume that $R$ is a commutative ring in
which $2$ is invertible and that $\Omega$ is $\bu$--admissible.  This
section constructs a cellular basis for~$\W_{r,n}=\W_{r,n}(\bu)$ using
the cellular bases of the algebras $\H_{r,n-2f}=\H_{r,n-2f}(\bu)$, for $0\le
f\le\floor{n2}$, together with a series of filtrations of~$\W_{r,n}$.
Our construction of a cellular basis of~$\W_{r,n}$ is modelled on
Enyang's work~\cite{Enyang} for the Brauer and BMW algebras.

Before we begin we need to fix some notation. Recall that the set
$\{S_1,\dots,S_{n-1}\}$ generates a subalgebra of $\W_{r,n}$ which is
isomorphic to the group ring of $\Sym_n$. For each permutation
$w\in\Sym_n$ we defined the corresponding braid diagram $\gamma(w)$ in
section~2; we now set $S_w=B_{\gamma(w)}$. Equivalently, if
$w=(i_1,i_1+1)\dots(i_k,i_k+1)$, where $1\le i_j<n$ for all $j$, then
$S_w=S_{i_1}\dots S_{i_k}$. Then $\set{S_w|w\in\Sym_n}$ is a basis for
the subalgebra of $\W_{r,n}$ generated by $\{S_1,\dots,S_{n-1}\}$.

Next, suppose that $f$ is an integer with $0\le f\le\floor{n2}$. It
follows from Theorem~\ref{H-basis} that we can identity $\H_{r, n-2f}$
with the subalgebra of $\H_{r,n}(\bu)$ generated by $Y_i$ and $T_j$, where
$1\le i\le n-2f$ and $1\le j\le n-2f-1$. Similarly, by
Proposition~\ref{parabolic}, we can identify $\W_{r, n-2f}$ with the
subalgebra of $\W_{r,n}$ generated by $X_i$, $S_j$ and $E_j$, where
$1\le i\le n-2f$ and $1\le j\le n-2f-1$.

\begin{Defn}\label{epsilon f} 
Suppose $0\le f<\floor{n2}$.  Let $\Ef=\W_{r,n-2f}E_1\W_{r,n-2f}$ be
the two-sided ideal of $\W_{r, n-2f}$ generated by~$E_1$.
\end{Defn}

\begin{Prop}\label{degen isom} Suppose that $0\le f<\floor{n2}$. Then
there is a unique $R$--algebra isomorphism
$\proj_f:\H_{r,n-2f}(\bu)\cong\W_{r,n-2f}/\Ef$ such that
$$\proj_f(T_i)=S_i+\Ef \text{ and } \proj_f(Y_j)=X_j+\Ef,$$
for $1\le i<n-2f$ and $1\le j\le n-2f$.
\end{Prop}

\begin{proof}We first show that $\W_{r,n-2f}/\Ef$ is a free
$R$--module of rank $r^{n-2f}(n-2f)!$. It follows from the
multiplication formulae for Brauer diagrams that an $r$-regular
monomial $X^\alpha B_{\gamma} X^\beta$ in $\W_{r,n-2f}$ belongs to
$\Ef$ whenever $\gamma$ has a horizontal edge (equivalently,
$\gamma\ne\gamma(w)$ for some $w\in\Sym_{n-2f}$). If
$\gamma=\gamma(w)$, for some $w\in\Sym_{n-2f}$, then $B_\gamma=S_w$
and $\gamma$ contains no horizonal edges, so the definition of
regularity (Definition~\ref{regular}), forces $\beta=0$. So, by
Theorem~\ref{W-basis}, $\W_{r,n-2f}/\Ef$ is spanned by the elements
$\set{X^\alpha S_w +\Ef|0\le \alpha_i<r, \text{ for }1\le i\le n-2f,
              \text{ and }w\in\Sym_{n-2f}}$.
Note that this set contains $r^{n-2f}(n-2f)!$ elements.

To see that the elements at the end of the last paragraph are linearly
independent we use the seminormal representations from section~4.
Using the arguments and the notation from the proof of
Theorem~\ref{W-basis}, it is enough to show that
$\dim_\R\W_\R(\bu')/\Ef\ge r^{n-2f}(n-2f)!$. Now a seminormal
representation $\Delta(\lambda)$ of~$\W_\R(\bu')$ is a representation
of $\W_\R(\bu')/\Ef$ if and only if $\Ef\Delta(\lambda)=0$, which
happens if and only if $\lambda$ is a multipartition of $n-2f$.
Therefore, by the arguments of section~5, $\dim_\R\W_\R(\bu')/\Ef\ge
r^{n-2f}(n-2f)!$. Hence, by the arguments used in the proof of
Theorem~\ref{W-basis} (compare, Theorem~\ref{H-basis}), the elements
above are a basis of $\W_{r,n-2f}/\Ef$ and, consequently,
$\W_{r,n-2f}/\Ef$ is free as an $R$--module of rank $r^{n-2f}(n-2f)!$
as claimed.

Inspecting the relations of $\H_{r,n-2f}(\bu)$ and $\W_{r,n-2f}$ shows that
there is a unique algebra homomorphism
$\proj_f\map{\H_{r,n-2f}(\bu)}\W_{r,n-2f}/\Ef$ such that
$\proj_f(T_i)=S_i+\Ef$ and $\proj_f(Y_j)=X_j+\Ef$. To see that
$\proj_f$ is an isomorphism observe that $\proj_f$ maps the
basis of $\H_{r,n-2f}(\bu)$ to the basis of $\W_{r,n-2f}/\Ef$.
Hence, it is an isomorphism with inverse determined by
$\proj_f^{-1}(X^\alpha S_w+\Ef)=Y^\alpha T_w$, for
$w\in\Sym_{n-2f}$ and $0\le\alpha_i<r$ where $1\le i\le n-2f$.
\end{proof}

\begin{Defn}\label{Ef}
Let $E^f=E_{n-1}E_{n-3}\cdots E_{n-2f+1}$ and let
$\W_{r,n}^f=\W_{r,n}E^f\W_{r,n}$ be the two-sided ideal of $\W_{r,n}$
generated by $E^f$.
If $f=\floor{n2}$ then we set $\H_{r,n-2f}(\bu)=R$ and $\W_{r,n}^{f+1}=0$.
\end{Defn}
Note that this gives a filtration of $\W_{r,n}$ by two--sided ideals:
$$\W_{r,n}=\W_{r,n}^0\supset\W_{r,n}^1\supset\dotsi
     \supset\W_{r,n}^{\floor{n2}}\supset\W_{r,n}^{\floor{n2}+1}=0.$$
For $0\le f\le\floor{n2}$ let
$\pi_f\map{\W_{r,n}^f}\W_{r,n}^f/\W_{r,n}^{f+1}$ be the corresponding
projection map of $\W_{r,n}$--bimodules.

For convenience we set $\N_r=\{0,1,\dots,r-1\}$ and define $\Nrf$ to be
the set of $n$--tuples $\kappa=(k_1,\dots,k_n)$ such that
$k_i\in\N_r$ and $k_i\ne0$ only for $i=n-1,n-3,\dots,n-2f+1$.  Thus,
if $\kappa\in \Nrf$ then
$X^\kappa=X_{n-1}^{k_{n-1}}X_{n-3}^{k_{n-3}}\dots X_{n-2f+1}^{k_{n-2f+1}}
   \in\W_{r,n}$.

\begin{Lemma}\label{if} Suppose that $0\le f\le\floor{n2}$ and
$\kappa\in \Nrf$. Then $E^f X^\kappa\Ef\subset\W_{r,n}^{f+1}$.
\end{Lemma}

\begin{proof}As $E^{f+1}=E^fE_{n-2f-1}$, this follows
because $\Ef=\W_{r,n-2f}E_{n-2f-1}\W_{r,n-2f}$ and every
element of $\W_{r,n-2f}$ commutes with $E^fX^\kappa$.
\end{proof}

Combining the last two results we have a well-defined $R$-module
homomorphism $\sigma_f\map{\H_{r, n-2f}}\W^f_{r,n}/\W_{r,n}^{f+1}$,
for each integer $f$, with $0\le f\le\floor{n2}$, given by
$$\sigma_f(h)=E^f\proj_f(h)+\W_{r,n}^{f+1},$$
for $h\in \H_{r, n-2f}$.\label{sigma f}

We will need the following subgroups in order to understand the ideals
$\W_{r,n}^f$.

\begin{Defn}Suppose that $0\le f\le \floor{n2}$. Let $\Bcal_f$ be the
subgroup of $\Sym_n$ generated by
$\{S_{n-1}, S_{n-2}S_{n-1}S_{n-3}S_{n-2}, \cdots, S_{n-2f+2}
         S_{n-2f+1}S_{n-2f+3} S_{n-2f+2}\}$.
\end{Defn}

The symmetric group $\Sym_n$ acts on the set of Brauer diagrams $\BDiag(n)$
from the right. Let $\gamma=\gamma_{n-1}\circ\dots\circ\gamma_{n-2f+1}$. Then
$E^f=B_\gamma$ and $\Bcal_f$ is the stablizer in $\Sym_n$ of the diagram
$\gamma$. The group $\Bcal_f$ is isomorphic to the hyperoctahedral group
$\Z/2\Z\wr\Sym_f$, a Coxeter group of type~$\Bcal_f$.

Given an integer $f$, with $0\le f\le \floor{n2}$, let
$\tau=((n-2f),(2^f)\)$ and define
$$\Dcal_f=\Set[80]d\in \Sym_n|
      $\t^{\tau} d=(\t_1, \t_2)$ is a row standard $\tau$-tableau
     and the first column of $\t_2$ is increasing from top to bottom|.$$

The following result is equivalent to \cite[Prop.~3.1]{Enyang}.
(Enyang considers a subgroup of $\Sym_n$ which is conjugate to $\Bcal_f$.)

\begin{Lemma}\label{cosets}
Suppose that $0\le f\le\floor{n2}$. Then $\Dcal_f$ is a complete set of
right coset representatives for $\Sym_{n-2f}\times \Bcal_f$ in $\Sym_n$.
\end{Lemma}

The point of introducing the subgroup $\Bcal_f$ is the following.

\begin{Lemma}\label{E^f}
Suppose that $0\le f\le\floor{n2}$, $w\in\Sym_{n-2f}$ and that
$b\in \Bcal_f$. Then $S_wE^f=E^fS_w$ and $E^fS_b=E^f=S_bE^f$.
\end{Lemma}

\begin{proof}
The first claim is obvious by (\ref{Waff relations})(d)(i). For the second claim it
is enough to consider the case when $b$ is a generator of $\Bcal_f$. In this case
the claim is easily checked using the tangle relations and the untwisting relations.
\end{proof}

Motivated by the definition of the elements
$m_{\s\t}\in\H_{r,n-2f}(\bu)$
from the previous section, and by the work of Enyang \cite{Enyang}, we
make the following definition.

\begin{Defn}\label{M_st} Suppose that $0\le f\le\floor{n2}$ and
$\lambda\in\Lambda_r^+(n-2f)$. Then for each pair $(\s,\t)$ of
standard $\lambda$--tableaux define
$$M_{\s\t}=S_{d(\s)^{-1}}\cdot \prod_{s=2}^r
      \prod_{i=1}^{|\lambda^{(1)}|+\dots+|\lambda^{(s-1)}|}(X_i-u_s)
        \sum_{w\in\Sym_\lambda}S_w \cdot S_{d(\t)}.$$
\end{Defn}

We remark that we will not ever really use this explicit formula for
the elements $M_{\s\t}$. In what follows all that we need is a
family of elements $\{M_{\s\t}\}$ in $\W_{r,n}$ which are related to
some cellular basis of $\H_{r,n-2f}(\bu)$ as in
Lemma~\ref{M_st properties}(d) below.

The following result follows easily using the relations of $\W_{r,n}$ and the
definitions.

\begin{Lemma}\label{M_st properties}
Suppose that $0\le f\le\floor{n2}$, $\lambda\in\Lambda_r^+(n-2f)$ and
that $\s,\t\in\Std(\lambda)$. Then:
\begin{enumerate}
\item $E^fM_{\s\t}=M_{\s\t}E^f\in\W_{r,n}^f$.
\item If $\kappa\in \Nrf$ then $M_{\s\t}X^\kappa=X^\kappa M_{\s\t}$.
\item If $w$ is a permutation of $\{n-2f+1,\dots,n\}$ then
$M_{\s\t}S_w=S_wM_{\s\t}$. In particular, $M_{\s\t}S_w=S_wM_{\s\t}$ if
$w\in \Bcal_f$.
\item We have $\sigma_f(m_{\s\t})=\pi_f(E^fM_{\s\t})$.
\end{enumerate}
\end{Lemma}

The filtration of $\W_{r,n}$ given by the ideals $\W_{r,n}^f$ is still
too coarse to be cellular.

\begin{Defn}\label{Wlambda}
Suppose that $\lambda$ is a multipartition of $n-2f$, where
$0\le f\le\floor{n2}$. Define $\Wlambda$ to be
the two--sided ideal of $\W_{r,n}$ generated by $\W_{r,n}^{f+1}$ and
the elements
$$\set{E^fM_{\s\t}|\s,\t\in \Std(\mu) \text{ and }
    \mu\in\Lambda_r^+(n-2f)\text{ with }\mu\trianglerighteq\lambda}.$$
We also set $\Wlam =\sum_{\mu\gdom\lambda}\W_{r,n}^{\gedom\mu}$, where
in the sum $\mu\in\Lambda_r^+(n-2f)$.
\end{Defn}

Observe that
$$\W_{r,n}^{f+1}\subseteq\Wlam\subset\Wlambda\subseteq\W_{r,n}^f$$
and that $\Wlambda\subset\W_{r,n}^{\gdom\mu}$ whenever
$\lambda\gdom\mu$. Consequently, the ideals $\{\Wlambda\}$
give a refinement of the filtration of $\W_{r,n}$ by the ideals
$\{\W_{r,n}^f\}$.

\begin{Defn} Suppose that $\s\in \Std(\lambda)$. We define
$\Delta_\s(f,\lambda)$ to be the $R$-submodule of
$\Wlambda/\Wlam$  spanned by the elements
$$\set{E^fM_{\s\t}X^\kappa S_d+\Wlam|(\t,\kappa,d)\in\delta(f,\lambda)},$$
where $\delta(f,\lambda)
  =\set{(\t,\kappa,d)|\t\in\Std(\lambda),\kappa\in \Nrf \text{ and }
                        d\in \Dcal_f}$.
\end{Defn}

We will see below that $\Delta_\s(f,\lambda)$ is a right $\W_{r,n}$--module and
that the spanning set in the definition is a basis of
$\Delta_\s(f,\lambda)$. Moreover, there is a natural isomorphism
$\Delta_\s(f,\lambda)\cong\Delta_\t(f,\lambda)$,
whenever $\s,\t\in\Std(\lambda)$.

Before we begin studying the modules $\Delta_\s(f,\lambda)$ it is
convenient to define a degree function on $\W_{r,n}$.  Recall from
Theorem~\ref{W-basis} that the set of $r$--regular monomials is a
basis of $\W_{r,n}$.

\begin{Defn}\label{degree}
   Suppose that $a=\sum r_{\alpha\gamma\beta}X^\alpha B_\gamma
   X^\beta\in\W_{r,n}$, where each of the monomials in the sum is
   $r$--regular.  Then the \textsf{degree} of $a$ is the
   integer
   $$\deg a=\max\seT{\sum_{i=1}^n(\alpha_i+\beta_i)|%
          r_{\alpha\gamma\beta}\ne0\text{ for some } \gamma\in\BDiag(n)}.$$
\end{Defn}

In particular, $\deg S_i=\deg E_i=0$, for $1\le i<n$, and $\deg X_j=1$,
for $1\le j\le n$. We note that the proof of
\cite[Lemma~4.4]{Nazarov:brauer} implies that
$$\deg(ab)\le\deg(a)+\deg(b),\quad \text{for all } a,b\in\W_{r,n}.$$

\begin{Lemma}\label{small}
Suppose that $1\le j<n$ and that $1\le k<r$. Then
$E_jX_j^kE_j=E_j\omega_j^{(k)}$, where $\omega_j^{(k)}$ is a central
element in $\W_{r,j-1}$ with $\deg\omega_j^{(k)}<k$.
\end{Lemma}

\begin{proof} We argue by induction on $j$. If $j=1$ then $\deg\omega_1^{(k)}=0$
because $\omega_1^{(k)}\in R$ by relation (\ref{Waff relations})(f). Suppose then
that $j>1$.

By Lemma~\ref{tilde W} $\omega_j^{(k)}$ is a central
element of~$\W_{r,j-1}$ in~$R[X_1,\dots,X_{j-1}]$ and $\deg\omega^{(k)}_j\le k$.
Consequently, if $\omega_j^{(k)}=\sum_\alpha r_\alpha X^\alpha$, for some
$r_\alpha\in R$, then $\omega_j^{(k)}E_j=\sum_\alpha r_\alpha X^\alpha E_j$ where
each of the monomials $X^\alpha E_j$ is $r$--regular. Hence,
$\deg(\omega_j^{(k)}E_j)=\deg\omega_j^{(k)}$.
Therefore, it is enough to prove that $\deg(\omega_j^{(k)}E_j)<k$.
By Lemma~\ref{SX^a},
\begin{align*}
\omega_j^{(k)}E_j &=E_{j} X_j^k E_j
                   =(-1)^k E_j X_{j+1}^k E_j
                   =(-1)^k E_jS_{j-1}X_{j+1}^k S_{j-1} E_j\\
                  &=(-1)^k E_jE_{j-1}S_{j} X_{j+1}^k S_{j}E_{j-1}E_j\\
                  &=(-1)^k E_jE_{j-1}(X_j^k+X)E_{j-1}E_j,
\end{align*}
where $X\in\W_{r,j+1}$ and $\deg X<k$ since $\deg(ab)\le \deg(a)+\deg(b)$. We have that
$\deg(E_jE_{j-1}XE_{j-1}E_j)\le \deg X<k$ and that
\begin{align*}
E_jE_{j-1} X_{j}^k E_{j-1}E_j
  &=(-1)^k E_jE_{j-1} X_{j-1}^k E_{j-1}E_j
   =(-1)^k\omega_{j-1}^{(k)} E_jE_{j-1}E_j\\
  &= (-1)^k\omega_{j-1}^{(k)}E_j.
\end{align*}
By induction $\deg\omega_{j-1}^{(k)}<k$, so this completes the proof
of the Lemma.
\end{proof}

Given integers $j$ and $k$, with $1\le j,k\le n$, let
$E_{j,k}=B_\gamma$ where $\gamma$ is the Brauer diagram with
horizontal edges $\{j,k\}$ and $\{\bar j,\bar k\}$, and with all other
edges being vertical. Thus,
$S_wE_iS_{w^{-1}}=E_{(i)w^{-1},(i+1)w^{-1}}$, for all $w\in\Sym_n$.
Finally, note that $E_i=E_{i,i+1}$.

Until further notice we fix an integer $f$, with
$0\le f\le\floor{n2}$, a multipartition $\lambda\in\Lambda_r^+(n-2f)$
and a standard $\lambda$--tableau $\s$ and consider
$\Delta(f,\lambda)=\Delta_\s(f,\lambda)$.  The next two Lemmas show
that $\Delta(f,\lambda)$ is a right $\W_{r,n}$--submodule of
$\Wlambda/\Wlam$ and that the action of $\W_{r,n}$ on
$\Delta(f,\lambda)$ does not depend on $\s$.

If $\kappa=(\kappa_1,\dots,\kappa_n)\in \Nrf$ we set
$|\kappa|=\kappa_{n-1}+\kappa_{n-3}+\dots+\kappa_{n-2f+1}
         =\deg X^\kappa$.

\begin{Lemma}\label{|k|=0}
Suppose that $\t\in\Std(\lambda)$ and $d\in \Dcal_f$. For $1\le i<n$ and
$1\le j\le n$ there exist scalars $a_{\v e},b_{\v e},c_{\v\rho e}\in R$, which do
not depend on $\s$, such that:
\begin{enumerate}
\item $\displaystyle E^fM_{\s\t}S_d\cdot S_i\equiv
  \sum_{\substack{\v\in\Std(\lambda)\\e\in\Dcal_f}}
              a_{\v e}E^fM_{\s\v} S_e \pmod{\Wlam}$,
\item $\displaystyle E^fM_{\s\t}S_d\cdot E_i\equiv
  \sum_{\substack{\v\in\Std(\lambda)\\e\in\Dcal_f}}
              b_{\v e}E^fM_{\s\v} S_e \pmod{\Wlam}$,
\item $\displaystyle E^fM_{\s\t}S_d\cdot X_j\equiv
  \sum_{\substack{(\v,\rho,e)\in\delta(f,\lambda)\\|\rho|\le1}}
              c_{\v\rho e}E^fM_{\s\v}X^\rho S_e \pmod{\Wlam}$.
\end{enumerate}
\end{Lemma}

\begin{proof} (a) Now, $S_dS_i=S_{d(i,i+1)}$ and by Lemma~\ref{cosets} we can write
$d(i,i+1)=abe$ where $a\in \Sym_{n-2f}$, $b\in \Bcal_f$ and $e\in \Dcal_f$; so
$S_dS_i=S_aS_bS_e$. By part (d) of Lemma~\ref{M_st properties},
respectively, we have
$$E^fM_{\s\t} S_a\equiv E^f\proj_f(m_{\s\t}) S_a
              \equiv E^f\proj_f(m_{\s\t}T_a)\pmod{\Wlam},$$
since $\W_{r,n}^{f+1}\subseteq\Wlam$.  As $m_{\s\t}$ is a cellular
basis element for $\H_{r, n-2f}$, we can write $m_{\s\t}T_a$ as a
linear combination of terms $m_{\s\v}$ plus an element of
$\H_{r,n}^{\gdom\lambda}$. Consequently,
$(E^fM_{\s\t}+\Wlam)S_a$ can be written
in the desired form. Hence,  we may now assume that $a=1$.

To complete this case, observe that if $\v\in\Std(\lambda)$ then,
by Lemma~\ref{M_st properties}(c) and Lemma~\ref{E^f},
$E^fM_{\s\v}S_bS_e=E^fS_bM_{\s\v}S_e=E^fM_{\s\v}S_e$ as required.

\medskip\noindent(b) We have to consider the product $E^fM_{\s\t}S_d E_i$. Let
$j=(i)d^{-1}$ and $k=(i+1)d^{-1}$. Then $S_dE_i=E_{j,k}S_d$ so that
$E^fM_{\s\t}S_dE_i=E^fM_{\s\t}E_{j,k}S_d$. By part(a) we may
assume that $d=1$. We can also assume that $j<k$ since
$E_{j,k}=E_{k,j}$. So we need to show that $E^fM_{\s\t}E_{j,k}+\Wlam$
has the required form. There are three cases to consider.

(1) First, suppose that $j<k\le n-2f$. Then
$E_{j,k}\in\W_{r,n-2f}$, so that $M_{\s\t}E_{j,k}\in\Ef$ and
$E^fM_{\s\t}E_{j,k}\in E^f\Ef\subseteq\W_{r,n}^{f+1}$ by
Lemma~\ref{if}. Hence,
$E^fM_{\s\t}S_dE_i\in\W_{r,n}^{f+1}\subseteq\Wlam$ and part~(b) is
true when $j<k\le n-2f$.

(2) Next, suppose that $j\le n-2f<k$. An easy exercise in
multiplying Brauer diagrams shows that
$$
E^fE_{j,k}=\begin{cases}
        E^fS_{(j, k-1)}, &\text{ if $n-k$ is even,}\\
        E^fS_{(j, k+1)}, &\text{ if $n-k$ is odd.}
\end{cases}
$$
So Lemma~\ref{M_st properties}(a) implies that
$E^fM_{\s\t}S_d E_i=M_{\s\t}E^fE_{j,k}S_d=E^fM_{\s\t}S_dS_{d^{-1}(j,k\pm1)d}$,
we again deduce the result from part~(a).

(3) Finally, suppose that $n-2f<j<k$. Then
$M_{\s\t}E_{j,k}=E_{j,k}M_{\s\t}$ and a Brauer diagram calculation
shows that $E^fE_{j,k}=E^fS_w$, where $w$ is a permutation of
$\{n-2f+1,\dots,n\}$. Consequently,
$$E^fM_{\s\t}S_dE_i=E^fM_{\s\t}E_{j,k}S_d=E^fE_{j,k}M_{\s\t}S_d
                   =E^fS_wM_{\s\t}S_d=E^fM_{\s\t}S_wS_d,$$
where the last equality follows from Lemma~\ref{M_st properties}(c).
As $S_wS_d=S_dS_{d^{-1}wd}$ we are done by part~(a).

\medskip\noindent(c) It follows from the skein relations that 
$S_dX_j=X_{(j)d}S_d+B$, for some
$B\in\B_n(\omega_0)$. Hence, by parts~(a) and~(b) it suffices to show that
$E^fM_{\s\t}X_i$ can be written in the required form, for
$1\le i\le n$. If $i\le n-2f$ then
$$E^fM_{\s\t}X_i+\Wlam=E^f\epsilon_f(m_{\s\t})X_i+\Wlam
                      =E^f\epsilon_f(m_{\s\t}Y_i)+\Wlam,$$
so the result follows because $m_{\s\t}$ is a cellular basis element
of $\H_{r,n-2f}(\bu)$. If $i>n-2f$ then the result is immediate if $n-i$ is
odd. If $n-i$ is even then $i-1>n-2f$, so the result follows because
$E_{i-1}X_i=-E_{i-1}X_{i-1}$ by (\ref{Waff relations})(h).

This completes the proof of the Lemma.
\end{proof}

\begin{Prop}\label{|k|>=0}
Suppose that $(\t,\kappa,d)\in\delta(f,\lambda)$. For $1\le i<n$ and
$1\le j\le n$ there exist scalars $a_{\v\rho e},b_{\v\rho e},c_{\v\rho e}\in R$,
which do not depend on $\s$, such that:
\begin{enumerate}
\item $\displaystyle E^fM_{\s\t}X^\kappa S_d\cdot S_i\equiv
  \sum_{\substack{(\v,\rho,e)\in\delta(f,\lambda)\\|\rho|\le|\kappa|}}
              a_{\v\rho e}E^fM_{\s\v}X^\rho S_e \pmod{\Wlam}$,
\item $\displaystyle E^fM_{\s\t}X^\kappa S_d\cdot E_i\equiv
  \sum_{\substack{(\v,\rho,e)\in\delta(f,\lambda)\\|\rho|\le|\kappa|}}
              b_{\v\rho e}E^fM_{\s\v}X^\rho S_e \pmod{\Wlam}$,
\item $\displaystyle E^fM_{\s\t}X^\kappa S_d\cdot X_j\equiv
  \sum_{\substack{(\v,\rho,e)\in\delta(f,\lambda)\\|\rho|\le|\kappa|+1}}
              c_{\v\rho e}E^fM_{\s\v}X^\rho S_e \pmod{\Wlam}$.
\end{enumerate}
\end{Prop}

\begin{proof}The case $|\kappa|=0$ is precisely Lemma~\ref{|k|=0}. We now
assume that $|\kappa|>0$ and argue by induction on $|\kappa|$.

(a) Write $S_dS_i=S_aS_bS_e$, where $a\in\Sym_{n-2f}$, $b\in \Bcal_f$ and
$e\in \Dcal_f$. As $E^fM_{\s\t}X^\kappa=E^fX^\kappa M_{\s\t}$ we may
assume that $a=1$ by repeating the argument from the proof part~(a) of
Lemma~\ref{|k|=0}. By the right handed version of Lemma~\ref{SX^a},
$X^\kappa S_b=S_bX^{\kappa b^{-1}}+X$, where $X$ is a linear combination
of monomials of the form $x_1\dots x_k$ with
$x_j\in\set{S_l,E_l,X_m|1\le l<n\text{ and }1\le m\le n}$
and $k<|\kappa|$.  For each summand $x_1\dots x_k$ of $X$ we have
$k<|\kappa|$ so by induction we can write
$(E^fM_{\s\t}+\Wlam)x_1\dots x_l$ in the required form, for
$l=1,\dots,k$; consequently, by induction, we can write
$(E^fM_{\s\t}+\Wlam)x_1\dots x_kS_e$ in the required form.  Hence, we
are reduced to showing that $E^fM_{\s\t}S_bX^{\kappa b^{-1}}S_e+\Wlam$
can be written in the required form. Now,
$E^fM_{\s\t}S_b=E^fS_bM_{\s\t}=E^fM_{\s\t}$ by Lemma~\ref{M_st
properties}(c) and Lemma~\ref{E^f}. Therefore, using Lemma~\ref{E^f}
once again,
$$E^fM_{\s\t}S_bX^{\kappa b^{-1}}S_e
    =E^fM_{\s\t}X^{\kappa b^{-1}}S_e
    =M_{\s\t}E^fX^{\kappa b^{-1}}S_e
    =\pm M_{\s\t}E^fX^{\kappa'}S_e,$$
where $\kappa'\in\Nrf$ because $b\in B^f$ and $E_jX_{j+1}=-E_jX_j$ by
the skein relations. Hence, $E^fM_{\s\t}S_bX^{\kappa b^{-1}}S_e=\pm
E^fM_{\s\t}X^{\kappa'}S_e$ and the inductive step of the Proposition
is proved when $h=S_i$.

\medskip\noindent(b) As in the proof of part~(b) of Lemma~\ref{|k|=0}, we have
$E^fM_{\s\t}X^\kappa S_d E_i=E^fM_{\s\t}X^\kappa E_{j,k}S_d$,
where $j=(i)d^{-1}$ and $k=(i+1)d^{-1}$. Further, as $E_{j,k}=E_{k,j}$
we may assume that $j<k$ and, by part~(a), we may assume that $d=1$. So
we need to show that $E^fM_{\s\t}X^\kappa E_{j,k}+\Wlam$ has the
required form. There are two cases to consider.

\Case{b1. $k=j+1$}
We must show that $E^fM_{\s\t}X^\kappa E_j$ can be written in the required form.

First suppose that  $j<n-2f$. Then
we may repeat the argument from the proof of
part~(b) of Lemma~\ref{|k|=0} to see that $M_{\s\t}E_j\in\Ef$, so that
$$E^fM_{\s\t}X^\kappa E_j=E^fX^\kappa M_{\s\t} E_j
    \in E^fX^\kappa\Ef.$$
Hence, $E^fM_{\s\t}X^\kappa S_dE_i\in\W_{r,n}^{f+1}\subseteq\Wlam$ by
Lemma~\ref{if}, and the Proposition is true when $j<n-2f$.

Next, suppose that $j\ge n-2f$. If $\kappa_j+\kappa_{j+1}=0$
then $X^\kappa E_j=E_jX^\kappa$ so the result follows by induction.
Suppose then that $\kappa_j+\kappa_{j+1}>0$.

If $j\equiv n-1\pmod 2$ then $E_j$ is a factor of $E_f$ and $\kappa_j>0$. By
Lemma~\ref{small} we have that $E_jX_j^{\kappa_j}E_j=E_j\omega_j^{(\kappa_j)}$,
where $\omega_j^{(\kappa_j)}$ is a central element of $\W_{r,j-1}$ with
$\deg\omega_j^{(\kappa_j)}<\kappa_j$. Write
$E^f=\dot E^f E_j$ and $X^\kappa=\dot X^\kappa X_j^{\kappa_j}$. Then
$$E^fM_{\s\t}X^\kappa E_j
        =\dot E^fM_{\s\t}\dot X^\kappa E_jX_j^{\kappa_j}E_j
        =\dot E^fM_{\s\t}\dot X^\kappa E_j\omega_j^{(\kappa_j)}
        =E^fM_{\s\t}\dot X^\kappa E_j\omega_j^{(\kappa_j)}.$$
As $\deg\dot X^\kappa=|\kappa|-\kappa_j$ and
$\deg\omega_j^{(\kappa_j)}<\kappa_j$, the result now follows by
writing $\omega_j^{(\kappa_j)}$ as a linear combination of terms of
the form $x_1\dots x_l$ and applying induction to each of the products
$E^fM_{\s\t}\dot X^\kappa E_j x_1\dots x_m$, for $1\le m\le l$
(compare the proof of part~(a)).

If $j\equiv n\pmod 2$ then $E_{j+1}$ is a factor of $E^f$ and $\kappa_{j+1}>0$.
Write $E^f=\dot E^f E_{j+1}$ and $X^\kappa=\dot X^\kappa X_{j+1}^{\kappa_{j+1}}$.
Then
\begin{align*}
E^fM_{\s\t}X^\kappa E_j
  &=\dot E^fM_{\s\t}\dot X^\kappa E_{j+1}X_{j+1}^{\kappa_{j+1}}E_j
   =\pm\dot E^fM_{\s\t}\dot X^\kappa E_{j+1}X_j^{\kappa_{j+1}}E_j\\
  &=\pm\dot E^fM_{\s\t}\dot X^\kappa X_j^{\kappa_{j+1}}E_{j+1}E_j
   =\pm\dot E^fM_{\s\t}\dot X^\kappa X_j^{\kappa_{j+1}}E_{j+1}S_jS_{j+1}\\
  &=\pm\dot E^fM_{\s\t}\dot X^\kappa E_{j+1} X_j^{\kappa_{j+1}} S_jS_{j+1}.
   = E^fM_{\s\t} X^\kappa S_jS_{j+1}.
\end{align*}
Hence, the result follows by part~(a).

\Case{b2. $k>j+1$}
Since $|\kappa|>0$ we can fix $l$ with $\kappa_l\ne0$ (so
$l\equiv n-1\pmod 2$). Write $E^f=\dot E^f E_l$ and
$X^\kappa=\dot X^\kappa X_l^{\kappa_l}$. Set $l'=l$ if $l\notin\{j,k\}$ and $l'=l+1$
if $l\in\{j,k\}$, and put $l''=l'$ if $l'\ne j+1$ and $l''=k$ if $l'=j+1$. Note
that $l'\notin\{j,k\}$ and $l''\notin\{j,j+1\}$ since $k>j+1$. We have
\begin{align*}
E^fM_{\s\t}X^\kappa E_{j,k}
  &=\pm\dot E^fM_{\s\t}\dot X^\kappa E_lX_{l'}^{\kappa_l}
             S_{(j+1,k)}E_{j}S_{(j+1,k)}\\
  &=\pm\dot E^fM_{\s\t}\dot X^\kappa E_l(S_{(j+1,k)}X_{l''}^{\kappa_l}+X)
             E_{j}S_{(j+1,k)}\\
  &=\pm E^fM_{\s\t}\dot X^\kappa(S_{(j+1,k)}X_{l''}^{\kappa_l}+X)
             E_{j}S_{(j+1,k)},
\end{align*}
where $\deg X<\kappa_l$. Hence, by induction and part~(a) it suffices to show that
$E^fM_{\s\t}\dot X^\kappa S_{(j+1,k)}X_{l''}^{\kappa_l}E_{j}$ can be written
in the required form. As $l''\notin\{j,j+1\}$
$$
 E^fM_{\s\t}\dot X^\kappa S_{(j+1,k)}X_{l''}^{\kappa_l}E_{j}
   =E^fM_{\s\t}\dot X^\kappa S_{(j+1,k)}E_{j}X_{l''}^{\kappa_l}.
$$
Therefore, $E^fM_{\s\t}X^\kappa E_{j,k}$ can be written in the required form
by induction.

\medskip\noindent(c) As in the proof of Lemma~\ref{|k|=0}, we may assume that $r>1$
and, by the skein relations, $S_dX_j=X_{(j)d}S_d+B$, for some
$B\in\B_n(\omega_0)$. Hence, by parts~(a) and~(b) it suffices to show that
$E^fM_{\s\t}X^\kappa\cdot X_i$ can be written in the required form.
If $i\le n-2f$ then
$$E^fM_{\s\t}X^\kappa X_i+\Wlam=E^fX^\kappa\sigma(m_{\s\t})X_i+\Wlam
                      =E^fX^\kappa\sigma(m_{\s\t}Y_i)+\Wlam,$$
so the result follows because $\{m_{\s\t}\}$ is a cellular basis of
$\H_{r,n-2f}(\bu)$. If $i>n-2f$ and $\kappa_i<r-1$ then
$E^fM_{\s\t}X^\kappa X_i$ is of the desired form. If $\kappa_i=r-1$ then
$X_i^{\kappa_i}X_i=X_i^r$ can be written as a linear combination of
$r$--regular monomials of degree less than or equal to $\kappa_i$ by
the proof of Theorem~\ref{W-basis}. Hence, using parts~(a) and~(b) and
induction for each of these $r$--regular monomials,
$E^fM_{\s\t}X^\kappa X_i+\Wlam$ can be written in the required form.

This completes the proof of the Proposition.
\end{proof}

Recall from (\ref{* involution}) that $\Waff$ has a unique $R$--linear
anti--automorphism $*\map\Waff\Waff$ which fixes all of the generators
of $\Waff$. This involution induces an anti--isomorphism of
$\W_{r,n}$, which we also call~$*$. Thus, $S_i^*=S_i$, $E_i^*=E_i$,
$X_j^*=X_j$ and $(ab)^*=b^*a^*$, for $1\le i<n$, $1\le j\le n$ and all
$a,b\in\W_{r,n}$. Observe that $S_w^*=S_{w^{-1}}$, for $w\in\Sym_n$,
and that $M_{\s\t}^*=M_{\t\s}$.

\begin{Prop}\label{Wlam spanning} Suppose $0\le f\le \floor{n2}$ and
$\lambda\in\Lambda_r^+(n-2f)$. Then $\Wlambda/\Wlam$ is
spanned by the elements
$$\set{S_e^*X^\rho E^fM_{\s\t}X^\kappa S_d+\Wlam|
            (\t,\kappa,d),(\s,\rho,e)\in\delta(f,\lambda)}.
$$
\end{Prop}

\begin{proof}Let $W$ be the $R$--submodule of
$\Wlambda/\Wlam$ spanned by the elements in the statement
of the Proposition. As the generators $\{E^fM_{\s\t}+\Wlam\}$ of
$\Wlambda/\Wlam$ are contained in~$W$, and
$W\subseteq\Wlambda/\Wlam$, it suffices to show that $W$
is a $\W_{r,n}$--bimodule.

First, by Proposition~\ref{|k|>=0}, $W$ is closed under right multiplication
by elements of $\W_{r,n}$. To see that $W$ is also closed under left
multiplication by elements of $\W_{r,n}$ note that $(\Wlam)^*=\Wlam$ as
the set of generators for $\Wlam$ is invariant under~$*$ because
$(E^fM_{\s\t})^*=M_{\t\s}(E^f)^*=M_{\t\s}E^f=E^fM_{\t\s}$. Therefore,
if $a\in\W_{r,n}$ then
$$a(S_e^*X^\rho E^fM_{\s\t}X^\kappa S_d+\Wlam)
    =\big((S_d^*X^\kappa E^fM_{\t\s}X^\rho S_e+\Wlam)a^*\big)^*\in W,$$
by Proposition~\ref{|k|>=0}. Hence, $W$ is closed under left multiplication
by elements of~$\W_{r,n}$.
\end{proof}

Let $\Lambda_r^+=\set{(f,\lambda)|0\le f\le\floor{n2}\text{ and }
                     \lambda\in\Lambda_r^+(n-2f)}$.
If $(f,\lambda)\in\Lambda_r^+$ and
$(\s,\rho,e),(\t,\kappa,d)\in\delta(f,\lambda)$ then we define
$$C^{(f,\lambda)}_{(\s,\rho,e)(\t,\kappa,d)}
              =S_e^*X^\rho E^fM_{\s\t}X^\kappa S_d.$$
\label{W cell basis}

We can now prove Theorem~B from the introduction.

\begin{Theorem}\label{W cellular}
Let $R$ be a commutative ring in which $2$ is invertible and let
$\bu\in R^r$. Suppose that $\Omega$ is $\bu$--admissible. Then
$$\mathscr C=\set{C^{(f,\lambda)}_{(\s,\rho,e)(\t,\kappa,d)}|
            (\s,\rho,e),(\t,\kappa,d)\in\delta(f,\lambda),
              \text{ where }(f,\lambda)\in\Lambda_r^+}$$
is a cellular basis of $\W_{r,n}(\bu)$.
\end{Theorem}

\begin{proof} Applying the definitions it is easy to check that
$(C^{(f,\lambda)}_{(\s,\rho,e)(\t,\kappa,d)})^*
          =C^{(f,\lambda)}_{(\t,\kappa,d)(\s,\rho,e)}$. Furthermore,
by Proposition~\ref{|k|>=0}, for each $h\in\W_{r,n}$ there exist scalars
$r_{(\t',\kappa',d')}(h)\in R$, which do not depend on
$(\s,\rho,e)$, such that
$$C^{(f,\lambda)}_{(\s,\rho,e)(\t,\kappa,d)}\cdot h
     =\sum_{(\t',\kappa',d')\in\delta(f,\lambda)} r_{(\t',\kappa',d')}(h)
       C^{(f,\lambda)}_{(\s,\rho,e)(\t',\kappa',d')} \pmod\Wlam.$$

To show that $\mathscr C$ is a cellular basis of
$\W_{r,n}$ it remains to check that $\mathscr C$ is a basis of~$\W_{r,n}$.  
Now, $\W_{r,n}=\W_{r,n}^{0}\supset \W_{r,n}^{1}\supset\cdots\supset
         \W_{r,n}^{\floor{n2}}$
is a filtration of $\W_{r,n}$ by two--sided ideals, and the two--sided
ideals $\W_{r,n}^{\gedom\lambda}$, where $\lambda\in\Lambda_r^+(n-2f)$,
induce a filtration of~$\W_{r,n}^f/\W_{r,n}^{f+1}$.
Therefore,~$\mathscr C$ spans $\W_{r,n}$ by
Proposition~\ref{Wlam spanning}. To complete the proof observe that
$\#\delta(f,\lambda)=\#\UPD_n(\lambda)$, by Lemma~\ref{updown},
and $\#\mathscr C=r^n(2n-1)!!$, by Corollary~\ref{(2n-1)!!}.
As $\W_{r,n}$ is a free $R$--module of rank $r^n(2n-1)!!$ by
Theorem~\ref{W-basis},  this implies that $\mathscr C$ is an
$R$--basis of $\W_{r,n}$.  Hence, $\mathscr C$ is a cellular basis of
$\W_{r,n}$ as required.
\end{proof}

The reader may check that the proof of Theorem~\ref{W cellular} does
not rely on the explicit definition of the elements
$M_{\s\t}\in\W_{t,n}(\bu)$. The important property of these elements,
as far as the proof of the Theorem is concerned, is that they are
related to a cellular basis of $\H_{r,n}(\bu)$ by the formula of
Lemma~\ref{M_st properties}(d). Consequently, for each cellular basis
of $\H_{r,n}(\bu)$ the argument of Theorem~\ref{W cellular} produces a
corresponding cellular basis of~$\W_{r,n}(\bu)$.

We now show that we can, in principle, construct all of the finite
dimensional irreducible representations of the affine Wenzl algebras
over an algebraically closed field. First
recall that $\Omega$ is rational if for $i\gg0$ it satisfies a
recurrence relation of the form
$\omega_{i+k}+a_1\omega_{i+k-1}+\dots+a_k\omega_i=0$, for some $k>0$
and some $a_1,\dots,a_k\in R$.

\begin{Lemma}\label{rational E} 
Suppose that $\Omega$ is admissible and that $R$ is an algebraically
closed field. Then $\Omega$ is rational if and only if there is a
finite dimensional $\Waff_n(\Omega)$--module upon which $E_1$ is
non--zero.
\end{Lemma}

\begin{proof} First suppose that $\Omega$ is rational. As in the proof 
of Proposition~\ref{rational irreds}, $\Omega$ is rational if and
only if 
$$
\tilde W_1(y)+y-\frac1{2}=(y+\frac1{2})\prod_{i=1}^s \frac{y+c_i}{y-c_i},
$$
for some $c_i\in R$ and some $s\ge0$. Hence, if $\Omega$ is rational
then $\Omega$ is $\bu$--admissible where 
$$\bu=\begin{cases}
   (c_1,\dots,c_s),&\text{if $s$ is odd,}\\
   (c_1,\dots,c_s,0),&\text{if $s$ is even.}
\end{cases}$$
Hence, $\W_{r,n}(\bu)$ is a finite dimensional
$\Waff_n(\Omega)$--module upon which the action of~$E_1$ is non--zero.

Conversely, suppose that there is a finite dimensional
$\Waff_n(\bu)$--module $M$ upon which $E_1$ is non--zero. Let
$c(t)=\det(tI-X_1)$ be the characteristic polynomial for $X_1$ acting
on~$M$, where $t$ is an indeterminate and $I$ is the identity matrix
on~$M$. Write $c(t)=\sum_{j=0}^ka_jt^{k-j}$, where $a_0=1$. Then
$\sum_{j=0}^ka_jX_1^{k-j}=0$ on~$M$ by the Cayley--Hamilton theorem.
Hence, $\sum_{j=0}^ka_jE_1X_1^{i+k-j}E_1=\sum_{j=0}^ka_j\omega_{i+k-j}E_1$ 
is zero on~$M$, for any $i\ge0$. Therefore, 
$\omega_{i+k}+a_1\omega_{i+k-1}+\dots+a_k\omega_i=0$, for $i\ge0$,
since~$E_1$ is non--zero on~$M$. Thus, $\Omega$ is rational as required.
\end{proof}

\begin{Theorem}\label{finite irreds}
Suppose that $R$ is an algebraically closed field. Then we can
construct all of the finite dimensional irreducible
$\Waff_n(\Omega)$--modules.
\end{Theorem}

\begin{proof} First suppose that $\Waff_n(\Omega)$ has a finite
    dimensional irreducible module upon which $E_1$ is non--zero. Then
    $\Omega$ is admissible. Then $\Omega$ is rational  by
    Lemma~\ref{rational E}. Hence, by Proposition~\ref{rational
    irreds} every finite dimensional irreducible
    $\Waff_n(\Omega)$--module can be considered as a finite
    dimensional $\W_{r,n}(\bu)$--module for some $\bu\in R^r$ such
    that $\Omega$ is $\bu$--admissible. By Theorem~\ref{W cellular}
    $\W_{r,n}(\bu)$ is a cellular algebra, so every irreducible
    $\W_{r,n}(\bu)$--module arises (in a unique way) as the head of
    some cell module. Hence, we have a construction of every finite
    dimensional irreducible $\Waff_n(\Omega)$--module when $\Omega$ is
    rational.

Finally, suppose that $E_1$ acts as zero on every finite dimensional
irreducible $\Waff_n(\Omega)$--module. Then every finite dimensional
irreducible module $M$ can be considered as an irreducible module for
the degenerate affine Hecke algebra of type~$A$. Therefore,~$M$ can be
considered as an irreducible module for some degenerate Hecke algebra
$\H_{r,n}(\bu)$ (cf.~the proof of Proposition~\ref{rational irreds}).
Now $\H_{r,n}(\bu)$ is a cellular algebra by Theorem~\ref{H cellular},
so we can again construct all finite dimensional
$\Waff_n(\Omega)$--modules.
\end{proof}

Note that any given irreducible $\Waff_n(\Omega)$--module can be
considered as an irreducible module for an infinite number of
cyclotomic Nazarov--Wenzl algebras.  Consequently, the classification
of the irreducible $\W_{r,n}(\bu)$--modules when $\omega_0\ne0$
(Theorem~C), does not give a classification of the finite dimensional
irreducible~$\Waff_n(\Omega)$--modules when $\Omega$ is admissible and
$\omega_0\ne0$.

%

\section{Classification of the irreducible $\W_{r,n}$-modules}

In this section we classify the irreducible $\W_{r,n}(\bu)$--modules,
for fields in which~$2$ is invertible, in terms of the irreducible
$\H_{r,n}(\bu)$--modules. As the involution $*$ induces a functorial
bijection between left $\W_{r,n}$--modules and right
$\W_{r,n}$--modules, we continue to work with right
$\W_{r,n}$--modules as in the previous section.

We begin by recalling a useful result of Wenzl's.

\begin{Lemma}[\protect{Wenzl~\cite[Propositions 2.1(a) and 2.2(a)]{Wenzl:ssbrauer}}]
\label{dec}\hspace*{20mm}
  \begin{enumerate}
      \item Any monomial $B_{\gamma}\in
          \B_n(\omega_0)$ is either in $\B_{n-1}(\omega)$ or it can
          be written in the form $a_1 \alpha a_2$, where $a_i\in
          \B_{n-1}(\omega_0)$ and $\alpha\in \{E_{n-1}, S_{n-1}\}$.
      \item $E_{n-1} \B_{n}(\omega_0) E_{n-1} =\B_{n-2}(\omega_0)
          E_{n-1}$.
  \end{enumerate}
\end{Lemma}

\begin{Lemma}\label{n1} Suppose that $n\ge 2$. Then
$S_{n-1}\B_{n-1}(\omega_0)E_{n-1}=\B_{n-1}(\omega_0)E_{n-1}$.
\end{Lemma}
\begin{proof}
If $a\in \B_{n-2}(\omega_0)$ then  $S_{n-1} a
E_{n-1}=aS_{n-1}E_{n-1}=aE_{n-1}$. Suppose $a\not\in
\B_{n-2}(\omega_0)$. By Lemma~\ref{dec}, we can write $a=a_1
\alpha a_2$ with $a_i\in \B_{n-2}(\omega_0)$ and $\alpha\in \{
E_{n-2}, S_{n-2}\}$. If $\alpha=E_{n-2}$, then $S_{n-1} a E_{n-1}
= a_1S_{n-1}E_{n-2}E_{n-1}a_2=a_1S_{n-2}a_2E_{n-1}$. If
$\alpha=S_{n-2}$ then $ S_{n-1} a E_{n-1} =
a_1S_{n-1}S_{n-2}E_{n-1}a_2 =a_1E_{n-2}a_2E_{n-1}$.
 In all cases we have $S_{n-1} a E_{n-1}\in \B_{n-1}(\omega_0) E_{n-1}$.
\end{proof}

\begin{Lemma}  Suppose  that $n\ge2$. Then for each $a\in\W_{r,n}$ there exists $h\in\W_{r,n-2}$ such that $\deg h\le\deg a$ and $E_{n-1}aE_{n-1}=hE_{n-1}$.
In particular, $E_{n-1} \W_{r,n} E_{n-1}=\W_{r, n-2} E_{n-1}$.
\end{Lemma}

\begin{proof} We argue by induction on $\deg a$. It is enough to consider the case
where $a=X^{\alpha}B_\gamma X^{\beta}$ is an $r$-regular
monomial in $\W_{r,n}$. Write
$X^\alpha=\dot X^\alpha X_{n-1}^{\alpha_{n-1}}X_n^{\alpha_n}$ and
$X^\beta=\dot X^\beta X_{n-1}^{\beta_{n-1}}X_n^{\beta_n}$ and define
$k=\alpha_{n-1}+\alpha_n+\beta_{n-1}+\beta_n$.  If $k=0$ then the result follows
from Lemma~\ref{dec}(b), so we may assume that $k>0$. We split the proof into
two cases.

\Case{1. $B_\gamma\in \B_{n-1}(\omega_0)$}
First suppose that $B_\gamma\in \B_{n-2}(\omega_0)$. Then we have
\begin{align*}
E_{n-1} X^\alpha B_\gamma X^\beta E_{n-1}
 &=\dot X^\alpha B_\gamma E_{n-1}X_{n-1}^{\alpha_{n-1}+\beta_{n-1}}
               X_n^{\alpha_n+\beta_n}E_{n-1}\dot X^\beta\\
 &=(-1)^{\alpha_n+\beta_n}\dot X^\alpha B_\gamma E_{n-1}X_{n-1}^kE_{n-1}\dot X^\beta.
\end{align*}
However, $E_{n-1}X_{n-1}^kE_{n-1}=\omega_{n-1}^{(k)}E_{n-1}$ by
Lemma~\ref{tilde W}, where $\omega_{n-1}^{(k)}$ is a central element
in $\W_{r, n-2}$. If $k<r$ then $\deg\omega_{n-1}^{(k)}<k$ by
Lemma~\ref{small}, so the result follows by induction. Suppose then
that $k\ge r$ then $X_{n-1}^k$ can be written as a linear combination
of $r$--regular monomials of degree strictly less than~$k$, so the
result again follows by induction if $B_\gamma\in \B_{n-2}(\omega_0)$.

Next, suppose that $B_\gamma\notin\B_{n-2}(\omega_0)$. Then
$B_\gamma=B_{\gamma'} z B_{\gamma''}$, where
$B_{\gamma'},B_{\gamma''}\in \B_{n-2}(\omega_0)$ and
$z\in \{E_{n-2}, S_{n-2}\}$. So
$E_{n-1}X^\alpha B_\gamma X^\beta E_{n-1}
   =B_{\gamma'}\dot X^\alpha E_{n-1}X_{n-1}^{\alpha_{n-1}}X_n^{\alpha_n}z
           X_{n-1}^{\beta_{n-1}}X_n^{\beta_n}E_{n-1}\dot X^\beta B_{\gamma''}$.

If $z=E_{n-2}$ then
\begin{align*}
E_{n-1}X_{n-1}^{\alpha_{n-1}}X_n^{\alpha_n}E_{n-2}
           X_{n-1}^{\beta_{n-1}}X_n^{\beta_n}E_{n-1}
    &=\pm E_{n-1}X_{n-1}^{\alpha_{n-1}+\alpha_n}E_{n-2}
           X_{n-1}^{\beta_{n-1}+\beta_n}E_{n-1}\\
    &=\pm E_{n-1}X_{n-2}^{\alpha_{n-1}+\alpha_n}E_{n-2}
           X_{n-2}^{\beta_{n-1}+\beta_n}E_{n-1}\\
    &=\pm X_{n-2}^{\alpha_{n-1}+\alpha_n}E_{n-1}
              E_{n-2}E_{n-1}X_{n-2}^{\beta_{n-1}+\beta_n}\\
    &=\pm X_{n-2}^{\alpha_{n-1}+\alpha_n}E_{n-1}
                  X_{n-2}^{\beta_{n-1}+\beta_n}\\
    &=\pm X_{n-2}^kE_{n-1}.
\end{align*}
This completes the proof when $z=E_{n-2}$.

Now suppose that $z=S_{n-2}$. Using the relations
(\ref{Waff relations}),
\begin{align*}
E_{n-1}X_{n-1}^{\alpha_{n-1}}X_n^{\alpha_n}S_{n-2}
           X_{n-1}^{\beta_{n-1}}X_n^{\beta_n}E_{n-1}
  &=\pm E_{n-1}X_n^{\alpha_{n-1}+\alpha_n}S_{n-2}X_n^{\beta_{n-1}+\beta_n}E_{n-1}\\
  &=\pm E_{n-1}S_{n-2}X_n^kE_{n-1}\\
  &=\pm E_{n-1}E_{n-2}S_{n-1}X_n^kE_{n-1}.
\end{align*}
If $k\ge r$ then we can write $X_n^k$ as a linear combination of $r$--regular
monomials each with degree strictly less than $k$. So by induction we may assume that
$k<r$. Then,  by Lemma~\ref{SX^a}, $S_{n-1}X_n^k=X_{n-1}^kS_{n-1}+X$, where
$X\in\W_{r,n}$ is a linear combination of terms each of which has
total degree in $X_n$ and $X_{n-1}$ strictly less than $k$. Hence, by
induction, $E_{n-1}E_{n-2}XE_{n-1}=h'E_{n-1}$, where $h'\in\W_{r,n-2}$ and
$\deg h'<k$. Further,
\begin{align*}
E_{n-1}E_{n-2}X_{n-1}^kS_{n-1}E_{n-1}
  &=E_{n-1}E_{n-2}X_{n-1}^kE_{n-1}
   =(-1)^kE_{n-1}E_{n-2}X_{n-2}^kE_{n-1}\\
  &=(-1)^kE_{n-1}E_{n-2}E_{n-1}X_{n-2}^k
   =(-1)^kX_{n-2}^kE_{n-1}.
\end{align*}
Consequently, $E_{n-1}X^\alpha B_\gamma X^\beta E_{n-1}=hE_{n-1}$, where
$h\in\W_{r,n-2}$ and $\deg h\le \deg a$.

\Case{2. $B_\gamma\not\in \B_{n-1}(\omega_0)$}
Once again by Lemma~\ref{dec} we can write
$B_\gamma=B_{\gamma'}z B_{\gamma''}$, where
$B_{\gamma'}, B_{\gamma''}\in\B_{n-1}(\omega_0)$ and
$z\in \{S_{n-1}, E_{n-1}\}$.

If $z=E_{n-1}$ then the result follows using Case~1 twice, so suppose that
$z=S_{n-1}$. Then
$$E_{n-1} X^\alpha B_\gamma  X^\beta E_{n-1}
    =\dot X^\alpha E_{n-1}X_{n-1}^{\alpha_{n-1}}X_n^{\alpha_n}B_{\gamma'}
    S_{n-1}B_{\gamma''}X_{n-1}^{\beta_{n-1}}X_n^{\beta_n}E_{n-1}\dot X^\beta.$$
If $\beta_{n-1}+\beta_n=0$ then $S_{n-1}B_{\gamma''}E_{n-1}=h E_{n-1}$, for some
$h\in\B_{n-1}(\omega_0)$ by Lemma~\ref{n1}, so the result follows from Case~1.
Hence, we may assume that $\beta_{n-1}+\beta_n>0$.
Similarly, we may assume that $\alpha_{n-1}+\alpha_n>0$.

Next, suppose that $B_{\gamma''}\in \B_{n-2}(\omega_0)$. Then
$$E_{n-1} X^\alpha B_\gamma  X^\beta E_{n-1}
    =\pm\dot X^\alpha E_{n-1} X_{n-1}^{\alpha_{n-1}+\alpha_n} B_{\gamma'}B_{\gamma''}
           S_{n-1}X_n^{\beta_{n-1}+\beta_n} E_{n-1}\dot X^\beta.$$
Once again, by induction we may assume that $\beta_{n-1}+\beta_n<r$. Then,
by Lemma~\ref{SX^a}
$S_{n-1}X_n^{\beta_{n-1}+\beta_n}=X_{n-1}^{\beta_{n-1}+\beta_n}S_{n-1}+X$,
where $\deg X<\beta_{n-1}+\beta_n$. As $S_{n-1}E_{n-1}=E_{n-1}$ it is enough
to consider
$E_{n-1} X_{n-1}^{\alpha_{n-1}+\alpha_n} B_{\gamma'}B_{\gamma''}
           X_{n-1}^{\beta_{n-1}+\beta_n}E_{n-1}$.
As $B_{\gamma'}B_{\gamma''}\in\B_{n-1}(\omega_0)$,
this can be written in the required form by Case~1.

Finally, suppose that $B_{\gamma''}\notin\B_{n-2}(\omega_0)$. Then either
$B_{\gamma''}=B_{\gamma''_1}E_{n-2}B_{\gamma''_2}$, or
$B_{\gamma''}=B_{\gamma''_1}S_{n-2}B_{\gamma''_2}$, where $B_{\gamma''_1},
B_{\gamma''_2}\in\B_{n-2}(\omega_0)$.
If $B_{\gamma''}=B_{\gamma''_1}E_{n-2}B_{\gamma''_2}$ then
$$\begin{array}{l}
E_{n-1}X_{n-1}^{\alpha_{n-1}}X_n^{\alpha_n}B_{\gamma'}S_{n-1}
       B_{\gamma''_1}E_{n-2}B_{\gamma''_2}X_{n-1}^{\beta_{n-1}}X_n^{\beta_n}E_{n-1}\\
 \qquad=\pm E_{n-1}X_n^{\alpha_{n-1}+\alpha_n}B_{\gamma'}
       B_{\gamma''_1}S_{n-1}E_{n-2}B_{\gamma''_2}X_{n-1}^{\beta_{n-1}+\beta_n}E_{n-1}\\
 \qquad=\pm E_{n-1}X_n^{\alpha_{n-1}+\alpha_n}B_{\gamma'}
      B_{\gamma''_1}S_{n-1}E_{n-2}X_{n-1}^{\beta_{n-1}+\beta_n}E_{n-1}B_{\gamma''_2}\\
 \qquad=\pm E_{n-1}X_n^{\alpha_{n-1}+\alpha_n}B_{\gamma'}
       B_{\gamma''_1}S_{n-1}E_{n-2}X_{n-2}^{\beta_{n-1}+\beta_n}E_{n-1}B_{\gamma''_2}\\
 \qquad=\pm E_{n-1}X_n^{\alpha_{n-1}+\alpha_n}B_{\gamma'}
       B_{\gamma''_1}S_{n-1}E_{n-2}E_{n-1}X_{n-2}^{\beta_{n-1}+\beta_n}B_{\gamma''_2}\\
 \qquad=\pm E_{n-1}X_n^{\alpha_{n-1}+\alpha_n}B_{\gamma'}
       B_{\gamma''_1}S_{n-2}E_{n-1}X_{n-2}^{\beta_{n-1}+\beta_n}B_{\gamma''_2}.
\end{array}$$
Now $\deg(X^\alpha B_{\gamma'}B_{\gamma''_1}S_{n-2})<\deg a$ since
$\beta_{n-1}+\beta_n>0$. Hence, the result now follows by induction.
If $B_{\gamma''}=B_{\gamma''_1}S_{n-2}B_{\gamma''_2}$ then
$$\begin{array}{l}
E_{n-1}X_{n-1}^{\alpha_{n-1}}X_n^{\alpha_n}B_{\gamma'}S_{n-1}
       B_{\gamma''_1}S_{n-2}B_{\gamma''_2}X_{n-1}^{\beta_{n-1}}X_n^{\beta_n}E_{n-1}\\
 \qquad=\pm E_{n-1}X_n^{\alpha_{n-1}+\alpha_n}B_{\gamma'}
       B_{\gamma''_1}S_{n-1}S_{n-2}X_n^{\beta_{n-1}+\beta_n}E_{n-1}B_{\gamma''_2}\\
 \qquad=\pm E_{n-1}X_n^{\alpha_{n-1}+\alpha_n}B_{\gamma'}
       B_{\gamma''_1}S_{n-1}X_n^{\beta_{n-1}+\beta_n}S_{n-2}E_{n-1}B_{\gamma''_2}.
\end{array}$$
By Lemma~\ref{SX^a} we can write
$S_{n-1}X_n^{\beta_{n-1}+\beta_n}=X_{n-1}^{\beta_{n-1}+\beta_n}S_{n-1}+X$, where
$\deg X<\beta_{n-1}+\beta_n$. Now,
$$\begin{array}{l}
E_{n-1}X_n^{\alpha_{n-1}+\alpha_n}B_{\gamma'}B_{\gamma''_1}
       X_{n-1}^{\beta_{n-1}+\beta_n}S_{n-1}S_{n-2}E_{n-1}B_{\gamma''_2}\\
  \qquad=E_{n-1}B_{\gamma'}B_{\gamma''_1} X_{n-1}^{\beta_{n-1}+\beta_n}
             X_n^{\alpha_{n-1}+\alpha_n}E_{n-2}E_{n-1}B_{\gamma''_2}\\
  \qquad=E_{n-1}B_{\gamma'}B_{\gamma''_1} X_{n-1}^{\beta_{n-1}+\beta_n}
             E_{n-2}X_n^{\alpha_{n-1}+\alpha_n}E_{n-1}B_{\gamma''_2}\\
  \qquad=E_{n-1}B_{\gamma'}B_{\gamma''_1} X_{n-1}^{\beta_{n-1}+\beta_n}
             E_{n-2}X_{n-2}^{\alpha_{n-1}+\alpha_n}E_{n-1}B_{\gamma''_2}\\
  \qquad=E_{n-1}B_{\gamma'}B_{\gamma''_1} X_{n-1}^{\beta_{n-1}+\beta_n}
             E_{n-2}E_{n-1}B_{\gamma''_2}X_{n-2}^{\alpha_{n-1}+\alpha_n}.
\end{array}$$
As $\alpha_{n-1}+\alpha_n>0$ we can write
$E_{n-1}B_{\gamma'}B_{\gamma''_1}X_{n-1}^{\beta_{n-1}+\beta_n} E_{n-2}E_{n-1}$
in the required form and so completes the proof of the case---and hence the Lemma.
\end{proof}

By iterating the Lemma we obtain the result that we really want.

\begin{Cor}\label{class2}Suppose $f>0$, $w\in \Sym_n$ and that
$\kappa,\rho\in\Nrf$. Then
$$E^fX^\rho S_w X^\kappa E^f =hE^f,$$
for some $h\in \W_{r, n-2f}$.
\end{Cor}

As we now briefly recall, by the general theory of cellular
algebras~\cite{GL,M:ULect}, every irreducible $\W_{r,n}$--module
arises in a unique way as the simple head of some cell module.
For each $(f,\lambda)\in\Lambda_r^+$ fix
$(\s,\rho,e)\in\delta(f,\lambda)$ and let
$C^{(f,\lambda)}_{(\t,\kappa,d)}
       =C^{(f,\lambda)}_{(\s,\rho,e)(\t,\kappa,d)}+\Wlam$.
By Theorem~\ref{W cellular} the cell modules of $\W_{r,n}$ are the
modules $\Delta(f,\lambda)$ which are the free $R$--modules with
basis
$\set{C^{(f,\lambda)}_{(\t,\kappa,d)}|(\t,\kappa,d)\in\delta(f,\lambda)}$.
The cell module $\Delta(f,\lambda)$ comes equipped with a natural
bilinear form $\phi_{f,\lambda}$ which is determined by the
equation
$$C^{(f,\lambda)}_{(\s,\rho,e)(\t,\kappa,d)}
           C^{(f,\lambda)}_{(\t',\kappa',d')(\s,\rho,e)}
 \equiv\phi_{f,\lambda}\big(C^{(f,\lambda)}_{(\t,\kappa,d)},
               C^{(f,\lambda)}_{(\t',\kappa',d')}\big)\cdot
        C^{(f,\lambda)}_{(\s,\rho,e)(\s,\rho,e)}\pmod{\Wlam}.$$
The form $\phi_{f,\lambda}$ is $\W_{r,n}$--invariant in the sense that
$\phi_{f,\lambda}(xa,y)=\phi_{f,\lambda}(x,ya^*)$, for
$x,y\in\Delta(f,\lambda)$ and $a\in\W_{r,n}$. Consequently,
$$\Rad\Delta(f,\lambda)
   =\set{x\in\Delta(f,\lambda)|\phi_{f,\lambda}(x,y)=0\text{ for all }
                          y\in\Delta(f,\lambda)}$$
is a $\W_{r,n}$--submodule of $\Delta(f,\lambda)$ and
$D(f,\lambda)=\Delta(f,\lambda)/\Rad\Delta(f,\lambda)$ is either zero
or absolutely irreducible.

In exactly the same way, for each multipartition
$\lambda\in\Lambda_r^+(n-2f)$ the corresponding cell module
$S(\lambda)$ for $\H_{r,n-2f}(\bu)$, the Specht module of section~6,
carries a bilinear form $\phi_\lambda$. The quotient module
$D(\lambda)=S(\lambda)/\Rad S(\lambda)$ is either zero or an
absolutely irreducible $\H_{r,n-2f}(\bu)$--module.

We can now prove Theorem~\ref{simple classification}.

\begin{Theorem}\label{W simples}
Suppose that $R$ is a field in which $2$ is invertible, that $\Omega$ is
$\bu$--admissible and that $\omega_0\neq 0$. Let
$(f,\lambda)\in\Lambda_r^+$. Then $D^{(f,\lambda)}\neq 0$ if and only
if~$D^\lambda\neq0$.
\end{Theorem}

\begin{proof}
It is enough to prove that $\phi_{f,\lambda}\ne0$ if and
only if $\phi_\lambda\neq 0$.

First, suppose that $\phi_{\lambda}\neq 0$. Recall that the Specht
module $S(\lambda)$ has basis
$\set{m_\t|\t\in\Std(\lambda)}$. Then
$\phi_\lambda(m_\t, m_\v)\ne0$,
for some $\t, \v\in \Std(\lambda)$; that is,
$m_{\s\t}m_{\v\s}\notin\H_{r,n-2f}^{\gdom\lambda}$. Let
$\0$ to the zero vector in $\Nrf$. Then
\begin{align*}
C^{(f,\lambda)}_{(\s,\rho,e)(\t,\0,1)}
  C^{(f,\lambda)}_{(\v,\0,1)(\s,\rho,e)}
   &=S_e^*X^\rho E^f M_{\s\t} E^fM_{\v\s}X^\rho S_e\\
   &=S_e^*X^\rho (E^f)^2 M_{\s\t}M_{\v\s}X^\rho S_e\\
   &\equiv\omega_0^f\phi_\lambda(m_\t,m_\v)
         S_e^*X^\rho E^f M_{\s\s}X^\rho S_e\\
   &\equiv\omega_0^f\phi_\lambda(m_\t,m_\v)
         C^{(f,\lambda)}_{(\s,\rho,e)(\s,\rho,e)}\pmod\Wlam.
\end{align*}
Hence,
$\phi_{f,\lambda}\big(
     C^{(f,\lambda)}_{(\t,\0,1)},C^{(f,\lambda)}_{(\v,\0,1)} \big)
       =\omega_0^f\phi_{\lambda}(m_\t,m_\v)\neq 0$,
so that $\phi_{f,\lambda}\ne0$.

Now suppose that $\phi_{f,\lambda}\neq 0$. Then there exist
$(\u,\alpha,u),(\v,\beta',v)\in\delta(f,\lambda)$ such that
$\phi_{f,\lambda}\big(C^{(f,\lambda)}_{\u,\alpha,u)},
                   C^{(f,\lambda)}_{(\v,\beta,v)}\big)\ne0$. That
is,
\begin{align*}
0&\ne C^{(f,\lambda)}_{(\s,\rho,e)(\u,\alpha,u)}\cdot
       C^{(f,\lambda)}_{(\v,\beta,v)(\s,\rho,e)}\\
     &=S_e^*X^\rho E^fM_{\s\u}X^\alpha S_u\cdot
            S_v^*X^\beta E^fM_{\v\s}X^\rho S_e\\
     &=S_e^*X^\rho M_{\s\u}E^fX^\alpha S_u
            S_v^*X^\beta E^fM_{\v\s}X^\rho S_e\\
     &=S_e^*X^\rho M_{\s\u}hM_{\v\s}E^fX^\rho S_e,
\end{align*}
for some $h\in\W_{r,n-2f}$ by Corollary~\ref{class2}. Now,
$M_{\s\u}\Ef M_{\v\s}E^f\subseteq\Ef E^f\subseteq\W_{r,n}^{f+1}$, by
Lemma~\ref{if}. Therefore, Proposition~\ref{degen isom} implies that
there is an $h'\in \H_{r,n-2f}(\bu)$ such that
$m_{\s\u} h' m_{\v\s}\neq 0\pmod{\H_{r,n-2f}^{\gdom\lambda}}$.
Consequently, $\phi_\lambda\ne0$. This completes the proof of the
Theorem.
\end{proof}

We remark that the irreducible representations of the Ariki-Koike
algebras are indexed by the $\bu$--\textit{Kleshchev multipartitions};
see \cite{Ariki:class,AM:simples}.  In the special case when
$u_i=d_i\cdot1_R$, for $1\le i\le r$ and where $0\le d_i<\Char R$,
Kleshchev~\cite{Klesh:book} has shown that the simple
$\H_{r,n}(\bu)$--modules are \textit{labelled} by a set of multipartitions
which gives the same Kashiwara crystal as the set of $\bu$--Kleshchev
multipartitions of~$n$. Hence, in this case, the
simple $\W_{r,n}(\bu)$--modules are labelled by the set
$\{(f,\lambda)\}$, where $0\le f\le\floor{n2}$ and $\lambda$ is a
$\bu$--Kleshchev multipartition of $n-2f$. By modifying the proof of
\cite[Theorem~1.1]{DM:Morita}, or~\cite[Theorem~1.3]{AM:simples},
one can show that under the assumptions of Theorem~\ref{W simples} the
simple $\W_{r,n}(\bu)$--modules are always labelled by the
$\bu$--Kleshchev multipartitions. (Note, however, that we are not
claiming that $D^{(f,\lambda)}\ne0$ for the multipartitions $\lambda$
which Kleshchev~\cite{Klesh:book} uses to label the irreducible
$\H_{r,n}(\bu)$--modules.)

We close by classifying the quasi--hereditary cyclotomic
Nazarov--Wenzl algebras with $\omega_0\ne0$. See~\cite{CPS:qh} for the
definition of a quasi--hereditary algebra.

\begin{Cor}
Suppose that $R$ is a field in which $2$ is invertible, that $\Omega$ is
$\bu$--admissible and that $\omega_0\neq 0$. Then $\W_{r,n}(\bu)$ is a
quasi--hereditary algebra if and only if $\Char R>n$ and $\bu$ is
generic for $\H_{r,n}(\bu)$ $($Definition~\ref{H-generic}$)$.
\end{Cor}

\begin{proof}By \cite[(3.10)]{GL}, a cellular algebra is
quasi--hereditary if and only if the bilinear form on each cell module
does not vanish.  Therefore, $\W_{r,n}$ is a quasi--hereditary algebra
if and only if $D^{(f,\lambda)}\ne0$ for all
$(f,\lambda)\in\Lambda_r^+$ and $\H_{r,n-2f}(\bu)$ is
quasi--hereditary if and only if $D^\lambda\ne0$ for all
$\lambda\in\Lambda_r^+(n-2f)$.  Hence, by Theorem~\ref{W simples},
$\W_{r,n}(\bu)$ is quasi--hereditary if and only if the algebras
$\H_{r,n-2f}(\bu)$ are all quasi--hereditary, for $0\le
f\le\floor{n2}$.  However, the degenerate Hecke algebras
are Frobenius algebras by \cite[Cor.~5.7.4]{Klesh:book}, so they are
quasi--hereditary precisely when they are semisimple---since Frobenius
algebras have infinite global dimension when they are not semisimple,
whereas quasi--hereditary algebras have finite global dimension (see
\cite[Prop.~A2.3]{Donkin:book}). Hence the result follows from
Theorem~\ref{H semisimple}.
\end{proof}


\providecommand{\bysame}{\leavevmode ---\ }
\providecommand{\og}{``}
\providecommand{\fg}{''}
\providecommand{\smfandname}{and}
\providecommand{\smfedsname}{\'eds.}
\providecommand{\smfedname}{\'ed.}
\providecommand{\smfmastersthesisname}{M\'emoire}
\providecommand{\smfphdthesisname}{Th\`ese}

\end{document}